\documentclass[12pt]{article}
\usepackage{amsmath,amssymb,amsfonts,}
\usepackage[upright]{fourier}
\newcommand{\R}{\mbox{$I\!\!R$}}

\newcommand{\N}{\mbox{$I\!\!N$}}          
\newtheorem{thm}{Theorem}[section]
\newtheorem{lem}[thm]{Lemma}
\newtheorem{Def}[thm]{Definition}
\newtheorem{prop}[thm]{Proposition}

\newtheorem{rem}[thm]{Remark}
\newtheorem{ex}[thm]{Example}
\newtheorem{coro}[thm]{Corollary}
\newtheorem{thmD}[thm]{Theorem and Definition} 
\title{A USEFUL UNDERESTIMATE FOR THE CONVERGENCE OF INTEGRAL FUNCTIONALS. }
\author{ Emmanuel GINER \\  \footnotesize\centerline{Institut de Math\'ematiques, Laboratoire MIP, Universit\'e Paul Sabatier,Toulouse, France.} \\\\
}

\date{ }   
\begin{document}           
\maketitle                 
\thispagestyle{empty}
 
\footnotetext{
Current Address: Institut de Math\'ematiques, Laboratoire MIP, 
Universit\'e Paul Sabatier, 118 route de Narbonne 31062 Toulouse cedex 04, France. e-mail: giner@math.univ-toulouse.fr
}
\begin{abstract} 
\noindent This article deals with the lower compactness property of a sequence of integrands and the use of this key notion in various domains: convergence theory, optimal control, non-smooth analysis. First about the interchange of the weak lower epi-limit and the symbol of integration for a sequence of integral functionals. These functionals are defined on a topological space $( \mathcal{X},\mathcal{T})$, where  $\mathcal{X}$ is a subset of measurable functions and the $\mathcal{T}$-sequential convergence is stronger than or equal to the convergence in the Biting sense. Given a sequence $(f_{n})_n$ of integrands, if the integrand $f$ is the  weak lower sequential epi-limit of the integrands $f_{n}$ one of the main results of this article asserts that under the Ioffe's criterion, the $\mathcal{T}$-lower sequential epi-limit of the sequence of integral functionals at the point $x$ is bounded below by the value of the integral functional associated to the Fenchel-Moreau biconjugate of $f$ at the point $x$. Then the strong-weak semicontinuity (respectively the  subdifferentiability) of integral functionals, are studied in relation with the Ioffe's criterion. This permits, with original proofs, to give new conditions for the sequential strong-weak lower semi continuity {\it at a given point}, and to obtain  necessary and sufficient criteria for the Fr\'echet and the (weak) Hadamard subdifferentiability of integral functionals on general spaces, particularly on Lebesgue spaces.
\end{abstract}

\section{Introduction} The notion of convergence plays a key role in the study of variational problems and in nonlinear analysis, see \cite{120}, \cite{121}, \cite{30}, \cite{0}, \cite{130}, \cite{7}, \cite{1}, \cite{15}, \cite{110}, \cite{114}. The convergence of a sequence of functions is often defined through the convergence of their epigraphs. Given a  topological space $(\mathcal{X}, \mathcal{T})$, one can consider on the space of subsets of  $\mathcal{X}$ the Painlev\'e-Kuratowski convergence. This permits to define epi-convergence of a sequence of functions $(f_{n})_n$ defined on  $\mathcal{X}$ with extended numerical values. In fact, the use of the epigraphs of functions and the existence of a lower limit and an upper limit of a sequence of subsets allows to define the $\mathcal{T}$-upper epi-limit, $\mathcal{T}-ls_{e} f_{n}$, (respectively the $\mathcal{T}$-lower epi-limit, $\mathcal{T}-li_{e} f_{n}$) of a sequence of extended real-valued functionals. Moreover the $\mathcal{T}$-epi-limit, $\mathcal{T}-lim_{e} f_{n}$, can be defined by the coincidence of these last two epi-limits. In this case one says that the sequence $(f_{n})_n$ $\mathcal{T}$-epi-converges to $\mathcal{T}-lim_{e} f_{n}$. When we consider two topologies on  $\mathcal{X}$, $\mathcal{T}_1$ finer or equal than $\mathcal{T}_2$, we can define also the $(\mathcal{T}_{1}, \mathcal{T}_{2})$ Mosco-convergence of a sequence of functionals with extended numerical values, by the coincidence of the associated  $\mathcal{T}_1$ upper epi-limit and the $\mathcal{T}_2$ lower epi-limit. The case considered initially by U. Mosco  is a reflexive Banach space endowed with its strong and weak topologies \cite{120}. In \cite{0}, \cite{130}, \cite{7},  \cite{15} the reader can find many examples and applications of the notion of epi-convergence.\\
For integral functionals it seems that it is of great interest not to study directly the convergence, but to make separate studies of the  cases of lower and upper epi-limits.  In this article, first it is considered sequential (weak) lower epi-limits of a sequence of such functionals.\\
 Given a complete $\sigma$-finite measure space $\Omega$, a separable Banach space $E$ (or in some cases a reflexive separable Banach space), a topological subset $( \mathcal{X}, \mathcal{T})$ of the space of classes (for almost everywhere equality) of measurable $E$-valued functions, and a sequence $(f_{n})_n$  of measurable extended real-valued non necessarily convex integrands defined on $\Omega\times E$, our first purpose is to give a `` best '' lower bound for the $\mathcal{T}$-sequential lower epi-limit of the associated sequence $(I_{f_{n}})_n$ of integral functionals at a given point $x\in\mathcal{X}$. The study of convergence of integral functionals is originally started by J- L. Joly and F. De Th\'elin \cite{22} in the convex case with $E$ finite dimensional,  $\mathcal{X}=L_{p}(\Omega, E)$ and for Mosco-convergence;  then by  A. Salvadori \cite{72} when $\Omega$ is a finite measured space, $E$ is a reflexive separable Banach space; more recently for the slice convergence by J. Couvreux \cite{5}, $\Omega$ being a probability space and  $E$ a Banach space with separable dual, with  $\mathcal{X}=L_{p}(\Omega, E)$ also in the convex case. In our approach it is important to avoid any convergence, properness, or global convexity assumptions on the sequence of integrands and of functionals.\\
The second and third section are devoted to some known or new preliminaries on the calculus of the Fenchel-Moreau conjugate of a sequential weak lower epi-limit, on the Biting Lemma and tightness. In Section $4$  a measurability property is first proved: Theorem \ref{thm 5.1}.  Then the chapter VII of C. Castaing and M. Valadier'book \cite{4} is used and a consequence of  M. Valadier results \cite{890} section 3, \cite{4} Theorem VII-7 on the calculus of the Fenchel-Moreau conjugate of an integral functional is put in light. Section 5 is devoted to a property of the upper epi limit $ls_{e} I_{f_{n}}$ of a sequence $(I_{f_{n}})_n$ of integral functionals defined on the dual of $L_{1}(\Omega,E)$ and for the topology of uniform convergence on weakly compact sets of $L_{1}(\Omega,E)$, Proposition \ref{lem301}, this result permits to give a proof of the main result of section 6 by a duality method. 
In Section 6, the property $(P)$ required on the topology $\mathcal{T}$ is trivially satisfied when  $\mathcal{X}=L_{p}(\Omega,E)$ is endowed with a topology stronger than or equal  the weak topology.  Given a sequence of integrands  $(f_{n})_n$ defined on $\Omega\times E$, let the integrand $f=seq\; \sigma- li_{e} f_{n}$ and its Fenchel-Moreau biconjugate $f^{**}$, the main result is Theorem \ref{thm6.01}; it proves that given a sequence $(x_{n})_n$ $\mathcal{T}$-converging to $x$ under mild assumptions  the following inequality is true
$$\liminf_{n}\; I_{f_{n}}(x_{n})\geq I_{f^{**}}(x)-\delta^{+}((-f_{n}(x_{n}))_{n})\;,$$ 
where $\delta^{+}(.)$ is an extension to the $\sigma$-finite case of Rosenthal's modulus of uniform equi-integrability (see \cite{67}).
Following \cite{70}, we will say that a sequence  $(f_{n})_n$ of integrands satisfies the Ioffe's criterion at $x\in\mathcal{X}$ ( with respect to $\mathcal{T}$) when the following lower compactness property holds:  for every subsequence $(f_{n_{k}})_k$ of $(f_{n})_n$ and any $\mathcal{T}$-converging sequence $(x_{k})_k$ to $x$ such that the sequence $(I_{f_{n_{k}}}(x_{k}))_k$ is bounded above, the sequence of negative parts $(f_{n_{k}}^{ -}(x_{k}))_k$ is relatively weakly compact in $L_{1}(\Omega, {\R})$. 
When the Ioffe's criterion holds, the term $\delta^{+}$ vanishes and the inequality announced in the abstract holds. The converse being often valid when $f(x)=f^{**}(x)$ (Theorem \ref{thm6.2}). In  Section 7 a first application is obtained: the classical problem of strong-weak lower semicontinuity of an integral functional is considered. In this case the integrand is defined on $\Omega\times E\times F$, where $E$ is a topological space and $F$ is a separable Banach space. The global semicontinuity problem was considered by many authors, notably C. Olech \cite{103}, \cite{104}, in the case $L_{1}(\Omega,E^{2})$, and completely solved by A. D. Ioffe \cite{70} for a large class of spaces of finite dimensional-valued measurable functions. C. Castaing and P. Clauzure \cite{118} deal, with a strengthening of Olech's techniques, the case of measurable functions infinite dimensional valued: $F$ is a Banach space with separable dual. E. J. Balder \cite{101}, \cite{102}, \cite{111}, gives a new proof of this result using the concept of seminormality for an integral functional defined on $L_{1}(\Omega,E\times F)$, where $F$ is a separable reflexive Banach space. As a significant contribution he introduces and uses the notion of Nagumo tightness. A. Bourass, B. Ferrahi and O. Kahlaoui in \cite{105} show that the Ioffe's  techniques and results can be extended to the case $F$ is a separable reflexive Banach space. Moreover, when $f$ is non negative, I. Fonseca and G. Leoni have weakened the assumptions on the topologies (\cite{78} Corollary 7.9) and made a remarkable  characterization of the associated relaxed energy functional \cite{78} Theorem 7.13, but $F$ is supposed finite dimensional. More recently C. Castaing, P.R. de Fitte, M. Valadier in their book \cite{3} with various notions of tightness in relation with the theory of Young measures, obtain a semicontinuity result in case $F$ is a separable Banach space, \cite{3} Theorem 8.1.6. 
In this article, Theorem \ref{thm 1011} is first a quantitative estimate on the lack of sequential strong-weak semicontinuity. Not only it gives other proofs of the known semicontinuity results at least when $F$ has a strongly separable dual, \cite{101}, \cite{102}, \cite{111}, \cite{118}, \cite{78} and \cite{105},  but  in case $F$ is reflexive, with weak assumptions on the topologies it permits to extend the  Ioffe's result at {\it a given point} without any global convexity assumptions on the integrand, Corollary \ref{thm 101}.  
In section $8$ the lower compactness property respect to a bornology is defined. Fr\'echet and Hadamard lower compactness properties of a sequence of integrands are considered. Using growth conditions, very concrete examples are presented on Orlicz spaces,  Propositions \ref{coro554}, \ref{propw.200} and  \ref{propw.2}, and on  Lebesgue spaces  Corollaries \ref{coro157},  \ref{coro154}, \ref{coro156},  \ref{coro155}, \ref{thm8990} and \ref{coro5050}. 
The section $9$ is first devoted to introduce the notions of Fr\'echet and weak Hadamard subdifferentiability. Related to the J. P Penot's characterization of Fr\'echet subdifferentiability of an integral functional on $L_{p}(\Omega, E)$, $1\leq p<\infty$, \cite{216} Theorems 12 and 22, two complete characterizations for the subdifferentiability of an integral functional respect to the Fr\'echet and weak Hadamard bornologies, are given in a general setting: Theorem \ref{prop 020} and Theorem \ref{prop 021}; they permit to reach criterions in the case of Orlicz spaces Corollary \ref{coro613}, and of Lebesgue spaces when the measure is atomless: Corollaries \ref{coro713} and \ref{coro523} which are equivalent to the characterizations given in \cite{216} Theorems 12 and 22, Corollary \ref{coro73} treats the case $p=\infty$. In section 10 some additional properties of the Fr\'echet subdifferentiability are reached: Theorem \ref{thm4.8}, Corollaries \ref{prop450} and \ref{prop451}, they are in relation with the Fr\'echet lower compactness property of the differential quotients associated to the integrand (for practical examples see the sections 8 and 9). 
The last section is a short study of the weak  Hadamard subdifferentiability of an integral functional defined on Lebesgue spaces when the Banach space $E$ is supposed to be reflexive separable.  After a reduction with the results of the previous sections it is proved that it suffices to  treat only the case $p=1$. Recall that in this case when the measure is atomless, the Fr\'echet subdifferential of an integral functional coincide with the Moreau-Rockafellar subdifferential (see \cite{216} and \cite{777}). Corollary \ref{coro822} gives a complete characterization of weak Hadamard subdifferentiability when the measure is atomless. The proof of Theorem \ref{thm6.17} uses Theorem \ref{thm6.01}, and it is unreachable with the A.D Ioffe's type convexity assumptions made in the yet known results on semicontinuity.  
  
\section{Preliminaries} 
We adopt the following notation: ${\R}$ is the set of real numbers  and $\overline{{\R}}={\R}\cup\{\pm\infty\}$.  
Given $(\mathcal{X}, \mathcal{T})$ a topological vector space, the set of all open neighbourhoods of $x$ in $X$ will be denoted by $\mathcal{N}(x)$. For a subset $X$ of  $\mathcal{X}$, the indicator function $\iota_{X}$ is defined by $\iota_{X}(x)=0$ if $x\in X$, $+\infty$ if not.
For an extended real-valued  function $f$ defined on   $\mathcal{X}$ we consider its effective domain, $domf=\{x\in\mathcal{X}:f(x)<+\infty\}$, 
its epigraph, $epif=\{(x,r)\in\mathcal{X}\times{\R}:f(x)\leq r\}$ and its sublevel set of height $r\in {\R}$ $f^{\leq r}=\{x\in\mathcal{X}: f(x)\leq r\}$. 
 The function $f$ is said to be $\mathcal{T}$-inf-compact (respectively sequentially-$\mathcal{T}$ inf-compact) when every sublevel set is $\mathcal{T}$-compact (respectively sequentially $\mathcal{T}$-compact). When  $\mathcal{X}$ is a locally convex topological space with $\mathcal{X}^{*}$ as topological dual we consider the duality pairing between  $\mathcal{X}^{*}$ and  $\mathcal{X}$ defined by $\langle x^{*}, x\rangle=x^{*}(x)$. The function $f$ is said inf-compact for every slope (respectively sequentially-$\mathcal{T}$ inf-compact for every slope) if for every $x^{*}\in\mathcal{X}^{*}$ $x\mapsto f(x)-\langle x^{*}, x\rangle$ is an inf-compact (respectively sequentially-$\mathcal{T}$ inf-compact) function. We will say that the function $f$ is proper if its domain is nonempty and if the function $f$ does not take the value $-\infty $.  Recall that the Fenchel-Moreau conjugate $f^{*}$ is defined on  $\mathcal{X}^{*}$ by the formula \cite{96}: 
$$\displaystyle f^{*}(x^{*})=\sup_{x\in\mathcal{\mathcal{X}}} \langle x^{*}, x\rangle -f(x)\;.$$
\begin{Def} \label{def2.10} Given a sequence ${(M_{n})}_{n}$ of subsets of  $\mathcal{X}$,  its  $\mathcal{T}$-lower limit  and  $\mathcal{T}$-upper limit in the sense of Kuratowski \cite{0}, \cite{15}, \cite{7} Definition 4.10, are defined as follows:
\begin{displaymath}\liminf_{n}M_{n}=\{x\in\mathcal{X}: \forall U\in \mathcal{N}(x), \exists m\in {\N}: \forall n\geq m,\; U\cap M_{n}\neq\emptyset\}\;,\end{displaymath}
\begin{displaymath}\limsup_{n}M_{n}=\{x\in\mathcal{X}: \forall U\in \mathcal{N}(x), \forall m\in {\N},  \exists n\geq m: U\cap M_{n}\neq\emptyset\}\;.\end{displaymath}
\end{Def}
\begin{rem}\label{rem1} Clearly if  $\mathcal{X}$ is  a metric space with a distance $d$, for a subset $M$ of  $\mathcal{X}$, setting $d(x,M)=\inf\{d(x,m), m\in M\}$ if $M\neq \emptyset$, $+\infty$ if $M=\emptyset$, then:
$$\displaystyle\liminf_{n}M_{n}=\{x\in\mathcal{X}: \lim _{n}d(x,M_{n})=0\}\;,\quad  \displaystyle\limsup_{n}M_{n}=\{x\in\mathcal{X}: \liminf_{n}d(x,M_{n})=0\}.$$ \end{rem}
\begin{Def} \label{def2.11} (\cite{0}, \cite{7}) Given a sequence ${(f_{n})}_{n}$ of extended real-valued functions,  the upper  epi-limit  (resp lower  epi-limit) (also called $\Gamma$-limits) is defined as the function $\displaystyle \mathcal{T}-ls_{e}f_{n}$ (resp  $\displaystyle \mathcal{T}-li_{e}f_{n}$) whose epigraph is the lower limit  (resp upper limit) of the sequence of epigraphs of the $f_n$'s.
\end{Def}
The following formulas give analytic means for these limits (see \cite{7} Definition 4.1):
$$\left.\begin{array}{l}
\displaystyle \mathcal{T}-ls_{e}f_{n}(x)=\sup_{V\in\mathcal{N}(x)}\limsup_{n}\inf_{x^{'}\in V}f_{n}(x^{'})\;,\\\displaystyle \mathcal{T}-li_{e}f_{n}(x)=\sup_{V\in\mathcal{N}(x)}\liminf_{n}\inf_{x^{'}\in V}f_{n}(x^{'})\;.
\end{array}\right.$$
In the sequel, $(x_{n})\xrightarrow{\mathcal{T}} x$ means that the sequence ${(x_{n})}_{n}$ $\mathcal{T}$-converges to $x$; we will use frequently the next sequential version of the definition of the $\mathcal{T}$-lower epi-limit (the two notions coincide when every point has a countable base of neighbourhoods \cite{7}). 
\begin{Def}\label{def6.1} Given  a sequence ${(f_{n})}_n$ of $\overline{{\R}}$-valued functions defined on $( \mathcal{X}, \mathcal{T})$, and $\mathcal{Z}= \mathcal{X}^{{\N}}$,  the sequential $\mathcal{T}$-lower  epi-limit $f=seq\;\mathcal{T}-li_{e} f_{n}$ is defined by:
$$f(x)=\displaystyle\inf_{\{{(x_{n})}_{n}\in \mathcal{Z}:\;\; (x_{n})\xrightarrow{\mathcal{T}} x \}}\liminf_{n}f_{n}(x_{n}).$$
\end{Def} 
The links between epi-limits and conjugacy have been studied in the convex case in \cite{121}, \cite{30} and \cite{120}, \cite{0}, \cite{15}; moreover in \cite{112}, \cite{113}, \cite{114}  there is links between epi-limits, variational convergences, usual operations and conjugacy. In order to prove the main results of this article (in section 6), we need to obtain few preliminar properties of the Fenchel-Moreau conjugate of the sequential weak-lower epi-limit of a sequence of nonconvex functions (Proposition \ref{prop6.1}, Theorem \ref{thm7.7}, Corollary \ref{coro7.8}), $\mathcal{X}_\mathcal{T}$ will denote the topological space $(\mathcal{X}, \mathcal{T})$ and the symbol $\mathcal{W}^{*}_{b}( \mathcal{X}^{*},\mathcal{X}_\mathcal{T})$ (in abbreviate form $\mathcal{W}^{*}_{b}$) denotes the topology on  $\mathcal{X}^{*}$ of uniform convergence on the symmetric sequentially $\mathcal{T}$-compact sets of  $\mathcal{X}$. By similarity with the case where  $\mathcal{X}_\mathcal{T}$ is a Banach space
 endowed with its strong topology (see \cite{97} \S 18 D, \cite{110} Theorem 1.13)  it is called {\it bounded weak star topology} or {\it bounded weak$^{*}$ topology} .
\begin{prop} \label{prop6.1} Let $(\mathcal{X}, \mathcal{T})$ be a locally convex topological linear space. 
If   $f=seq\; \mathcal{T}-li_{e} f_{n}$ then:
$$f^{*}\leq  \mathcal{W}^{*}_{b}-ls_{e}f_{n}^{*}.$$
  \end{prop}
Proof of Proposition \ref{prop6.1}. 
 \begin{lem} \label{lem6.1} Let $\mathcal{K}$ be the set of all sequentially $\mathcal{T}$-compact sets.
Setting for $K\in \mathcal{K}$, $\iota_{K}(x)=0$ if $x\in K$, $+\infty$ if $x\notin K$, $f_{n}^{K}=f_{n}+\iota_{K}$, $f=seq\; \mathcal{T}-li_{e} f_{n}$, then:
$$f=\displaystyle\inf_{K\in \mathcal{K}} seq\; \mathcal{T}-li_{e} f_{n}^{K}.$$ 
\end{lem}
Proof of Lemma \ref{lem6.1}. Let $g=\displaystyle\inf_{K\in \mathcal{K}} seq\; \mathcal{T}-li_{e} f_{n}^{K}$. Since for every 
$K\in \mathcal{K}$, $f_{n}\leq f_{n}^{K}$, we have $f\leq seq\; \mathcal{T}-li_{e} f_{n}^{K}$, therefore $f\leq g$. Conversely, suppose that$f(x)< r$. Then there exists a  sequence $(x_{n})_n$   $\mathcal{T}$-converging to $ x $ such that $\displaystyle\liminf_{n}f_{n}(x_{n})<r.$ Setting $K=\{x\}\cup\{x_{n},\; n\in {\N}\}\in\mathcal{K}$, it is clear that
$g(x)\leq seq\; \mathcal{T}-li_{e} f_{n}^{K}(x)\leq\displaystyle\liminf_{n}f_{n}(x_{n})<r$, hence $g\leq f$.
\begin{lem} \label{lem6.3}    If $f=\displaystyle seq\; \mathcal{T}-li_{e} f_{n}$, then
  $f^{*}(x^{*})\leq \displaystyle \limsup_{n} f_{n}^{*}(x^{*})$ for all $x^{*}\in\mathcal{X}^{*}$.
\end{lem}
Proof of Lemma \ref{lem6.3}.  Let us suppose that $\displaystyle \limsup_{n} f_{n}^{*}(x^{*})<r$. Then for $n$ sufficiently large: $\langle x^{*},.\rangle-r<f_{n}$ and as a consequence: $\langle x^{*},.\rangle-r\leq f$ or equivalently $f^{*}\leq r$.$\;\Box$ 
 \begin{lem} \label{lem6.4} Let $\mathcal{K}$ and $f_{n}^{K}$ be as in Lemma \ref{lem6.1}. For every $K\in \mathcal{K}$, one has:
 $$\displaystyle \limsup_{n} (f_{n}^{K})^{*}\leq   \mathcal{W}^{*}_b-ls_{e} f_{n}^{*}.$$
\end{lem}
Proof of Lemma \ref{lem6.4}. Recall that given two extended real valued functions $f,\;g$ defined on  $\mathcal{X}$, the infimal convolution $f\Box g$ is (classically) defined by the formula:
$$(f\Box g) (x)=\inf_{y\in\mathcal{X}}f(x-y)+g(y)\;,$$
 one has $(f\Box g)^{*}=f^{*}+g^{*}$ and $(f+g)^{*}\leq f^{*}\Box g^{*}$.
For each sequentially $\mathcal{T}$-compact set $K$, setting $L=K\cup-K$, we get $$(f_{n}^{K})^{*}=(f_{n}+\iota_{K})^{*}\leq f_{n}^{*}\Box {\iota_{K}}^{*}\leq f_{n}^{*}\Box {\iota_{L}}^{*}\;.$$ Moreover the family of sets $V(L,\epsilon)=\{v^{*}:\;{\iota_{L}}^{*}(v^{*})\leq \epsilon \}$, $K\in \mathcal{K}$, $\epsilon>0$, is a base of $ \mathcal{W}^{*}_{b}$-neighbourhoods of the origin. For every $v^{*}\in V(L,\epsilon)$, we obtain:
$$(f_{n}^{K})^{*}(x^{*})\leq (f_{n}^{*}\Box {\iota_{L}}^{*})(x^{*})\leq f_{n}^{*}(x^{*}+v^{*})+\epsilon,$$
thus
$$(f_{n}^{K})^{*}(x^{*})\leq \displaystyle \inf_{v^{*}\in V(L,\epsilon)} f_{n}^{*}(x^{*}+v^{*})+\epsilon\;\;$$
and
$$\displaystyle \limsup_{n} (f_{n}^{K})^{*}(x^{*})\leq \displaystyle \limsup_{n}\displaystyle \inf_{v^{*}\in V(L,\epsilon)} f_{n}^{*}(x^{*}+v^{*})+\epsilon\leq   \mathcal{W}^{*}_b-ls_{e} f_{n}^{*}(x^{*})+\epsilon,$$
and since  $\epsilon$ is arbitrary, the proof of Lemma \ref{lem6.4} is complete.\;$\;\Box$ \\\
End of the proof of Proposition \ref{prop6.1}. From Lemma \ref{lem6.1} $f=\displaystyle\inf_{K\in \mathcal{K}} seq\; \mathcal{T}-li_{e}f_{n}^{K}$, hence due to  Lemmas \ref{lem6.3}, \ref{lem6.4} we get:
$$f^{*}=\displaystyle\sup_{K\in \mathcal{K}}(seq\; \mathcal{T}-li_{e} f_{n}^{K})^{*}\leq \displaystyle \sup_{K\in \mathcal{K}}\limsup_{n} (f_{n}^{K})^{*}\leq   \mathcal{W}^{*}_b-ls_{e} f_{n}^{*}\;.\;\;\;\Box$$
\begin{thm}\label{thm7.7} Let $(\mathcal{X}, \mathcal{T})$ be a locally convex topological linear space, $ \mathcal{W}^{*}_{b}$ be the topology on  $\mathcal{X}^{*}$ of uniform convergence on $\mathcal{T}$-compact sets of  $\mathcal{X}$.  Given a sequence ${(f_{n})}_n$ of elements of ${\overline{\R}}^{\mathcal{X}}$ bounded below by a function $h\in {\overline{\R}}^{\mathcal{X}}$ inf-sequentially  $ \mathcal{T}$-compact for every slope, then setting $f=seq\; \mathcal{T}-li_{e} f_{n}$, we have:
$$\displaystyle f^{*}=\limsup_{n}f_{n}^{*}= \mathcal{W}^{*}_b-ls_{e} f_{n}^{*}\;.$$
\end{thm}
\begin{coro}\label{coro99} Under the assumptions of Theorem \ref{thm7.7}, for any topology $\mathcal{\tau}$ stronger or equal than $\mathcal{W}^{*}_{b}$ we have:
$\displaystyle \mathcal{\tau}-ls_{e} f_{n}^{*}=f^{*}=\limsup_{n}f_{n}^{*}$.
\end{coro}
Proof of Corollary \ref{coro99}. From Theorem \ref{thm7.7},  we obtain the chain of inequalities
$$\limsup_{n}f_{n}^{*}= \mathcal{W}^{*}_{b}-ls_{e} f_{n}^{*}\leq \tau-ls_{e} f_{n}^{*}\leq \limsup_{n}f_{n}^{*}=f^{*}\;.\;\Box $$
Proof of Theorem \ref{thm7.7}. Let us first prove the following lemmas:
\begin{lem} \label{lem6.2} Let $K$ be a sequentially $ \mathcal{T}$-compact set. If ${(g_{n})}_n$ is a sequence of $\overline{{\R}}$-valued functions defined on $K$, and $g=seq\; \mathcal{T}-li_{e}g_{n}$,
then: $$\displaystyle\inf_{x\in K}g(x)\leq \liminf_{n}\inf_{x\in K}g_{n}(x)$$
\end{lem}
Proof of Lemma \ref{lem6.2}. Let $x_{n}\in K$ such $r_{n}=\displaystyle\inf_{x\in K}g_{n}(x)\geq g_{n}(x_{n})-\frac{1}{n}$, and let $\displaystyle r>\liminf_{n} r_{n}$. Extracting subsequences we can find  subsequences $(r_{n_{k}})_k$  such $\displaystyle\lim_{k}r_{n_{k}}<r$ and $(x_{n_{k}})_k$ $ \mathcal{T}$-converging to $ x \in K$. Then
$g(x)\leq \displaystyle\liminf_{k}g_{n_{k}}(x_{n_{k}})\leq \displaystyle\lim_{k} r_{n_{k}}+\frac{1}{n_{k}}<r.$$\;\Box$ 

\begin{lem} \label{lem6.66} Let $h:\mathcal{X}\to {\R}\cup\{\infty\}$ be a sequentially $ \mathcal{T}$-inf-compact function. Given a sequence ${(f_{n})}_n$ of elements of ${\overline{\R}}^{ \mathcal{X}}$ bounded below by $h$, setting $f=seq\; \mathcal{T}-li_{e} f_{n}$, we have:
$$f^{*}(0)=\displaystyle\limsup_{n} f_{n}^{*}(0)=  \mathcal{W}^{*}_b-ls_{e} f_{n}^{*}(0).$$
\end{lem}
Proof of Lemma \ref{lem6.66}. Let us first give the proof with the additional assumption: 
$$(H)\;\;\;\displaystyle\limsup_{n}\inf_{x\in\mathcal{X}}f_{n}(x)<\infty\;.$$
With $(H)$ let $r\in{\R}$ and ${(x_{n})}_{n}$ be a sequence such ${(f_{n}(x_{n}))}_{n}$ is bounded above by $r$. Since $f_{n}\geq h$ for every integer $n$, one has $\displaystyle\inf_{x\in\mathcal{X}}f_{n}(x)=\displaystyle\inf_{x\in h^{\leq r}}f_{n}(x)$ and  $h^{\leq r}$ is  nonempty sequentially $ \mathcal{T}$-compact. Extracting from  ${(x_{n})}_{n}$ a $ \mathcal{T}$-converging subsequence, we remark that $f^{\leq r}\neq \emptyset$. Moreover since $h$ is $ \mathcal{T}$ sequentially lower semicontinuous, we have $f\geq h$ and also $$\displaystyle\inf_{x\in\mathcal{X}}f(x)=\displaystyle\inf_{x\in f^{\leq r}}f(x)=\displaystyle\inf_{x\in h^{\leq r}}f(x)\;.$$ 
From Lemma \ref{lem6.2} we obtain: 
$$\displaystyle\inf_{x\in\mathcal{X}}f(x)=\displaystyle\inf_{x\in h^{\leq r}}f(x)\leq \liminf_{n}\inf_{x\in h^{\leq r}}f_{n}(x)=\liminf_{n}\inf_{x\in\mathcal{X}}f_{n}(x).$$
Therefore:
$$f^{*}(0)=-\displaystyle\inf_{x\in\mathcal{X}}f(x)\geq -\liminf_{n}\inf_{x\in\mathcal{X}}f_{n}(x)= \limsup_{n}-\inf_{x\in\mathcal{X}}f_{n}(x)=\limsup_{n}f_{n}^{*}(0)\;,$$
hence, from Proposition \ref{prop6.1}: 
$$f^{*}(0)\geq \limsup_{n}f_{n}^{*}(0)\geq  \mathcal{W}^{*}_b-ls_{e} f_{n}^{*}(0)\geq f^{*}(0).$$
Proof of Lemma \ref{lem6.66} without assumption $(H)$.
 Let us consider a subsequence  $(n_{k})_k$ such that $ \displaystyle\limsup_{n} f_{n}^{*}(0)=\lim_{k}f_{n_{k}}^{*}(0)$.   If $\infty=\displaystyle\limsup_{k}\inf_{x\in\mathcal{X}}f_{n_{k}}(x)=-\lim_{k}f_{n_{k}}^{*}(0)$,  using Proposition \ref{prop6.1} the conclusion stems from the relations 
 $$f^{*}(0)\geq \displaystyle\lim_{k} f_{n_{k}}^{*}(0)=\displaystyle\limsup_{n} f_{n}^{*}(0)=-\infty= \mathcal{W}^{*}_b-ls_{e}f_{n}^{*}(0)\geq f^{*}(0)\;.$$
Suppose now that $\displaystyle\limsup_{k}\inf_{x\in\mathcal{X}}f_{n_{k}}(x)<\infty$. Define
$\displaystyle g:=seq\; \mathcal{T}-li_{e} f_{n_{k}}\geq f$. Since the sequence $(f_{n_{k}})_k$ satisfies the assumption $(H)$, then from the first part of the proof and Proposition \ref{prop6.1} we have:
$$f^{*}(0)\geq g^{*}(0)=\displaystyle\limsup_{k} f_{n_{k}}^{*}(0)=\displaystyle\limsup_{n} f_{n}^{*}(0) \geq \mathcal{W}^{*}_b-ls_{e} f_{n}^{*}(0)\geq f^{*}(0)\;.$$ 
This ends the proof of Lemma \ref{lem6.66}. $\;\Box$ \\\\
End of the proof of Theorem \ref{thm7.7}. Given $x^{*}$ in  $\mathcal{X}^{*}$, set $g_{n}=f_{n}-\langle x^{*},.\rangle$ and $l=h-\langle x^{*},.\rangle$. Then $l$ is a sequentially $ \mathcal{T}$  inf-compact function, and for every $n$, $g_{n}\geq l$.  Moreover since $g_{n}^{*}(y^{*})=f_{n}^{*}(x^{*}+y^{*})$, $g_{n}^{*}(0)=f_{n}^{*}(x^{*})$, $\mathcal{W}^{*}_b -ls_{e} g_{n}^{*}(0)=  \mathcal{W}^{*}_b-ls_{e} f_{n}^{*}(x^{*})$ and $g=seq\; \mathcal{T}-li_{e} g_{n}=f- \langle x^{*},.\rangle$. Applying Lemma \ref{lem6.66} we obtain:
$$f^{*}(x^{*})=g^{*}(0)=\limsup_{n}g_{n}^{*}(0)=\limsup_{n}f_{n}^{*}(x^{*})= \mathcal{W}^{*}_b-ls_{e}g_{n}^{*}(0)= \mathcal{W}^{*}_b-ls_{e} f_{n}^{*}(x^{*}).$$
The proof of Theorem \ref{thm7.7} is complete.$\;\Box$ \\\\
The symbol  $\sigma$ denotes the weak topology  $\sigma(\mathcal{X}, \mathcal{X}^{*})$ on  $\mathcal{X}$. When  $\mathcal{X}$ is a Banach space, considering the case $\mathcal{T}=\sigma(, \mathcal{X},  \mathcal{X}^{*})$, due to  Eberlein-Smulian's theorem \cite{938} IV \S 5 section 3* Theorem 2, \cite{97} \S 18 Corollary A and Theorem B, the $\sigma$-compact sets are the sequentially $\sigma$-compact sets.  Krein's result (\cite{938} IV \S 5 section 5* Theorem 1, \cite{97} \S 19 Theorem E) asserts that the closed convex hull of a weakly compact set is weakly compact too. Therefore the topology $ \mathcal{W}^{*}_b$ on  $\mathcal{X}^{*}$  coincide with the Mackey topology $\tau( \mathcal{X}^{*},\mathcal{X})$ on  $\mathcal{X}^{*}$ (of uniform convergence on weakly convex compact sets of  $\mathcal{X}$) which is coarser than or equal to the topology of the dual norm $\Vert .\Vert_*$. An immediate consequence of Theorem \ref{thm7.7} and  Corollary \ref{coro99} is: 
\begin{coro} \label{coro7.8} Let $(\mathcal{X}, \Vert.\Vert)$ be a Banach space, $\tau^*$ be the Mackey topology $\tau(\mathcal{X}^{*},\mathcal{X})$ on $X^*$ and let $h:\mathcal{X}\to {\R}\cup\{\infty\}$ be an inf $\sigma$-compact function for every slope. Given a sequence ${(f_{n})}_n$ of elements of ${\overline{\R}}^{\mathcal{X}^{*}}$ bounded below by $h$, setting $f=seq\;\sigma-li_{e} f_{n}$ then 
$$\displaystyle f^{*}=\limsup_{n}f_{n}^{*}=\tau^{*}-ls_{e} f_{n}^{*}=\Vert.\Vert_{*}-ls_{e} f_{n}^{*}\;. $$ 
\end{coro}
The symbol  $\sigma^{*}$ denotes the weak star topology  $\sigma(\mathcal{X},\mathcal{X})$ on  $\mathcal{X}^{*}$. The case of $(\mathcal{X}^{*}, \sigma^{*})$ for the sequential lower epi-limit in convex case is treated in \cite{1} Theorem 7.5.1 and the  real lower epi-limit in convex case is considered in \cite{113} Theorem 2.
\begin{coro} \label{coro7.81} Let $(\mathcal{X}, \Vert.\Vert)$ be a separable Banach space, and a sequence ${(f_{n})}_n$ of elements of ${\overline{\R}}^{\mathcal{X}^{*}}$ bounded below by an inf $\sigma^{*}$-compact function for every slope, setting $f=seq\;\sigma^{*}-li_{e} f_{n}$ then for the duality between $\mathcal{X}^{*}$ and $\mathcal{X}$, with $f=seq\;\sigma^{*}-li_{e} f_{n}$, 
$$\displaystyle f^{*}=\limsup_{n}f_{n}^{*}=\Vert.\Vert-ls_{e} f_{n}^{*}\;. $$ 
\end{coro}
Proof of Corollary \ref{coro7.81}.  We take $\mathcal{T}=\sigma(\mathcal{X}^{*},\mathcal{X})$ on $X^{*}$. When  $\mathcal{X}$ is a separable Banach space it is well known that the  $\sigma^{*}$ compact sets are metrisable, thus the $\sigma^{*}$ sequentially compact sets are the $\sigma^{*}$ compact sets. Since the unit ball of $\mathcal{X}^{*}$ is a $\sigma^{*}$ compact set, the associated topology $\mathcal{W}^{*}_b$ on $\mathcal{X}$ is the topology of uniform convergence on bounded sets of $\mathcal{X}^{*}$, that is the norm topology of $\mathcal{X}$. $\;\Box$ 
 
\section{Biting Lemma, Biting convergence and tightness.}
\noindent In the sequel $(\Omega, \mathbb{T},\mu)$ is a measure space endowed with a $\sigma$-finite positive measure $\mu$ and with  a  tribe $\mathbb{T}$.
For a measurable subset $A\in\mathbb{T}$, set $A^{c}=\{\omega\in\Omega, \omega\notin A\}$ and $1_{A}$ stands for the characteristic function of $A$: $1_{A}(\omega)=1 $ if $\omega\in A$, 0 if $\omega\notin A$. 
Given two measurable $\overline{{\R}}$-valued functions $u$ and $v$ denote by $\{u\geq v\}$ the set $\{\omega\in \Omega: u(\omega)\geq v(\omega)\}$. 
Given a topological space $(E, \tau)$  with Borel tribe  $\mathcal{B}(E)$, the space $\mathcal{L}_{0}(E)$ is the space of $\tau$-measurable $E$-valued functions. We will say that the sequence $(x_{n})_n$ of elements of  $\mathcal{L}_{0}(E)$ converges almost everywhere to $x$ if there exists a negligible set $N$ such for every $\omega\in N^{c}$, the sequence $(x_{n}(\omega))_n$ $\tau$-converges to $x(\omega)$. Let $L_{0}(\Omega,E)$ be the space of classes of measurable functions (for  $\mu$-almost everywhere equality) defined on $\Omega$ and with values in $E$. It is customary to use the abuse of notation which consists to identify $x$ and its class $[x]$, we will do it. When the topology $\tau$ is associated to a distance $d$ the function defined on  $L_{0}(\Omega, E)^2$ by $\displaystyle d_{\mu}(x, y)=\int_{\Omega} \frac{d(x(\omega), y(\omega))}{1+d(x(\omega), y(\omega))} \alpha(\omega) d\mu(\omega)$ (where $ \alpha$ is any positive valued integrable function) is a distance, and the topology associated is the topology of convergence in local $\mu$-measure, that is convergence in measure on each set of finite $\mu$-measure. Hereafter $E$ is a separable Banach space. 
Let $L_{p}(\Omega,E, \mu)$, $1\leq p\leq\infty$, be the Lebesgue-Bochner space of classes of $p$-$\mu$ integrable functions ($\mu$-essentially bounded functions if $p=\infty$) defined on $\Omega$ with values in $E$ and endowed with its strong natural topology. When there is no ambiguity with respect to the measure we denote it by $L_{p}(\Omega,E)$. $\Vert x \Vert_{p}$ is the usual norm of an element $x$ of $L_{p}(\Omega,E)$, where for $1\leq p<\infty$, $\displaystyle\Vert x\Vert_{p}^{p}=\int_{\Omega}\Vert x\Vert^{p} d\mu$. 
Given a map $v:\Omega\to\overline{{\R}}$,  set $v^{+}=\sup(0,v)$ and  $v^{-}=(-v)^{+}$.
The upper integral $I_v$ or $\int_{\Omega}^{*}v d\mu$ of $v$ is defined by:
\begin{displaymath} I_{v}=\int_{\Omega}^{*}v d\mu=\inf\{\int_{\Omega}u  d\mu,
u\in L_{1}(\Omega,{\R}),u\geq v \;\; \mu-a.e\} =I_{v^{+}}-I_{v^{-}}\end{displaymath}
with the convention $+\infty -\infty=+\infty$.
A function $\phi: E\to \overline{{\R}}_{+}$ is a Young function if it is convex even continuous at $0$ with $\phi(0)=0$ and verifies $\lim_{\Vert e\Vert_{E}\to \infty}\phi(e)=+\infty$. When $\lim_{\Vert e\Vert_{E}\to \infty}\frac{\phi(e)}{\Vert e\Vert}=+\infty$, $\phi$ is said strongly coercive.
In the sequel we refer to a function $f:\Omega\times E\to \overline{{\R}}$ as an integrand. For every  $\omega\in \Omega$ let $f_{\omega}=f(\omega,\;.)$. 
An integrand $f$ is said to be {\it convex} respectively {\it even} if for every  $\omega\in \Omega$ the function $f_{\omega}$ is convex (respectively even). The integrand $f$ is said to be {\it measurable} if it is measurable when $\Omega\times E$ is endowed with the tribe $\mathbb{T}\otimes \mathcal{B}(E)$. Given two integrands $f$ and $g$ we write $f\leq g$ when there exists a negligible set $N$ such that for every $(\omega, e)\in N^{c}\times E$, $f(\omega, e)\leq g(\omega, e)$.\\A Young integrand $\phi$ is an integrand such that for every $\omega\in \Omega$, $\phi(\omega, .)$ is a Young function. If $\alpha$ is a positive valued integrable function, an integrand $\phi$ is said to be an $\alpha$-Young integrand if there exists a non decreasing convex strongly coercive function $\psi: {\R}_{+}\to {\R}_{+}$ verifying $\psi(0)=0$ and: $\phi(\omega,e)=\alpha(\omega)\psi(\alpha^{-1}(\omega).\Vert e\Vert)$. Given a Young integrand $\phi: \Omega\times E \to \overline{\R}_+$,  for every $t>0$ we denote $\phi_{t}(\omega,e)=\phi(\omega, te)$.
We will say that an integrand $f$ is {\it of Nagumo type}  if  $f$ is non negative, measurable and if  for every $\omega\in\Omega$ $f_{\omega}$ {\it is  inf-$\sigma$-compact for every slope}. When $E$ is reflexive every  $\alpha$-Young integrand is a  Nagumo integrand. Given $x\in L_{0}(\Omega,E)$, and an integrand $f$, we will denote by $f(x)$ the map $\omega\mapsto f(\omega,x(\omega))$. 
The integral functional $I_f$ associated to an integrand $f$ is the functional defined at some point $x$ of $L_{0}(\Omega,E)$ by: 
$$I_{f}(x)=I_{f(x)}= \int_{\Omega}^{*}f(x) d\mu\;.$$  
Let us recall the following definitions. 
\begin{Def} \label{def3.1} (\cite{2}) Let $X$ be a subset of $L_{p}(\Omega,E)$ with $1\leq p<\infty$. 
$X$ is said to be $p$-equi-integrable (equi-integrable if $p=1$), if for every positive number  $\epsilon$, there exist a positive constant $\eta$ and a measurable set $K$ of finite measure such that:\\
(1) $\displaystyle\sup_{ x\in X}\Vert x1_{A}\Vert_{p}<\epsilon$ for every  measurable set $A$ satisfying $\mu(A)\leq\eta$. \\
(2) $\displaystyle\sup_{ x\in X}\Vert x1_{K^{c}}\Vert_{p}<\epsilon$.
\end{Def}
Some authors use only the first property as definition of equi-integrability see \cite{78}, Definition 2.23 for example. In order to measure the lack of equi-integrability let us introduce the following notion:
\begin{Def}\label{defw6.102} (see \cite{67}) Let $\Sigma$ be the collection of all decreasing sequences $\sigma=(S_{k})_{k}$ of  measurable sets $S_{k}$ with a negligible intersection.
Given a sequence $(u_{n})_n$ of measurable $\overline{\R}$-valued functions it is convenient to use the index of equi-integrability $\delta^{+}_{\mu}((u_{n})_n)$ defined by:
$$ \delta_{\mu}^{+}((u_{n})_{n})=\displaystyle\sup_{\sigma\in \Sigma, \;\sigma=(S_{k})_{k}}\limsup_{k}\sup_{n\geq k}\int^{*}_{S_{k}} u_{n} d\mu\;.$$
When there is no ambiguity on the measure we note $\displaystyle\delta^{+}(( u_{n})_{n}):=\delta_{\mu}^{+}((u_{n})_{n})$
\end{Def}
Recall the following result about the index (of equi-integrability):
\begin{thm}\label{prop991} (see \cite{67} Proposition 1.7 and Corollary 1.9)\\
$(a)$ A sequence $(u_{n})_n$ of integrable functions is equi-integrable in the sense of Definition \ref{def3.1} if and only if $\displaystyle \delta^{+}((\Vert u_{n}\Vert )_{n})=0$.\\
$(b)$ Suppose the measure $\mu$ is atomless. A sequence $(u_{n})_n$ of measurable functions is eventually uniformly integrable if and only if $\displaystyle \delta^{+}((\Vert u_{n}\Vert )_{n})=0$.
\end{thm}
\begin{thmD} \label{def33.2}  De la Vall\'ee Poussin's Theorem ( \cite{37} page 33, \cite{98}, \cite{9} II-5,\cite{10} Chapter 2, Theorems 22 and 25, see \cite{38} Theorem 16.8). A subset $X$ of $L_{1}(\Omega,E)$ is said uniformly integrable when one of the following equivalent conditions is satisfied:\\ 
$(a)$ $X$ is bounded and equi-integrable in  $L_{1}(\Omega,E)$, \\
$(b)$ There exists a positive valued integrable function $\alpha$ and a $\alpha$-Young integrand $\phi=\psi_{\alpha}$ such that $\sup_{x\in X} I_{\phi}(x)\leq 1$.\\
$(c)$ For every integrable positive valued function $\beta$, $\displaystyle\lim_{n}\sup_{x\in X} \int_{\{\Vert x\Vert\geq n\beta\}}\Vert x\Vert d\mu=0\;.$\\
When the measure $\mu$ is finite,(see \cite{9} II-5), in the assertion $(b)$ and $(c)$ it suffices to consider the case $\alpha=1=\beta$.
\end{thmD}
Proof. The equivalences $(a)\Leftrightarrow (b)\Leftrightarrow (c)$ are well-known when the measure is finite and $\alpha=1$. In the $\sigma$-finite case, for every positive valued integrable function $\alpha$, first let us remark that since the $\mu$-negligible sets are the $\alpha\mu$-negligible sets, taking the assertion $(a)$ has definition of uniform integrability, from Thorem \ref{prop991} $(a)$, we deduce that $X$ is uniformly integrable in  $L_{1}(\Omega,E, \mu)$ if and only if  $\alpha^{-1} X$ is uniformly integrable in $L_{1}(\Omega,E, \alpha\mu)$. Therefore $(a)\Leftrightarrow (b)$. Moreover  For every integrable positive valued function $\beta$, and for every integer $n$ and $x\in X$,
$$\int_{\{\Vert x\Vert\geq n\beta\}}\Vert x\Vert d\mu=\int_{\{\Vert\beta^{-1} x\Vert\geq n\}}\Vert\beta^{-1} x\Vert \beta d\mu,$$
thus  $X$ verifies $(c)$ in $L_{1}(\Omega,E, \mu)$ if and only if  $\beta^{-1} X$ verifies $(c)$ with "$\beta=1$" in $L_{1}(\Omega,E, \beta\mu)$. Therefore  $(a)\Leftrightarrow (c)$. $\;\Box$\\\\
Let us give now the statement of the Biting Lemma (valid in the $\sigma$-finite case):
\begin{thm} \label {thm 1} (\cite{705} Lemma 3.4, \cite{3} Theorem 6.1.4).  Let $(x_{n})_n$ be a bounded sequence in $L_{1}(\Omega,E)$. There exists a subsequence 
$(x_{n_{k}})_k$ and a decreasing sequence  $(A_{k})_k$ of measurable sets with a negligible intersection such that the sequence $(x_{n_{k}}1_{A^{c}_{k}})_k$ is equi-integrable 
\end{thm}
Proof of Theorem \ref{thm 1}. The result is true when $\mu$ is finite valued (\cite{705} Lemma 3.4, \cite{3} Theorem 6.1.4). If $\mu$ is $\sigma$-finite there exists a positive valued integrable function $\alpha$. The measure $\nu(A)=\int_{A} \alpha d\mu$ is finite valued.  Let $(x_{n})_n$ be a bounded sequence in $L_{1}(\Omega,E, \mu)$. The sequence  $(\alpha^{-1}x_{n})_n$ is a bounded sequence in $L_{1}(\Omega,E, \nu)$. There exists a subsequence 
$(\alpha^{-1}x_{n_{k}})_k$ and a decreasing sequence  $(A_{k})_k$ of measurable sets with a $\nu$-negligible intersection such that the sequence $(\alpha^{-1}x_{n_{k}}1_{A^{c}_{k}})_k$ is $\nu$-equi-integrable. The $\nu$-negligible sets being exactly the $\mu$-negligible sets we deduce with Theorem \ref{prop991} $(a)$ that the sequence $(x_{n_{k}}1_{A^{c}_{k}})_k$ is $\mu$-equi-integrable. $\;\Box$
\begin{Def} \label{def6.0} Given a sequence ${(x_{n})}_n$ of $E$-valued measurable functions defined on $\Omega$, we will say that $(x_{n})_n$ converges to a measurable $E$-valued  function $x$ in the Biting sense if there exists an increasing covering (up to a negligible set) ${(\Omega_{k})}_k$ of $\Omega$  by measurable sets such that for all $k$ the sequence ${(x_{n}\vert\Omega_{k})}_n$  of restrictions to $\Omega_{k}$ $\sigma(L_{1}(\Omega_{k},E, \mu), L_{1}(\Omega_{k},E,\mu)^{*})$-converges to $x\vert\Omega_{k}$.  
\end{Def}
\begin{coro} \label{coro 658} (\cite{78}, \cite{3} Remark 6.1.5), Suppose $E$ is reflexive. Then any bounded sequence in $L_{1}(\Omega,E)$ admits a Biting converging subsequence to an integrable function.
\end{coro}
Proof: It is well-known (\cite{150}) that in this case $L_{1}(\Omega,E)^{*}= L_{\infty}(\Omega,E^{*})$, and in $L_{1}(\Omega,E)$ any bounded equi-integrable sequence admits a weakly converging subsequence. Applying the first part of Theorem \ref{thm 1} and  extracting from $(x_{n_{k}}1_{A^{c}_{k}})_k$ a weakly converging sequence $(x_{n_{k_{l}}}1_{A^{c}_{k_{l}}})_l$ to an integrable function $x$, the sequence $(x_{n_{k_{l}}})_l$ converges in the Biting sense: the restrictions on each $A^{c}_{k_{l}}$ converge weakly to the restriction of the function $x$. And for every integer $k$ 
$$\Vert x1_{\Lambda_{k}}\Vert_{1}\leq \liminf_{n}\Vert x_{n}1_{\Lambda_{k}}\Vert_{1}\leq \liminf_{n}\Vert x_{n}\Vert_{1}\;,$$
Since the norm is weakly semicontinuous and the sequence is norm bounded:
$$\Vert x\Vert_{1}\leq \liminf_{k} \Vert x1_{\Lambda_{k}}\Vert_{1}\leq \liminf_{n}\Vert x_{n}\Vert_{1}\;<\infty\;.\;\; \Box$$
\begin{coro} \label{lem3.11} Suppose $E$ is reflexive. Let $1\leq p\leq \infty$ and $\beta$ be a measurable positive valued function. Every bounded sequence in $L_{p}(\Omega, E, \beta\mu)$ admits a Biting converging subsequence with a Biting limit in $L_{p}(\Omega, E, \beta\mu))$.
\end{coro}
Proof of Corollary \ref{lem3.11}. First let us give a proof when $\beta=1$. When $p=1$ it is exactly the Biting Lemma in $\sigma$-finite case above. When $\Omega$ is of finite measure since $L_{p}(\Omega, E)$ is topologically included in $L_{1}(\Omega,E)$, the result of the existence of a Biting-converging subsequence is an immediate consequence of the classical Biting Lemma. First let us show that the sequence admits a Biting convergent subsequence in $\sigma$-finite case. There exists a bounded integrable function $\alpha$ positive valued. Let a bounded sequence  $(x_{m})_m$ of $p$-integrable elements. Set $y_{m}= \alpha^{\frac{-1}{p}}x_{m}$, then the sequence  $(y_{m})_m$ is bounded in $L_{p}(\Omega, E, \alpha\mu)$ thus in $L_{1}(\Omega, E, \alpha\mu)$. Due to the Biting Lemma there exists $y\in L_{1}(\Omega_{k}, E,\alpha\mu)$  an increasing covering $(\Omega_{k})_k$ of $\Omega$ by measurable sets, a subsequence  $(y_{m_n})_n$ such for all integer $k$ the sequence of the restrictions of the $y_{m_n}$ to $\Omega_{k}$ weakly converges in $L_{1}(\Omega_{k}, E,\alpha\mu)$ to the restriction to $\Omega_{k}$ of  $y$. That is for all $x^{*}\in L_{\infty}(\Omega_{k}, E^{*} )$,
$$\lim_{n}\int_{\Omega_{k}} \langle y_{m_n}, x^{*}\rangle\alpha d\mu=\int_{\Omega_{k}} \langle y, x^{*}\rangle\alpha d\mu$$
equivalently:
$$\lim_{n}\int_{\Omega_{k}} \langle x_{m_n}\alpha^{1-\frac{1}{p}}, x^{*} \rangle d\mu=\int_{\Omega_{k}}\langle  y\alpha, x^{*}\rangle d\mu.$$
Setting
$\Lambda_{k}=\{\omega\in \Omega_{k}:\alpha(\omega)\geq \frac{1}{k}\}$, it is clear that $(\Lambda_{k})_k$ is an increasing covering of $\Omega$ and for every $x^{*}\in L_{\infty}(\Lambda_{k}, E^{*})$, $y^{*}=\alpha^{\frac{1}{p}-1}.x^{*}$ is an element of $L_{\infty}(\Lambda_{k}, E^{*})$, thus:
$$\lim_{n}\int_{\Lambda_{k}} \langle x_{m_n}, x^{*}\rangle  d\mu=\lim_{n}\int_{\Lambda_{k}} \langle x_{m_n}, \alpha^{1-\frac{1}{p}}y^{*}\rangle  d\mu=\int_{\Lambda_{k}} \langle y\alpha, y^{*}\rangle d\mu=\int_{\Lambda_{k}} \langle y\alpha^{\frac{1}{p}}, x^{*}\rangle d\mu\;.$$
This proves that the sequence  $(x_{m_n})_n$ Biting converges to $x=y\alpha^{\frac{1}{p}}$. The assertion $x\in L_{p}(\Omega, E)$  is a consequence of the Corollary \ref{3.11} bellow. The proof of Corollary \ref{lem3.11} is complete in the case $\beta=1$. \\ Let  $\beta$ be a measurable positive valued function, the measure $\nu=\alpha\mu$ is $\sigma$-finite. If the sequence  $(x_{m})_m$ is bounded in $L_{p}(\Omega, E, \nu)$ applying the first part of the proof we deduce the existence of a Biting-converging subsequence $(x_{m_n})_n$ to an element of $L_{p}(\Omega, E, \nu)$. That is there exists an increasing covering ${(\Omega_{k})}_k$ of $\Omega$  by measurable sets such that for all $k$ the sequence ${(x_{n}\vert\Omega_{k})}_n$  of restrictions to $\Omega_{k}$ $\sigma(L_{1}(\Omega_{k},E, \nu), L_{1}(\Omega_{k},E,\nu)^{*})$-converges to $x\vert\Omega_{k}$. Define for all positive integer $k$, $\Lambda_{k}=\Omega_{k}\cap\{\omega: k^{-1}\leq \beta(\omega)\leq k\}$. Since $\beta$ is positive valued, the $\nu$-negligible sets are the $\mu$-negligible sets. Moreover ${(\Lambda_{k})}_k$ is an increasing covering of $\Omega$  by measurable sets such that for all $k$ the sequence ${(x_{n}\vert\Lambda_{k})}_n$  of restrictions to $\Lambda_{k}$ $\sigma(L_{1}(\Lambda_{k},E, \mu), L_{1}(\Lambda_{k},E,\mu)^{*})$-converges to $x\vert\Lambda_{k}$. $\;\Box$
\begin{prop}\label{919}  Any convex $C\subset L_{0}(\Omega, E)$ closed for the convergence in local measure is Biting sequentially closed.
\end{prop}
The proof of Proposition \ref{919} is an immediate consequence of the following Lemma:
\begin{lem} \label{l1} If $(x_{m})_m$ Biting-converges to $x$ there exists a sequence  $(\overline{x}_{m})_m$ converging in local measure to $x$ such that for all integer $m$,
$\overline{x}_{m}\in co\{x_{n},\;n\geq m\}$
\end{lem}
Proof of Lemma\ref{l1}. Given $\alpha$ a positive valued integrable function less than $1$, the topology of convergence in local measure is defined by the distance:
$$\displaystyle d_{\mu}(x, y)=\int_{\Omega} \frac{\Vert x(\omega)- y(\omega)\Vert}{1+\Vert x(\omega)- y(\omega)\Vert} \alpha(\omega) d\mu(\omega)\;.$$
Since the sequence  $(x_{m})_m$ Biting converges to $x$ there exists an increasing covering  $(\Omega_{k})_k$ of $\Omega$ such that for every integer $k$ the sequence  $(x_{m}1_{\Omega_{k}})_m$ weakly converges in $L_{1}(\Omega, E)$ to $x1_{\Omega_{k}}$.
Pick an integer $k_m$ such that $\int_{\Omega^{c}_{k_{m}}} \alpha d\mu< \frac{1}{3m}$. Due to Mazur's Lemma, $x1_{\Omega_{k_{m}}}$ is in the strong closure of $co\{x_{n}1_{\Omega_{k_{m}}}, \;n\geq m\}$, thus there exists an element
$\overline{x}_{m}\in co\{x_{n}, \;n\geq m\}$ such that $\Vert (x-\overline{x}_{m})1_{\Omega_{k_{m}}}\Vert_{1}<\frac{1}{3m}$. But
$$d_{\mu}(x, \overline{x}_{m})\leq d_{\mu}(x, x1_{\Omega_{k_{m}}})+ d_{\mu}(x1_{\Omega_{k_{m}}}, \overline{x}_{m}1_{\Omega_{k_{m}}})+d_{\mu} (\overline{x}_{m}1_{\Omega_{k_{m}}}, \overline{x}_{m})$$
therefore:
$$d_{\mu}(x, \overline{x}_{m})\leq 2\sup_{z}d_{\mu}(z1_{\Omega^{c}_{k_{m}}}, 0)+ d_{\mu}(x1_{\Omega_{k_{m}}}, \overline{x}1_{\Omega_{k_{m}}})\;,$$
and we get:
$$d_{\mu}(x, \overline{x}_{m})\leq 2\int_{\Omega^{c}_{k_{m}}} \alpha d\mu +\int_{\Omega_{k_{m}}} \frac{\Vert x -\overline{x}_{m} \Vert}{1+\Vert x -\overline{x}_{m} \Vert} \alpha d\mu<2\frac{1}{3m}+\Vert (x-\overline{x}_{m})1_{\Omega_{k_{m}}}\Vert_{1}<3\frac{1}{3m}=\frac{1}{m}\;.$$
This proves that the sequence $(\overline{x}_{m})_m$ converges in local measure to $x$. The proof of Lemma \ref{l1} is complete.$\;\Box$\\\\
Since by Fatou's Lemma every norm closed ball of $L_{p}(\Omega, E)$ is closed for the convergence in local measure it follows:
\begin{coro} \label{3.11}  Let $1\leq p\leq \infty$. The norm closed balls of $L_{p}(\Omega, E)$ are sequentially Biting closed.
\end{coro}
\begin{Def}\label{def6.202} (see \cite{111}) A subset $X$ of $L_{0}(\Omega,E)$ is said to be Nagumo tight if there exists a Nagumo integrand $h$ such $\;\displaystyle\sup_{x\in X} \int_{\Omega}h(x)d\mu<\infty$.
\end{Def}
\begin{rem} \label{rem 1}\end{rem} If $M$ is a $E$-valued multifunction with nonempty weakly compact values and with  $\mathbb{T}\otimes \mathcal{B}(E)$-measurable graph, the integrand $f(\omega,e)=\iota_{M(\omega)}(e)$ is measurable, then $f$ is an example of Nagumo integrand. As a consequence, for a such multifunction, the set of measurable almost everywhere selections of  $M$ is an example of Nagumo tight set.\\
Let us give another useful example of Nagumo  tight sets. 
\begin{prop}\label{prop00} If $E$ is a reflexive Banach space, then every converging sequence in the Biting sense is Nagumo tight.
\end{prop}
Proof of the Proposition \ref{prop00}. When $E$ is reflexive, every integrand of $\alpha$-Nagumo type is a Nagumo integrand.  By the Dunford-Pettis Theorem, every weakly compact set $X$ in $L_{1}(\Omega,E)$ is uniformly integrable, thus from Theorem \ref{def33.2} $(b)$, $X$ is Nagumo tight. Therefore every weakly converging sequence in $L_{1}(\Omega,E)$ is Nagumo tight.
Let a sequence ${(x_{n})}_n$ converging in the Biting sense. If ${(\Omega_{k})}_k$ is the sequence of measurable sets  appearing in 
Definition \ref{def6.0},  let $\Lambda_{0}=\Omega_{0}$ and for $k\geq 1$, $\Lambda_{k}=\Omega_{k}\backslash\Omega_{k-1}$.
Since every weakly compact set $X$ in $L_{1}(\Lambda_{k},E)$ is uniformly integrable, for each integer $k$ there exists a Nagumo integrand $h_k$ defined on $\Lambda_{k}\times E$ such that: $\displaystyle \sup_{n}\int_{\Lambda_{k}} h_{k}(x_{n})d\mu\leq 2^{-k-1}$. Define $\overline{h_{k}}(\omega,e)=h_{k}(\omega,e)$ if $\omega\in \Lambda_{k}$ and $0$ if not; $\displaystyle h=\sum_{k} \overline{ h_{k}}$, then $h$ is a 
Nagumo integrand which satisfies\\ 
\centerline{$\displaystyle \sup_{n}\int_{\Omega} h(x_{n})d\mu\leq \sum_{k}2^{-k-1}=1\;.$} This proves that the sequence ${(x_{n})}_n$ is Nagumo tight and ends the proof of the Proposition.$\;\Box$
\begin{Def} \label{def733} (see \cite{3} Section 6 and Lemma 6.1.1) A subset $X$ of $L_{0}(\Omega,E)$ is said to be weakly flexibly tight if  for every set $A$ of finite measure, for every $\epsilon>0$, there exists a measurable multifunction (see Definition \ref{def737}) $M_{\epsilon}$ with nonempty weakly compact values in $E$ such that for all $x\in X$,
$$\mu(A\cap\{\omega\in \Omega: x(\omega)\notin M_{\epsilon}(\omega)\})\leq \epsilon\;.$$
\end{Def}
The following result makes a link between  Nagumo tightness and the notion of weak flexible tightness used in \cite{3} Theorem 8.1.6.
\begin{prop} \label{prop 919} Let $E$ be a separable Banach space. A countable weakly  flexibly tight subset of $L_{0}(\Omega,E)$ is Nagumo tight.
\end{prop}
Proof of Proposition \ref{prop 919}.  Suppose that $X=\{x_{n}, n\in {\N}\}$ is a countable weakly  flexibly tight subset.  First, let us prove the following lemma:
\begin{lem} \label{lem 19} The following assertions hold:\\
$(a)$ For every measurable set $A$ of finite measure, for each $\epsilon>0$ there exist a measurable multifunction $L_{\epsilon}$ with weakly compact values a measurable set $A_{\epsilon}\subset A$ with $\mu(A_{\epsilon})\leq \epsilon$ and for all $x\in X$, for all $\omega\in A_{\epsilon}^{c}$, $x(\omega)\in L_{\epsilon}(\omega)$.\\
$(b)$ Proposition \ref{prop 919} holds if $\Omega$ is of finite measure. More precisely there exists a measurable multifunction $M$ with weakly compact values such that every $x\in X$ is an almost everywhere selection of $M$.
\end{lem}
Proof of Lemma \ref{lem 19}. Let $\epsilon>0$. For every integer $n$ keep the multifunction  $M_{\frac{\epsilon }{2^{n+1}}}$ of Definition  \ref{def733}. The measurable sets $A_{n}=\{\omega\in A: x_{n}(\omega)\notin M_{\frac{\epsilon }{2^{n+1}}}\}$ and $A_{\epsilon}=\bigcup_{n}A_{n}$ satisfies $\mu(A_{n})\leq \frac{\epsilon }{2^{n+1}}$ and
$\mu(A_{\epsilon})\leq \epsilon$. If $L_{\epsilon}=\bigcap_{n} M_{\frac{\epsilon }{2^{n+1}}}$, then $L_{\epsilon}$ is measurable with weakly compact values. $L_{\epsilon}$ and $A_{\epsilon}$ satisfy the assertion $(a)$ of Lemma \ref{lem 19}.\\
Suppose $\Omega$ has finite measure. In order to prove the second assertion, let us build a multifunction $M$ in the following way. Set $A_{0}=\Omega$.
From the first assertion, there exist a measurable set $A_{1}$ such $\mu(A_{1})\leq \frac{\mu(\Omega)}{2}$ a measurable multifunction $K_{1}$ with weakly compact values such for every $\omega\in A_{0}\backslash A_{1}$, 
for every $x\in X$, $x(\omega)\in K_{1}(\omega)$. Suppose build $A_{n}$ and $K_{n}$ such $\mu(A_{n})\leq \frac{\mu(\Omega)}{2^{n}}$, $K_{n}$ is measurable with weakly compact values such for every $\omega\in A_{n-1}\backslash A_{n}$, for every $x\in X$, $x(\omega)\in K_{n}(\omega)$. From  assertion $(a)$, applied to the measure space $A_{n}$, there exist a measurable set $A_{n+1}\subset A_{n}$ a measurable multifunction $K_{n+1}$ with weakly compact values defined on $A_{n}$ such $\mu(A_{n+1})\leq \frac{\mu(\Omega)}{2^{n+1}}$ and such for all $\omega\in A_{n}\backslash A_{n+1}$,  for every $x\in X$, $x(\omega)\in K_{n+1}(\omega)$. Set $A_{\infty}=\bigcap_{n}A_{n}$ (then $\mu(A_{\infty})=0)$, 
$M(\omega)=K_{n}(\omega)$ if $\omega\in A_{n-1}\backslash A_{n}$, $\{0\}$ if $\omega\in A_{\infty}$. Then $M$ has nonempty weakly compact values,  is measurable in sense of Definition \ref{Def 117}, and by construction for every $\omega \in A^{c}_{\infty}$, for every  $x\in X$,  $x(\omega)\in M(\omega)$. Due to the Remark \ref{rem 1}, the set $X$ is Nagumo tight. $\;\Box$\\
End of the proof of Proposition \ref{prop 919}. When the measure $\mu$ is $\sigma$-finite, let $(\Omega_{p})_p$ be an increasing covering of  $\Omega$ by measurable sets of finite measure. Using the above lemma, for each integer $p$ there exists a measurable  multifunction $M_{p}$ with nonempty weakly compact values defined on $\Omega_{p}$ such for every  $x\in X$, for almost every $\omega\in \Omega_{p}$, $ x(\omega)\in M_{p}(\omega)$. Define then the multifunction $M$ by $M(\omega)=M_{0}(\omega)$ if $\omega\in \Omega_{0}$, $M(\omega)=M_{p}(\omega)$ if $\omega\in \Omega_{p}\backslash \Omega_{p-1}$. By construction the multifunction $M$ is measurable with weakly compact values. Moreover since $X$ is countable, every element of $X$ is an almost everywhere selection of
$M$ and  Remark \ref{rem 1} allows to conclude.$\;\Box$

\section{Measurability and polarity }
In the sequel the tribe $\mathbb{T}$ is supposed to be $\mu$-complete and $E$ is a locally convex space.
The epigraph multifunction of an integrand $f:\Omega\times E\to{\overline{\R}}$ is the multifunction $epi\;f$ defined by $epi\;f(\omega)=\displaystyle epi f_{\omega}$. Let $E^{*}$ be the topological dual of $E$. The Fenchel-Moreau conjugate of $f$ is the integrand $f^{*}$ defined on $\Omega\times E^{*}$ by $\displaystyle f^{*}(\omega,e^{*})={(f_{\omega})}^{*}(e^{*})$. 
A Suslin space is a continuous image of a metrisable separable complete space. Let $(E,\tau)$ be a Suslin locally convex space often denoted by $E_{\tau}$.  Let $\sigma=\sigma(E, E^{*})$ be the weak topology on $E$. On the topological dual $E^{*}$, let  $\sigma^{*}=\sigma(E^{*}, E)$ be the weak star topology on $E^{*}$. Notice that if $E$ is a  separable  Banach space, then  $E_{\sigma}$ and $E^{*}_{\sigma^{*}}$ are examples of Suslin locally convex spaces. The tribe $\mathcal{B}(E)$ is the Borel tribe of $E_\tau$,  $\mathcal{B}(E_{\sigma})$ is the Borel tribe of $E_{\sigma}$ and $\mathcal{B}(E^{*}_{\sigma^{*}})$ is the Borel tribe of $E^{*}_{\sigma^{*}}\;$.  A Suslin space is hereditary Lindel\"of (that is with the property that every covering of an open set by open sets admits a countable subcovering). Therefore if $E$, $F$ are Suslin spaces then $\mathcal{B}(E\times F)=\mathcal{B}(E)\times \mathcal{B}(F)$. Recall the following notions. 
\begin{Def} (\cite{100}, \cite{4}  VII. 1) \label{def737} Let $F$ be a Suslin space. A  closed and $F$-valued multifunction $M$ defined on $\Omega$ is said to be measurable if its graph is $\mathbb{T}\otimes \mathcal{B}(F)$ measurable.
\end{Def}
\begin{Def} \label {Def 117} ( \cite{4} VII. 1) Let $(E,\tau)$ be a Suslin space and let $f:\Omega\times E\to{\overline{\R}}$ be an integrand. The integrand $f$ is called normal on $\Omega\times E_{\tau}$ if it is $\mathbb{T}\otimes \mathcal{B}(E_{\tau})$ measurable and for every $\omega\in\Omega$ the function $f_{\omega}$ is $\tau$-lower semicontinuous.
\end{Def}
\begin{prop} \label{prop 31} (see \cite{100} Proposition 2 and \cite{4} Corollary VII-2) Let $E$ be a separable  Banach space and let $f:\Omega\times E\to{\overline{\R}}$ be a $\mathbb{T}\otimes \mathcal{B}(E)$ measurable integrand. The Fenchel-Moreau conjugate integrand $f^{*}$ is a convex normal integrand on $\Omega\times E^{*}_{\sigma^{*}}$, and  the restriction of the biconjugate $f^{**}$ to $\Omega\times E$ is a convex normal integrand on $\Omega \times E_{\sigma}$.
\end{prop}
Proof of Proposition \ref{prop 31}. It is an immediate consequence of the following lemma. 
\begin{lem} \label{prop 315} (see \cite{890} Lemma 8, \cite{4} Corollary VII-2) Let $(E, \tau)$ be a Suslin locally convex space with topological dual $E^{*}$
and let $f:\Omega\times E\to{\overline{\R}}$ be a  $\mathbb{T}\otimes \mathcal{B}(E_{\tau})$-measurable integrand. The Fenchel-Moreau conjugate integrand $f^{*}$ is a convex normal integrand on $\Omega\times E^{*}_{\sigma^{*}}\;$. 
\end{lem}
By the use of the projection Theorem and a Castaing representation of the epigraph multifunction of $f$, the proofs of \cite{890} Lemma 8, \cite{4} Corollary VII-2 are valid even when $f$ is only supposed to be a $\mathbb{T}\otimes \mathcal{B}(E_{\tau})$-measurable integrand.
When $E$ is equipped with the norm topology it is a Suslin locally convex space. Moreover $E_{\sigma^{*}}^{*}$ is a Suslin locally convex space space too with topological dual $E$. Therefore the two parts of Proposition \ref{prop 31} are easy consequences of Lemma \ref{prop 315}. The proof of Proposition \ref{prop 31} is complete.$\;\Box$

\begin{Def}\label{def6.201} A sequence  $(f_{n})_n$ of extended real valued integrands is said to be quasi inf-$\sigma$-compact for every slope if for every  $\omega\in \Omega$ the sequence of functions $(f_{n}(\omega, .))_n$ is eventually bounded below by a function that is inf-$\sigma$-compact for every slope. 
\end{Def} 
\begin{thm} \label{thm 5.1}  Let $E$ be a separable Banach space and ${(f_{n})}_n$ be a sequence of $\mathbb{T}\otimes \mathcal{B}(E)$ measurable integrands defined on $\Omega\times E$.  The  pointwise weak lower sequential epi-limit integrand $f$ is defined for all $(\omega,e)\in \Omega\times E$ by $f(\omega,e)=seq\;\sigma-li_{e}{f_{n}}_{\omega}(e)$. Under one of the following assumptions the integrand  $f^{*}$ is normal on $\Omega\times E^{*}_{\sigma^{*}}$ (and moreover $f^{**}$ is normal on $\Omega\times E_{\sigma}$): \\
$(a)$ The sequence ${(f_{n})}_n$ is quasi inf-$\sigma$-compact for every slope in the sense of Definition \ref{def6.201}.\\
$(b)$ The Banach $E$ is reflexive and separable.
\end{thm}
Proof of Theorem \ref{thm 5.1}.
Due to Proposition \ref{prop 31} it suffices to proves that $f^{*}$ is normal on $\Omega\times E^{*}_{\sigma^{*}}$. We begin with the following lemmas:
\begin{lem} \label{lem 5.3} Let $(E, \Vert\;.\Vert)$ be a separable Banach space  and let $(f_{n})_n$ be a sequence of $\mathbb{T}\otimes \mathcal{B}(E)$ measurable integrands defined on $\Omega\times E$ which is quasi inf-$\sigma$-compact for every slope. If  $f=seq\;\sigma-li_{e}f_{n}$, then  $f^{*}$ is normal on $\Omega\times E^{*}_{\sigma^{*}}$. 
\end{lem}
Proof of Lemma \ref{lem 5.3}. Due to Proposition \ref{prop 31} the integrands $f_{n}^{*}$ are normal integrands.  From  Corollary \ref{coro7.8} and Definition \ref{def6.201} we have    $\displaystyle\limsup_{n} f_{n}^{*}=f^{*}$. From  Definition \ref{Def 117} $f^{*}$ is a normal integrand on $\Omega\times E^{*}_{\sigma^{*}}$. $\;\Box$ 
\begin{lem} \label{lem 5.4} Let $(E, \Vert\;.\Vert)$ be a separable reflexive Banach space and let $(f_{n})_n$ be a sequence of $\mathbb{T}\otimes \mathcal{B}(E)$ measurable integrands defined on $\Omega\times E$ and bounded below by some measurable function $h: \Omega\to {\R}$. If $f=seq\;\sigma-li_{e}f_{n}$, then $f^{*}$ is a normal integrand on $\Omega\times E^{*}_{\sigma^{*}}$.
\end{lem}
Proof of Lemma \ref{lem 5.4}. The non negative integrands $g_{n}=f_{n}-h$ are $\mathbb{T}\otimes \mathcal{B}(E)$ measurable. For each integer $k$ define $\displaystyle g_{n,\;k}=g_{n}+k^{-1}\Vert .\Vert^{2}$ and  $g_{k}=seq\;\sigma-li_{e}g_{n\;k}$. The integrands $\displaystyle g_{n,\;k}$ are $\mathbb{T}\otimes \mathcal{B}(E)$ measurable, bounded below by $k^{-1}\Vert .\Vert^{2}$, thus from Lemma \ref{lem 5.3}, for each integer $k$, $g_{k}^{*}$ is normal on $\Omega\times E^{*}_{\sigma^{*}}$. Moreover   $$f-h=\inf_{k}seq\;\sigma-li_{e}g_{n,\;k}=\inf_{k} g_{k}\;,$$ therefore 
$\displaystyle f^{*}+h=\sup_{k} g_{k}^{*}$. The  integrands $g_{k}^{*}$ being normal,  Definition \ref {Def 117} shows that the integrand $\displaystyle f^{*}=\sup_{k} g_{k}^{*}-h$ is a normal integrand on  $\Omega\times E^{*}_{\sigma^{*}}$. The proof of Lemma \ref{lem 5.4} is complete. $\;\Box$\\\\
End of the proof of Theorem \ref{thm 5.1}. The first assertion is exactly  Lemma \ref{lem 5.3}. Suppose now $E$ is reflexive and separable. For each integer $k$ define $f_{n,\;k}=\sup(f_{n}\;\;,\;-k)$.  The $f_{n,\;k}$ are measurable and bounded below by $k$. If $\displaystyle f_{k}= seq\;\sigma-li_{e}f_{n,\;k}$ applying  Lemma \ref{lem 5.4} we deduce that each $\displaystyle {f_{k}}^{*}$ is a normal integrand on  $\Omega\times E^{*}_{\sigma^{*}}$. But $f=\inf_{k} f_{k}$, then the integrand $\displaystyle f^{*}=\sup_{k}{f_{k}}^{*}$ is a normal integrand on $\Omega\times E^{*}_{\sigma^{*}}$. The proof of Theorem \ref{thm 5.1} is complete.$\;\Box$\\\\
Assume as precedently that the Banach $E$ is separable.  When $E^{*}$ is equipped  with the tribe $\mathcal{B}(E_{\sigma^{*}}^{*})$, it is known that a function $x^{*}:\Omega\to E^{*}$ is measurable if and only if it is scalarly measurable, that is for every $e\in E$, the map $\omega\mapsto\langle x^{*}(\omega), e\rangle$ is $\mathbb{T}$-measurable. Indeed as a function $x^{*}:\Omega\to E^{*}$ is scalarly measurable if and only if the inverse image of each weak* open translated half-space of $E^{*}$ is  $\mathbb{T}$-measurable, it suffices to prove that the tribe $\mathcal{B}^{'}$ generated by the weak* open translated half-spaces coincide with  $\mathcal{B}(E_{\sigma^{*}}^{*})$.
Since $E$ is  strongly separable the dual norm $\Vert x^{*}\Vert_{*}$ of any  scalarly measurable  $E^{*}$-valued function is $\mathbb{T}$-measurable. Let  $L_{p}(\Omega,E_{\sigma^{*}}^{*})$ $1\leq p\leq \infty$ be the set of equivalence classes for the equality almost everywhere of scalarly  $E^{*}$-valued functions $x^{*}$ such   $\Vert x^{*}\Vert_{*}\in L_{p}(\Omega,{\R})$.
Then it is known (see \cite{706} for $p=1$, and \cite{78} Theorem 2.112 ), that the strong dual of $L_{p}(\Omega,E)$ is $L_{q}(\Omega,E_{\sigma^{*}}^{*})$ with the duality pairing
$$\langle x^{*}, x\rangle=\int_{\Omega} \langle x^{*}(\omega), x(\omega)\rangle d\mu(\omega)\;.$$ 
The following useful duality result is an immediate consequence of the main Theorem  in \cite{890} and \cite{4} Theorem VII-7 (which is valid when the decomposable spaces considered are subspaces of scalarly measurable functions instead of subspaces of scalarly integrable functions). 
\begin{thm} \label{thm 483} (see \cite{890} main Theorem, \cite{4} Theorem VII-7) Let $E$ be a separable Banach space and $g:\Omega\times E^{*}\to \overline{{\R}}$ be a normal integrand on $\Omega \times E_{\sigma^{*}}^{*}$. If $I_g$ is finite at at least one point of $L_{\infty}(\Omega,E_{\sigma^{*}}^{*})$  then for the duality between  $L_{\infty}(\Omega,E_{\sigma^{*}}^{*})$ and $L_{1}(\Omega,E)$,
for every $x\in L_{1}(\Omega,E)$ we have $I_{g}^{*}(x)=I_{g^{*}}(x)$.
\end{thm}

\section{Two properties of the Mackey topology $\tau(L_{\infty}(\Omega,E_{\sigma^{*}}^{*}), L_{1}(\Omega,E))$}

In this section we suppose that the Banach $E$ is separable. As above,  $\mathcal{B}(E^{*}_{\sigma^{*}})$  is the Borel tribe of $(E^{*},\sigma(E^{*}, E))$. 
 Recall that the Mackey topology $\tau^{*}=\tau(L_{\infty}(\Omega,E_{\sigma^{*}}^{*}), L_{1}(\Omega,E))$ is the topology of uniform convergence on the $\sigma (L_{1}(\Omega,E), L_{\infty}(\Omega,E_{\sigma^{*}}^{*}))$-compacts convex sets of $L_{1}(\Omega,E)$.
\begin{prop} \label{prop 704} The $\sigma (L_{1}(\Omega,E), L_{\infty}(\Omega,E_{\sigma^{*}}^{*}))$-compact sets are uniformly integrable.
\end{prop}
Proof of Proposition \ref{prop 704}. Since $L_{1}(\Omega,E)^{*}=L_{\infty}(\Omega,E_{\sigma^{*}}^{*})$, due to Eberlein-Smulian Theorem, each  $\sigma (L_{1}(\Omega,E), L_{\infty}(\Omega,E_{\sigma^{*}}^{*}))$-compact set is a $\sigma (L_{1}(\Omega,E), L_{\infty}(\Omega,E_{\sigma^{*}}^{*}))$-sequentially compact set therefore it is a $\sigma (L_{1}(\Omega,E), L_{\infty}(\Omega,E_{\sigma^{*}}^{*}))$-conditionally compact set (every sequence has a weakly Cauchy subsequence), and we use \cite{150} IV Theorem 4. Other proof: see also \cite{703}.$\;\Box$ 
\begin{coro} \label{coro 73} If $E$ is a separable Banach space, for every bounded sequence $(x_{n}^{*})$ in $L_{\infty}(\Omega,E_{\sigma^{*}}^{*})$, for every decreasing sequence $(A_{n})_n$ of measurable sets with a negligible intersection, the sequence  $(x_{n}^{*}1_{A_{n}})_n$ $\tau^{*}$-converges to the origin.
\end{coro}
Proof of Corollary \ref{coro 73}. The topology  $\tau^{*}$ being the topology of uniform convergence on the convex $\sigma (L_{1}(\Omega,E), L_{\infty}(\Omega,E_{\sigma^{*}}^{*}))$-compact sets,  Proposition  \ref{prop 704} gives the result. Indeed, considering a  $\sigma (L_{1}(\Omega,E), L_{\infty}(\Omega,E_{\sigma^{*}}^{*}))$-compact set $K$, it suffices to use Theorem \ref{prop991} $(a)$ and to remark that:
$$\sup_{x\in K}\vert\int_{A_{n}}\langle x, x_{n}^{*}\rangle d\mu\vert \leq (\sup_{n}\Vert x_{n}^{*}\Vert_{\infty})\sup_{x\in K}\Vert x1_{A_{n}}\Vert_{1}\;.\;\;\Box $$

\begin{prop}\label{lem301} Let $\tau^{*}$ be the Mackey topology $\tau(L_{\infty}(\Omega,E_{\sigma^{*}}^{*}), L_{1}(\Omega,{E}))$. Given a sequence $(f_{n})_n$ of extended real-valued $\mathbb{T}\otimes \mathcal{B}_{\sigma^{*}}(E^{*})$ measurable integrands defined on $\Omega\times E^{*}$ such that there exist an element  $x_{0}^{*}\in L_{\infty}(\Omega,E_{\sigma^{*}}^{*})$, an element $u_{0}\in L_{1}(\Omega, {\R})$ verifying for every integer $n$ $f_{n}(x_{0}^{*})\leq u_{0}$, then with $f=\Vert .\Vert_{*}-ls_{e} f_{n}$ one has
$$\tau^{*}- ls_{e} I_{f_{n}}\leq I_{f}.$$
\end{prop} 
Proof of Proposition \ref{lem301}. Clearly considering $g_{n}(\omega, e)=f_{n}(\omega,x_{0}^{*}+e)-u_{0}(\omega)$, one can suppose that $x_{0}^{*}=0$ and $u_{0}=0$. Let us endows $E^{*}$ with the dual norm $\Vert. \Vert_{*}$ and $E^{*}\times {\R}$ with the product norm and the tribe $\mathcal{B}(E_{\sigma^{*}}^{*}\times{\R})$ . We denote the unit ball of  $E^{*}$ by $B^{*}$ and the unit ball of $L_{\infty}(\Omega,E_{\sigma^{*}}^{*})$ by $B_{\infty}$. Let $x^{*}\in L_{\infty}(\Omega,E_{\sigma^{*}}^{*})$, $r\in{\R}$ verifying $I_{f}(x^{*})<r<\infty$. Then there exists $u\in epi f$, integrable such $f(x^{*})\leq u$ and $\int u d\mu<r$. Due to Definition \ref{def2.11} for every $\omega\in \Omega$, $epi f_{\omega}=\liminf_{n} epi{f_{n}}_{\omega}$, hence the function $s_{n}(\omega)=d((x^{*}(\omega),u(\omega)), epi{f_{n}}_{\omega})$ converges simply to $0$. 
\begin{lem} \label{lem 149} The integrand $g(\omega,e^{*},r)=\Vert x^{*}(\omega)-e^{*}\Vert_{*}+\vert u(\omega)-r\vert$ is $\mathbb{T}\otimes\mathcal{B}(E_{\sigma^{*}}^{*}\times{\R})$ measurable. 
\end{lem}
Proof of Lemma \ref{lem 149}. Let $x^{*}$ be a scalarly measurable function.
 Since the map $(\omega,e^{*})\mapsto e^{*}$ is  $\mathbb{T}\otimes\mathcal{B}(E_{\sigma^{*}}^{*})$-measurable and the map $(\omega,e^{*})\mapsto x^{*}(\omega)$ is measurable too, the map $(\omega,e^{*})\mapsto x^{*}(\omega)-e^{*}$ is $\mathbb{T}\otimes\mathcal{B}(E_{\sigma^{*}}^{*})$-measurable. The dual norm  $\Vert. \Vert_{*}$ being $\sigma(E^{*}, E)$-lower semicontinuous, the map $(\omega,e^{*})\mapsto \Vert x^{*}(\omega)-e^{*}\Vert_{*}$
is $\mathbb{T}\otimes\mathcal{B}(E_{\sigma^{*}}^{*})$-measurable and the integrand  $$h:(\omega,e^{*},r)\mapsto \Vert x^{*}(\omega)-e^{*}\Vert_{*}$$
is  $\mathbb{T}\otimes\mathcal{B}(E_{\sigma^{*}}^{*}\times{\R})$-measurable. The Caratheodory integrand $(\omega,r)\mapsto \vert u(\omega)-r\vert$ being $\mathbb{T}\otimes\mathcal{B}({\R})$-measurable, the integrand  $$k:(\omega,e^{*},r)\mapsto \vert u(\omega)-r\vert$$ is $\mathbb{T}\otimes\mathcal{B}(E_{\sigma^{*}}^{*}\times{\R})$-measurable. Then $g=h+k$ is the sum of $\mathbb{T}\otimes\mathcal{B}(E_{\sigma^{*}}^{*}\times{\R})$ measurable integrands hence it is $\mathbb{T}\otimes\mathcal{B}(E_{\sigma^{*}}^{*}\times{\R})$ measurable. The proof of Lemma \ref{lem 149} is complete.\\ 
Since $f_n$ is $\mathbb{T}\otimes\mathcal{B}(E_{\sigma^{*}}^{*})$-measurable, the multifunction $epi f_{n}$ has a  $\mathbb{T}\otimes\mathcal{B}(E_{\sigma^{*}}^{*}\times{\R})$-measurable graph and the above lemma ensures that
$g$ is $\mathbb{T}\otimes\mathcal{B}(E_{\sigma^{*}}^{*}\times{\R})$-measurable. 
Therefore, due to \cite{4} Lemma III 39, the function $s_n$ is  $\mathbb{T}$-measurable. Let $\alpha$ be a function in $L_{\infty}(\Omega,{\R})\cap L_{1}(\Omega,{\R})$ with positive values such $\int \alpha d\mu=1$ and let $\epsilon=r-\int u d\mu$. Since every $\sigma(L_{1}(\Omega,E), L_{\infty}(\Omega,E_{\sigma^{*}}^{*}))$-compact set is bounded, the strong topology of $ L_{\infty}(\Omega,E_{\sigma^{*}}^{*})$ is stronger or equal that  $\tau^{*}$. Let $V$ and $W$ be symmetric $\tau^{*}$-neighbourhoods of $0$ such $W+W\subset V$, and $0<\beta\leq \frac{\epsilon}{2}$   such $2\beta\alpha B_{\infty}\subset W$.\\
Define $\Omega_{n}=\cap_{k\geq n}\{s_{k}< \beta\alpha\}$. Then $(\Omega_{n})_n$ is an increasing covering of $\Omega$ by measurable sets. Due to Corollary \ref{coro 73}, the sequence $(x^{*}1_{\Omega_{n}^{c}})_n$ $\tau^*$-converges to $0$. Thus there exists an integer $n_{W}$ such $n\geq n_{W}\Rightarrow x^{*}1_{\Omega_{n}^{c}}\in W$. 
Since $g$ is $\mathbb{T}\times\mathcal{B}(E_{\sigma^{*}}^{*}\times{\R})$-measurable, the multifunction $$\Gamma_{n}(\omega)=\{(\omega,e^{*},r): g(\omega, e^{*},r)< s_{n}(\omega)+\frac{1}{n}\beta \alpha(\omega)\}$$ 
has a $\mathbb{T}\times\mathcal{B}(E_{\sigma^{*}}^{*}\times{\R})$ measurable graph. On each $\Omega_{n}$, $\Gamma_{n} \cap epi f_{n}$ has a $\mathbb{T}\times\mathcal{B}(E_{\sigma^{*}}^{*}\times{\R})$ measurable graph and nonempty values.  
Due to \cite{4} Theorem III 22, on each $\Omega_{n}$, there exists a scalarly-measurable selection $(x^{*}_{n}, u_{n})$ of $epif_{n}$ satisfying 
 $$\quad\quad\quad\quad\quad\quad\quad\quad\quad\quad\quad\quad\quad\quad\quad \Vert x^{*}-x^{*}_{n}\Vert_{*}+\vert u-u_{n}\vert\leq 2\beta \alpha\;.\;\;\quad\quad\quad\quad\quad \quad\quad\quad\quad \quad\quad    (2)$$ 
 Define $\overline{x^{*}_{n}}(\omega)=x^{*}_{n}(\omega)$ if $\omega\in \Omega_{n}$, $0$ otherwise; and $\overline{u_{n}}(\omega)=u_{n}(\omega)$ if $\omega\in \Omega_{n}$, $0$ otherwise. 
Then due to $(2)$: $$n\geq n_{W}\Rightarrow x^{*}-\overline{x^{*}_{n}}=(x^{*}-\overline{x^{*}_{n}})1_{\Omega_{n}}+x^{*}1_{\Omega_{n}^{c}}\in 2\beta\alpha B_{\infty}+W\subset W+W\subset V.\quad\quad\quad\quad \quad\quad\quad\quad  (3)$$
Since $f_{n}(0)\leq 0$ we get $f_{n}(\overline{x_{n}})\leq \overline{u_{n}}$ and using $(2)$:
$$I_{f_{n}}(\overline{x^{*}_{n}})\leq \int_{\Omega_{n}} u_{n} d\mu\leq \int_{\Omega_{n}} (u +2\beta \alpha) d\mu\leq  \int_{\Omega_{n}} (u +\epsilon\alpha) d\mu.$$
Hence with $(3)$
$$\limsup_{n}\inf_{y^{*}\in x^{*}+V} I_{f_{n}}(y^{*})\leq  \int_{\Omega} (u +\epsilon\alpha) d\mu=r\;,$$
and since $V$ is arbitrary:
$$\tau^{*}-ls_{e}I_{f_{n}}(x^{*})\leq r.$$
Therefore for every $x^*$,
$$\tau^{*}- ls_{e} I_{f_{n}}(x^{*})\leq I_{f}(x^{*}).\;\;\Box $$

\section{Convergence results}
In this section $E$ is a separable Banach space and we consider now a topological subspace $(\mathcal{X}, \mathcal{T})$ of $L_{0}(\Omega,E)$, satisfying at $x\in \mathcal{X}$, as the usual topologies on the spaces  $L_{p},\;1\leq p\leq\infty$, the following property:\\\\
$(P)$ {\it Every sequence $\mathcal{T}$-converging to  $x\in \mathcal{X}$ converges to $x$ in the Biting sense }. 
\begin{Def} Let $(I, \leq)$ be a totally ordered set. We will say that  a family $(f_{i})_{i\in I}$ of functions  satisfies eventually a property ($\mathcal{P}$) if there exists $i_{0}\in I$ such that the family $(f_{i})_{i\geq i_{0}}$ satisfies this property.
\end{Def}
\begin{Def} \label{def6.2} Given a sequence ${(f_{n})}_n$ of $\overline{{\R}}$-valued integrands  defined on $\Omega\times E$, we will say that a  sequence $(x_{n})_n$ satisfies the boundedness property (with respect to ${(f_{n})}_n$) if the sequence ${(f_{n}^{-}(x_{n}))}_{n}$ is eventually bounded in $L_{1}(\Omega,{\R})$. The sequence $(x_{n})_n$ is said to satisfy the lower compactness property (with respect to ${(f_{n})}_n$) if the sequence ${(f_{n}^{-}(x_{n}))}_{n}$ is eventually uniformly integrable in $L_{1}(\Omega,{\R})$. 
\end{Def}
\begin{thm} \label{thm6.01} Let $E$ be a separable Banach space and let  ${(f_{n})}_n$ be a sequence of $\mathbb{T}\otimes \mathcal{B}(E)$ measurable and  $\overline{{\R}}$-valued integrands.  Suppose that property $(P)$ holds.  Given $x\in \mathcal{X}$, a Nagumo tight sequence $(x_{n})_n$  $\mathcal{T}$-converging to $x$ and satisfying the boundedness property, then  for $f=seq\;\sigma-li_{e}f_{n}$ one has:
$$ \liminf_{n}I_{f_{n}}(x_{n})\geq I_{f^{**}}(x)-\delta^{+}((-f_{n}(x_{n}))_{n})\;.$$
\end{thm}
\begin{coro} \label{thm6.0}  Let $E$ be a separable Banach space and let ${(f_{n})}_n$ be a sequence of $\mathbb{T}\otimes \mathcal{B}(E)$ measurable and  $\overline{{\R}}$-valued integrands.      Suppose that property $(P)$ holds.  Given $x\in \mathcal{X}$, a Nagumo tight sequence $(x_{n})_n$  $\mathcal{T}$-converging to $x$ and satisfying the lower compactness property, then  with $f=seq\;\sigma-li_{e}f_{n}$ one has
$$ \liminf_{n}I_{f_{n}}(x_{n})\geq I_{f^{**}}(x).$$  
\end{coro}
When $E$ is is reflexive and separable, the assumptions on tightness are satisfied (Proposition \ref{prop00}) and we obtain:
\begin{coro} \label{thm6.1} Suppose $E$ is reflexive and separable, $x\in \mathcal{X}$, and property $(P)$ holds.  Given a sequence ${(f_{n})}_n$ of $\mathbb{T}\otimes \mathcal{B}(E)$ measurable integrands, 
then for every sequence  $(x_{n})_n$ $\mathcal{T}$-converging to $x\in \mathcal{X}$ and satisfying the boundedness property, the following inequality holds  with $f=seq\;\sigma-li_{e}f_{n}$, 
$$ \liminf_{n}I_{f_{n}}(x_{n})\geq I_{f^{**}}(x)-\delta^{+}((-f_{n}(x_{n})).$$
If in addition $(x_{n})_n$ satisfies the lower compactness property then,
$$ \liminf_{n}I_{f_{n}}(x_{n})\geq I_{f^{**}}(x).$$
\end{coro}
The following definition, when the sequence ${(f_{n})}_n$ is constant,  is exactly the criterion put in light by A. D. Ioffe in the study of the strong weak-semicontinuity \cite{70}. 
\begin{Def} \label{def6.3} Given a sequence ${(f_{n})}_n$ of $\overline{{\R}}$-valued integrands defined on $\Omega\times E$, we will say that it satisfies the 
$\mathcal{T}$-Ioffe's criterion at $x\in \mathcal{X}$ if for every subsequence ${(f_{n_{k}})}_k$ and every sequence $(x_{k})_k$ $\mathcal{T}$-converging to $x$ such the sequence $(I_{f_{n_{k}}}(x_{k}))_k$ is bounded above, the sequence  $(x_{k})_k$  has the lower compactness property with respect to ${(f_{n_{k}})}_k$. 
\end{Def}
\begin{coro}\label{thm6.1111} Suppose $E$ is reflexive and separable, $x\in \mathcal{X}$, and property $(P)$ holds.  Given a sequence ${(f_{n})}_n$ of measurable integrands satisfying 
the $\mathcal{T}$-Ioffe's criterion at $x\in \mathcal{X}$, then:
$$ seq\;\mathcal{T}-li_{e}I_{f_{n}}(x)\geq I_{f^{**}}(x).$$ 
\end{coro}
There exists a converse of Corollary \ref{thm6.1111}.
Recall that a  subset $X$ of $L_{0}(\Omega,E)$ is said to be decomposable, see \cite{8}, if given  elements $x$, $y$ of $X$, then for every measurable set $A$, $z=x1_{A}+y1_{A^{c}}$ is also an element of $X$. We will suppose that the topological space $(\mathcal{X}, \mathcal{T})$ satisfies the following property:\\\\
{\it$ (Q)$ Given a sequence $(A_{n})_n$ of measurable subsets such for every measurable set $A$ of finite measure $\displaystyle \lim_{n} \mu(A\cap A_{n})=0$, and two sequences $(x_{n})_n$ and $(y_{n})_n$ $\mathcal{T}$-converging to $x\in \mathcal{X}$, then the sequence $(z_{n})_n$, with $z_{n}=x_{n}1_{A_{n}}+y_{n}1_{A^{c}_{n}}$, is $\mathcal{T}$-converging to $x$.}
\begin{thm} \label{thm6.2} (see \cite{67} Theorem 4.2) Suppose the measure $\mu$ is atomless, the topological space $(\mathcal{X}, \mathcal{T})$ is decomposable and satisfies 
 property $ (Q)$ at some point $x\in \mathcal{X}$. Given a sequence ${(f_{n})}_n$ of $\overline{{\R}}$-valued measurable integrands defined on $\Omega\times E$ and a $\mathbb{T}\otimes \mathcal{B}(E)$-measurable integrand $\overline{{\R}}$-valued $f$ satisfying:\\
for every $\epsilon>0$ there exists a sequence $(y_{n})_n$ $\mathcal{T}$-converging to $x$ such  $\delta^{+}((f_{n}^{-}(y_{n}))_n)\leq\epsilon$ and  $\;\displaystyle \limsup_{n} I_{f_{n}}(y_{n})\leq I_{f}(x)+\epsilon$.\\
Then if $\;-\infty<I_{f}(x)$, the Ioffe's criterion \ref{def6.3} is necessary to get  the inequality
$$ I_{f}(x)\leq seq\;\mathcal{T}-li_{e}I_{f_{n}}(x)\;. $$
\end{thm}
Notice the following simple statements of the Ioffe's criterion:
\begin{lem}(see \cite{67} Lemma 4.1) Suppose that there exists a $\mathcal{T}$-converging sequence $(y_{n})_n$ to $x\in X$ such the sequence $(I_{f_{n}}(y_{n}))_n$ is bounded above. \\
$(a)$ the Ioffe's criterion is equivalent to the following:\\
"for every sequence $(x_{n})_n$ $\mathcal{T}$-converging to $x$ such that the sequence $(I_{f_{n}}(x_{n}))_n$ is bounded above, the sequence  $(x_{n})_n$  has the lower compactness property with respect to ${(f_{n})}_n$''.\\
$(b)$ Moreover if the measure is atomless and if the topological space $\mathcal{X}$ satisfies the property $(Q)$, then the Ioffe's criterion is equivalent to\\
" every sequence $(x_{n})_n$ $\mathcal{T}$-converging to $x$ has the lower compactness property with respect to ${(f_{n})}_n$''.
\end{lem}
Proof of Theorem \ref{thm6.01}. We begin with a sequence of preliminary propositions.
\begin{prop} \label{lem6.8}  Given a sequence ${(f_{n})}_n$ of measurable integrands bounded below by a Nagumo integrand, we have on $\mathcal{X}=L_{1}(\Omega,E)$, when $\mathcal{T}$ is the weak topology and $f=seq\;\sigma-li_{e}f_{n}\;$:  
$$ seq\;\mathcal{T}-li_{e}I_{f_{n}}\geq I_{f^{**}}.$$ 
\end{prop}
Proof of Proposition \ref{lem6.8}. 
Let $g=\displaystyle\Vert.\Vert_{*}-ls_{e} {f^{*}_{n}}$ be the integrand pointwise upper epi-limit of the ${f^{*}_{n}}$ for the topology of the dual norm $\Vert. \Vert_{*}$ defined by $g(\omega,e)=\displaystyle\Vert.\Vert_{*}-ls_{e} {f^{*}_{n}}_{\omega}(e)$. The sequence ${(f_{n})}_n$ being  bounded below by a Nagumo integrand and the topology of the dual norm on $E^{*}$ being stronger than or equal to the Mackey topology $\tau(E^{*}, E)$, then from Corollary  \ref{coro7.8}, we have for (almost) every $\omega$, $g_{\omega}=\limsup_{n} {f^{*}_{n_{\omega}}}=f_{\omega}^{*}$ therefore $g=f^{*}$ and from Theorem \ref{thm 5.1}, $g$ is a normal integrand on $\Omega\times E^{*}_{\sigma^{*}}$.
Since the $f_n$'s are non negative, then for every integer $n$, $f_{n}^{*}(0)\leq 0$. From  Proposition \ref{lem301} with the Mackey topology $\tau^{*}=\tau(L_{\infty}(\Omega,E^{*}_{{\sigma}^{*}}), L_{1}(\Omega,E))$ (in case $E$ is a Banach space with strongly separable dual, from \cite{106} Corollary 4.8), 
$$\displaystyle \tau^{*}-ls_{e}I_{f_{n}^{*}}\leq I_{g}\;.$$
Thus with Proposition \ref{prop6.1}, 
$$(seq\;\sigma-li_{e}I_{f_{n}})^{*}\leq \displaystyle \tau^{*}-ls_{e}I_{f_{n}^{*}}\leq I_{g}\;.$$
As a consequence:
$$(seq\;\sigma-li_{e} I_{f_{n}})^{**}\geq (I_{g})^{*}\;.$$
and since $g(0)=\displaystyle\limsup {f_{n}}^{*}(0)\leq 0$, with \cite{4} Theorem VII-7 (see Theorem \ref{thm 483}) we  get on $L_{1}(\Omega,E)$:
$$(I_{g})^{*}=I_{g^{*}}$$ and therefore:
$$(seq\;\sigma-li_{e} I_{f_{n}})^{**}\geq I_{g^{*}}.$$
But  $g^{*}=f^{**}$, and as a consequence:
$$seq\;\sigma-li_{e} I_{f_{n}}\geq(seq\;\sigma-li_{e} I_{f_{n}})^{**}\geq I_{f^{**}}.$$
This ends the proof of Proposition \ref{lem6.8}.$\;\Box$ 
\begin{prop} \label{lem6.88} Suppose ${(f_{n})}_n$ is a sequence of non negative measurable integrands.  A Nagumo tight sequence $(x_{n})_n$ of elements of $\mathcal{X}=L_{1}(\Omega,E)$ which weakly-converges to $x$ satisfies
$$ \liminf_{n}I_{f_{n}}(x_{n})\geq I_{f^{**}}(x).$$ 
\end{prop}
Proof of Proposition \ref{lem6.88}.
There exists a Nagumo integrand $h$ such: $\displaystyle\sup_{n} I_{h}(x_{n})\leq 1$. Let us consider now for every $\epsilon>0$ the integrand $g^{\epsilon}_{n}=f_{n}+\epsilon h$. Clearly, $f_n$ being non negative, each $g^{\epsilon}_{n}$ is bounded below by $\epsilon h$. 
Therefore from Proposition \ref{lem6.8},  since we have with $g^{\epsilon}=seq\;\sigma-li_{e} g^{\epsilon}_{n}\geq seq\;\sigma-li_{e}f_{n}=f$, we get
$$\liminf_{n} I_{g^{\epsilon}_{n}}(x_{n})\geq seq\;\sigma-li_{e}I_{g^{\epsilon}_{n}}(x)\geq I_{{g^{\epsilon}}^{**}}(x)\geq I_{f^{**}}(x).$$
But since $ I_{f_{n}}(x_{n})\geq  I_{g^{\epsilon}_{n}}(x_{n})-\epsilon$, we deduce:
$$\displaystyle\liminf_{n} I_{f_{n}}(x_{n})\geq\liminf_{n} I_{g^{\epsilon}_{n}}(x)-\epsilon\geq  I_{f^{**}}(x)-\epsilon.$$
The above inequality being valid for every $x\in L_{1}(\Omega,E)$, for every tight sequence ${(x_{n})}_n$ weakly converging to $x$, and for every $\epsilon>0$, the proof of Proposition \ref{lem6.88} is complete.$\;\Box$ 
\begin{prop} \label{lem6.9} For any  measurable set $K$ of finite measure the conclusion of Corollary \ref{thm6.0} holds when  $\mathcal{X}=L_{1}(K,E)$ and $\mathcal{T}$ is the weak topology of $\mathcal{X}$.
\end{prop}
Proof of Proposition \ref{lem6.9}. Let a Nagumo tight sequence $(x_{n})_n$ weakly-converging to $x\in L_{1}(K,E)$ and satisfying the lower compactness property and $$\displaystyle\liminf_{n} I_{f_{n}}(x_{n})<r<\infty\;.$$
 Due to the lower compactness property of $(x_{n})_n$, the sequence
${(f_{n}^{-}(x_{n}))}_{n\geq m}$ is, for $m$ large enough, uniformly integrable.\\
 Setting  $A_{n,m}=\{f_{n}^{-}(x_{n})\geq m\}$, and $\displaystyle r_{m}=\sup_{n\geq m}\int_{A_{n,m}} f_{n}^{-}(x_{n})d\mu$, we have
$\displaystyle\lim_{m\to \infty}r_{m}=0.$\\
For every integer $m$, define $f_{n}^{m}=\sup(f_{n},-m)$, then
$$f_{n}(x_{n})=f_{n}(x_{n})1_{A^{c}_{n,m}} +f_{n}(x_{n})1_{A_{n,m}}= f_{n}^{m}(x_{n})1_{A^{c}_{n,m}} +f_{n}(x_{n})1_{A_{n,m}}\;.$$
Since $\displaystyle f_{n}^{m}(x_{n})1_{A_{n,m}}=-m1_{A_{n,m}}$, $\displaystyle f_{n}^{m}(x_{n})1_{A^{c}_{n,m}}\geq f_{n}^{m}(x_{n})1_{A^{c}_{n,m}}+f_{n}^{m}(x_{n})1_{A_{n,m}}=f_{n}^{m}(x_{n})$, one has
$$f_{n}(x_{n})\geq  f_{n}^{m}(x_{n})-f_{n}^{-}(x_{n})1_{A_{n,m}}\;.$$
Therefore for $n\geq m$, 
$$I_{ f_{n}}(x_{n})\geq I_{ f_{n}^{m}}(x_{n}) -r_{m}\;.$$ 
Using the above inequality and Proposition \ref{lem6.88} with $\displaystyle f^{m}=seq\;\sigma-li_{e} f_{n}^{m}$ we deduce
$$r>\displaystyle \liminf_{n}I_{ f_{n}}(x_{n})\geq \displaystyle \liminf_{n} I_{ f_{n}^{m}}(x_{n}) -r_{m}\geq I_{{f^{m}}^{**}}(x) -r_{m}\geq I_{f^{**}}(x) -r_{m} ,$$
and since $\displaystyle \lim_{m\to \infty} r_{m}=0$, we obtain the desired inequality: $I_{f^{**}}(x)\leq r$.
The proof of Proposition \ref{lem6.9} is complete.$\;\Box$ \\\\
End of the proof of Theorem \ref{thm6.01}.\\
Let  a Nagumo tight sequence  ${(x_{n})}_n$ $\mathcal{T}$-converging to $x\in \mathcal{X}$ and satisfying the boundedness property. Suppose $\displaystyle\liminf_{n} I_{f_{n}}(x_{n})<r<\infty$. We can extract a subsequence ${(x_{n_{k}})}_k$ such 
$$\displaystyle\limsup_{k} I_{f_{n_{k}}}(x_{n_{k}})<r<\infty\;.$$ 
From $(P)$ the sequence $(x_{n_{k}})_k$ is Biting converging to $x$. The sequence $(f^{-}_{n_{k}}(x_{n_{k}}))_k$ being eventually bounded in $L_{1}(\Omega,{\R})$, using the other form of the Biting Lemma, Corollary \ref{coro 658} in the $\sigma$-finite case with $E={\R}$, we can suppose (up to an extraction of a subsequence) that the sequences $(f^{-}_{n_{k}}(x_{n_{k}}))_k$ and $(x_{n_{k}})_k$ are converging in the Biting sense.
Let ${(\Omega_{p})}_p$ be a common increasing covering of $\Omega$ by measurable sets  appearing in 
Definition \ref{def6.0}. Since $\mu$ is $\sigma$-finite, we may suppose that the $\Omega_{p}$'s are of finite measure.
By the Dunford-Pettis Theorem \cite{78} Theorem 2.54 (see also \cite{150}, \cite{98}) the sequence of restrictions $(f^{-}_{n_{k}}(x_{n_{k}})\vert\Omega_{p})_k$ is uniformly integrable in the sense of Definition \ref{def33.2} in $L_{1}(\Omega_{p}, E)$. This proves that for each integer $p$ the sequence of restrictions $(x_{n_{k}}\vert \Omega_{p})_k$ has the lower compactness property with respect to the sequence $(f_{n_{k}}\vert\Omega_{p})_k$ of the restrictions to $\Omega_{p}$ of the integrands $f_{n_{k}}$.\\
The sequence ${(f_{n}^{-}(x_{n}))}_n$ being eventually bounded in $L_{1}(\Omega,{\R})$,  we deduce that the sequence $(I_{f_{n_{k}}^{+}}(x_{n_{k}}))_k$ is eventually bounded above by some positive real number $s$. Clearly $f\leq g=seq\;\sigma-li_{e} f_{n_{k}}^{+}$ and from Proposition \ref{lem6.88} applied to the sequence $(f_{n_{k}}^{+})_k$, for every integer $p$ on the space $L_{1}(\Omega_{p}, E)$,
$$\displaystyle \int^{*}_{\Omega_{p}}{f^{**}(x)}^{+} d\mu\leq \int^{*}_{\Omega_{p}}g^{**}(x) d\mu\leq \liminf_{k} I_{f_{n_{k}}^{+}}(x_{n_{k}})\leq s<\infty\;.$$
The Monotone Convergence Theorem for the upper integral gives:
$$\displaystyle \int^{*}_{\Omega}{f^{**}(x)}^{+} d\mu\leq s<\infty\;.$$
For each integer $p$ we have for $k\geq p$
$$I_{f_{n_{k}}}(x_{n_{k}})=\int^{*}_{\Omega_{p}}f_{n_{k}}(x_{n_{k}})d\mu + \int^{*}_{\Omega^{c}_{p}}f_{n_{k}}(x_{n_{k}})d\mu\;,$$
thus
$$I_{f_{n_{k}}}(x_{n_{k}})\geq\int^{*}_{\Omega_{p}}f_{n_{k}}(x_{n_{k}})d\mu + \inf_{k\geq p}\int^{*}_{\Omega^{c}_{p}}f_{n_{k}}(x_{n_{k}})d\mu\;,$$
hence
$$\limsup_{k} I_{f_{n_{k}}}(x_{n_{k}})\geq \liminf_{k}\int^{*}_{\Omega_{p}}f_{n_{k}}(x_{n_{k}})d\mu+\inf_{k\geq p}\int^{*}_{\Omega^{c}_{p}}f_{n_{k}}(x_{n_{k}})d\mu\;.$$
Therefore the use of Proposition \ref{lem6.9} gives
$$r>\limsup_{k} I_{f_{n_{k}}}(x_{n_{k}})\geq \int^{*}_{\Omega_{p}} f^{**}(x)d\mu+\inf_{k\geq p}\int^{*}_{\Omega^{c}_{p}}f_{n_{k}}(x_{n_{k}})d\mu\;,$$
and we obtain
$$r> I_{f^{**}}(x)-\int^{*}_{\Omega^{c}_{p}}{f^{**}(x)}^{+}d\mu+\inf_{k\geq p}\int^{*}_{\Omega^{c}_{p}}f_{n_{k}}(x_{n_{k}})d\mu\;.$$
Since $I_{{f^{**}(x)}^{+}}<\infty$ we have $\displaystyle\lim_{p}\int^{*}_{\Omega^{c}_{p}}{f^{**}(x)}^{+}d\mu=0$, and we deduce:
$$r\geq I_{f^{**}}(x)+\liminf_{p}\inf_{k\geq p}\int^{*}_{\Omega^{c}_{p}}f_{n_{k}}(x_{n_{k}})d\mu \;.$$
Therefore:
$$r\geq I_{f^{**}}(x)+\displaystyle\inf_{\sigma\in \Sigma, \;\sigma=(S_{p})_{p}}\liminf_{p}\inf_{k\geq p}\int^{*}_{S_{p}}f_{n_{k}}(x_{n_{k}})d\mu\;,$$
or
$$r\geq I_{f^{**}}(x)-\sup_{\sigma\in \Sigma, \;\sigma=(S_{p})_{p}}\limsup_{p}\sup_{k\geq p}\int^{*}_{S_{p}}-f_{n_{k}}(x_{n_{k}})d\mu\;,$$
equivalently:
$$r\geq I_{f^{**}}(x)-\delta^{+}((-f_{n_{k}}(x_{n_{k}}))_{k}).$$ 
Therefore, observing that if $(v_{n})_n$ is a subsequence of $(u_{n})_n$ we have $\delta^{+}((v_{n})_{n})\leq \delta^{+}((u_{n})_{n})$, we get
$$r\geq I_{f^{**}}(x)-\delta^{+}((-f_{n}(x_{n}))_{n}).$$ 
The proof of Theorem \ref{thm6.01} is complete.$\;\Box$ \\\\
Corollary \ref{thm6.0} is an immediate consequence of Theorem \ref{thm6.01} and Theorem \ref{prop991}. Indeed since $-f_{n}(x_{n})\leq f^{-}_{n}(x_{n})$ we have
$0\leq \displaystyle\delta^{+}((-f_{n}(x_{n}))_{n})\leq \delta^{+}((f^{-}_{n}(x_{n}))_{n})=0$.\\
The proof of Corollary \ref{thm6.1} is an immediate consequence of Theorem \ref{thm6.01}, Corollary \ref{thm6.0} and Proposition \ref{prop00}.$\;\Box$ \\ 
Proof of  Corollary \ref{thm6.1111}.\\
Let us suppose now that the sequence $(f_{n})_n$ satisfies the $\mathcal{T}$-Ioffe's criterion and: $$\displaystyle r> seq\;\mathcal{T}-li_{e}I_{f_{n}}(x)\;.$$ 
Let a subsequence $(f_{n_{k}})_k$ and  a sequence $(x_{k})_k$ $\mathcal{T}$-converging to $x$ such $r>\displaystyle\sup_{k} I_{f_{n_{k}}}(x_{k})$. Since $(f_{n})_n$ satisfies the $\mathcal{T}$-Ioffe's criterion at $x$, the sequence $(x_{k})_k$ has the lower compactness property with respect to $(f_{n_{k}})_k$ . If $g=seq\;\sigma-li_{e} f_{n_{k}}$ the last part of the Corollary  \ref{thm6.0} applied to the sequence $(f_{n_{k}})_k$ gives:
$$r>\liminf_{k} I_{f_{n_{k}}}(x_{k})\geq I_{g^{**}}(x)\geq I_{f^{**}}(x)\;.$$
This last inequality being valid for any $\displaystyle r> seq\;\mathcal{T}-li_{e}I_{f_{n}}(x)\;,$ the proof of Corollary \ref{thm6.1111} is complete.$\;\Box$\\\\
The following result is of interest even if $dim(E)=1$, in the sequel we will use it. It can be proved directly using the properties of uniform integrability and of convergence in local measure. 
\begin{lem}\label{890} Suppose the Banach $E$ is separable. Let a sequence $(u_{n})_n$ of real valued measurable functions converging in local measure to $0$ and an uniformly integrable sequence  $(z_{n})_n$ of $E$-valued measurable functions such that $z_{n}= u_{n}.y_{n}$.  If the sequence $(y_{n})_n$ is bounded in some $L_{p}(\Omega,E,\beta\mu)$ for some measurable positive valued function $\beta$ and $1\leq p\leq\infty$,  then the sequence $(z_{n})_n$ strongly converges to the origin of $L_{1}(\Omega, E)$.
\end{lem} 
Alternative proof of Lemma \ref{890}. It suffices to prove that every subsequence $(z_{n_{k}})_k$ of $(z_{n})_n$ admits a strongly converging subsequence to the origin. The sequence  $(\Vert y_{n}\Vert)_n$ is bounded in some $L_{p}(\Omega,{\R},\beta\mu)$ $1\leq p\leq\infty$, due to Corollary \ref{lem3.11}, it is possible from any subsequence  $(\Vert y_{n_{k}}\Vert)_k$ to extract a subsequence  $(\Vert y_{n_{k_{l}}}\Vert)_l$ which converges in the Biting sense to some $\xi\in L_{p}(\Omega,{\R}, \beta\mu)$. Moreover by extraction of  a subsequence one may suppose that $(u_{n_{k_{l}}})_l$ converges almost everywhere to $0$.  Let  $z^{*}\in  L_{\infty}(\Omega, {\R}, \mu)$, and let the integrands defined on $\Omega\times{\R}$ by
$$f_{l}(\omega,r)=u_{n_{k_{l}}}. r. z^{*}(\omega) \;\;\mbox{and}\;\;f=0= li_{e} f_{l}\;,\;\mbox{then}\;\;\displaystyle f_{l}(\Vert y_{n_{k_{l}}})\Vert= u_{n_{k_{l}}}\Vert y_{n_{k_{l}}}\Vert. z^{*}= \Vert z_{n_{k_{l}}}\Vert.z^{*}\;.  $$ 
Moreover for every measurable set $A\;$,
$\int_{A}\vert f_{l}(\Vert y_{n_{k_{l}}}\Vert)\vert d\mu\leq \Vert z^{*}\Vert_{\infty}.\Vert z_{n_{k_{l}}}1_{A}\Vert_{1}$. The sequence $(z)_n$ being uniformly integrable, the sequence  $(f_{l}(\Vert y_{n_{k_{l}}}\Vert))_l$ is uniformly integrable. Therefore the sequence  $(\Vert y_{n_{k_{l}}}\Vert)_l$ has the lower compactness property with respect to $(f_{l})_l$. Due to Corollary  \ref{thm6.1},
$$\liminf_{l}\int_{\Omega} \Vert z_{n_{k_{l}}}\Vert.z^{*} d\mu=\liminf_{l} I_{f_{l}}(\Vert y_{n_{k_{l}}}\Vert)\geq I_{f}(\xi)=0.$$
This inequality being valid for every $z^{*}\in  L_{\infty}(\Omega, {\R})$ we deduce that $(\Vert z_{n_{k_{l}}}\Vert)_n$ weakly (therefore strongly) converges to $0$. The proof of Lemma \ref{890} is complete. $\;\Box$ 

\section{Sequential strong-weak lower semicontinuity.}
In this section we consider a topological space  $(E,\tau)$, a separable Banach space $F$, with norm $\Vert .\Vert_{F}$, and its weak topology  $\sigma_{F}$. We will use  two topological spaces $(\mathcal{X}, \mathcal{S})$ and $(\mathcal{Y}, \mathcal{T})$ such  that $\mathcal{X}\subseteq L_{0}(\Omega, E)$ and  $\mathcal{Y}\subseteq L_{0}(\Omega, F)$. 
The topologies $\mathcal{S}$ and $\mathcal{T}$ verify the following assumptions: \\\\
$(H_{1})$ Every sequence $\mathcal{S}$-converging to $x\in \mathcal{X} $ admits a subsequence which converges to $x$ almost everywhere.\\\\
$(H^{'}_{1})$ Every sequence $\mathcal{S}$-converging to $x\in \mathcal{X}$ converges to $x$ in local measure.\\\\
$(H_{2})$ Every sequence $\mathcal{T}$-converging  to $y\in \mathcal{Y}$ is a converging sequence to $y$ in the Biting sense (see Definition \ref{def6.0}).\\\\
$(H_{3})$  The sets $\mathcal{X}$ and $\mathcal{Y}$ are decomposable. Moreover given a sequence $(A_{n})_n$ of measurable subsets such for every measurable set $A$ of finite measure $\displaystyle \lim_{n}\mu(A\cap A_{n})=0$, and a sequence $(x_{n},y_{n})_n$ of elements of $\mathcal{X}\times \mathcal{Y} $, $\mathcal{S}\times\mathcal{T}$-converging to $(x,y)\in \mathcal{X}\times \mathcal{Y}$, then the sequence $(w_{n})_n$, with $w_{n}=(x,y)1_{A^{c}_{n}}+(x_{n},y_{n}) 1_{A_{n}}$  $\mathcal{S}\times\mathcal{T}$-converges to $(x,y)$.\\\\
Given an integrand  $f:\Omega\times E\times F\to \overline{{\R}}$, the  integrand $f^{**}: \Omega\times E\times F\to \overline{{\R}}$ is in this section, the partial Fenchel-Moreau biconjugate defined for every $(\omega, e, e^{'})\in  \Omega\times E\times F$  by:
$$ f^{**}_{\omega}(e,\; e^{'})=f_{\omega}(e,\; .)^{**}(e^{'})\;.$$
\begin{thm}  \label{thm 1011}  Suppose that assumptions $(H_{1})$ and $(H_{2})$ hold. Let $x\in \mathcal{X}$ and $y\in \mathcal{Y}$, and let $f:\Omega\times E\times F\to \overline{{\R}}$ be a $\mathbb{T}\otimes \mathcal{B}(E\times F)$-measurable integrand.
Given a sequence $((x_{n}, y_{n}))_n$ $\mathcal{S}\times\ \mathcal{T}$-converging to $(x,y)$ and satisfying:\\
$(a)$ the sequence $( y_{n})_n$ is Nagumo tight,\\ 
$(b)$ the sequence $(f^{-}(x_{n}, y_{n}))_n$ is eventually bounded in $L_{1}(\Omega, {\R})$, \\
then with  the $\tau\times \sigma_{F}$ sequential lower semicontinuous regularization $g$ of $f$ defined for every $\omega\in \Omega$ by $g_{\omega}=seq\;\tau\times \sigma_{F}- li_{e}f_{\omega}$,
$$\quad\quad\quad\quad\quad\quad\quad\quad\quad  \liminf_{n}\; I_{f}(x_{n},y_{n})\geq I_{{g}^{**}}(x,y)-\delta^{+}((-f(x_{n}, y_{n}))_{n})\;.\;\quad\quad\quad\quad\quad\quad\quad\quad\quad (4)$$ 
One has in addition $\liminf_{n}\; I_{f}(x_{n},y_{n})\geq I_{f}(x,y)\;,$ when the additional conditions are fulfilled:\\
$(c)$ the sequence   $((x_{n}, y_{n}))_n$ satisfies the lower compactness property with respect to $f$,\\
$(d)$ for almost every $\omega\in \Omega$, $g_{\omega}^{**}(x(\omega),y(\omega))=f_{\omega}(x(\omega),\;y(\omega))$. 
\end{thm}
\begin{rem} \label{rem10} Remark (since in this case $f=g$), that the condition $(d)$ is satisfied at $(x, y)$ when the following two conditions hold for every  $\omega\in \Omega$,\\
$(i)$  the function $f_{\omega}$ is $\tau\times \sigma_{F}$-sequentially lower semicontinuous at each point of $\{x(\omega)\}\times F$.\\
$(j)$ $f^{**}_{\omega}(x(\omega),\;y(\omega))=f_{\omega}(x(\omega),\;y(\omega))$.\\
Moreover the condition $(d)$ is satisfied at any $(x,z)$ for $z\in Y$ when $(i)$ holds and we replace the condition $(j)$ by:\\ 
$(k)$ $f_{\omega}(x(\omega), .)$ is proper convex.
\end{rem}
(Indeed in this last case $f_{\omega}(x(\omega), .)$ is proper convex norm lower semicontinuous, hence we have $f_{\omega}(x(\omega), .)=f_{\omega}(x(\omega), .)^{**}$, \cite{110} Theorem 3.44). 
Conditions $(i)$ and $(k)$ first appear in Balder's work on seminormality and sequential semicontinuity \cite{111} Theorem 4.9 and Theorem 4.12. The authors in \cite{3} chapter 8.1  Theorem  8.1.6 obtain a semicontinuity result with a mild assumption of type $(i)$, conditions $(k)$ and $(c)$,  a Nagumo tightness assumption too (see Proposition \ref{prop 919}) and the topology $\mathcal{T}$ considered being $\sigma(L_{1}(\Omega,F), L_{\infty}(\Omega,F^{*}))$, moreover the proof is very distinct. 

\begin{rem} In the condition $(i)$ one cannot replace the weak topology of $F$ by the strong topology as show the example 8.1.8 in \cite{3}.
\end{rem}
It may happen that condition $(d)$ holds for $f$ and fails for $f^{**}$: because the semicontinuity condition $(i)$ fails as shows the following example.
\begin{ex} \label{ex1} $\Omega=(0, 1)$ endowed with the Lebesgue tribe and the Lebesgue measure $dt$, with $E=F={\R}$ and the continuous integrand $f$ is defined on $(0, 1)\times {\R}^{2}$ by
$$f(\omega,s,t)=\max(-\vert s\vert.\vert t\vert\;, -1)\;.$$
\end{ex} 
Then the lower semicontinuous regularization $g$ of $f$ is equal to $f$, hence $f$ verifies $(i)$ at any point $(x,y)$. Moreover when $s=0$, $f(s, .)^{**}=f(s, .)=0$:  $f$ verifies $(j)$ (therefore $(d)$) at any point $(0,t)\in\{0\}\times{\R}$; but $f(s, .)^{**}=-1$ when $s\neq 0$. Therefore $f^{**}$ is not lower semicontinuous at any point $(0,t)\in\{0\}\times{\R}$. The lower semicontinuous regularization $h$ of $f^{**}$ is $h(s,t)=-1=h^{**}(s,t)$, therefore for any  $(0,t)\in\{0\}\times{\R}, f^{**}(0,t) \neq h^{**}(0,t)$, this proves that the condition $(d)$ fails for $f^{**}$ at any point  $(0,t)\in\{0\}\times{\R}$.\\\\
Proof of Theorem \ref{thm 1011}. Let $r\in {\R}$ such that $\displaystyle\liminf_{n}\; I_{f_{n}}(x_{n},y_{n})<r$.  Since  $(H_{1})$ holds, extracting  subsequences we may reduce to the case when the sequence  $(x_{n})_n$ converges almost everywhere to $x$ and $\sup_{n}\; I_{f_{n}}(x_{n},y_{n})<r$. Setting $\displaystyle f_{n}(\omega,e)=f(\omega,x_{n}(\omega),e)$, the sequence $(y_{n})_n$ is supposed Nagumo tight and with $(b)$ is supposed to verify the boundedness property with respect to the sequence of measurable integrands $(f_{n})_n$.  Due to assumption $(H_{2})$, the topological space $(\mathcal{Y},\mathcal{T})$ verifies the property $(P)$ of the previous section. The use of Theorem \ref{thm6.01} gives, with $h=seq\; \sigma_{F}-li_{e} f_{n}$:
$$r>\liminf_{n}\; I_{f_{n}}(y_{n})\geq I_{h^{**}}(y)-\delta^{+}((-f(x_{n}, y_{n}))_{n})\;.$$ 
The above inequality is valid for every $r>\liminf_{n}\; I_{f}(x_{n}, y_{n})$. Moreover we remark that in addition, $g(x,.)\leq h$, hence we have almost everywhere $g^{**}(x,y)\leq h^{**}(y)$.  This proves the validity of the inequality $(4)$.\\
If  assumption $(c)$ holds, then $0\leq \delta^{+}((-f(x_{n}, y_{n}))_{n})\leq \delta^{+}((f^{-}(x_{n}, y_{n}))_{n})=0$. The Assumption $(d)$ and formula $(4)$ gives the second semicontinuity result. The proof of Theorem \ref{thm 1011} is complete.$\;\Box$ 
\begin{Def} (see \cite{70}) A measurable integrand $f:\Omega\times E\times F\to \overline{{\R}}$ satisfies the Ioffe's criterion at $(x,y)\in \mathcal{X}\times \mathcal{Y}$ if for every sequence $((x_{n}, y_{n}))_n$ $\mathcal{S}\times \mathcal{T}$-converging to $(x,y)$ such that the sequence $(I_{f}(x_{n},y_{n}))_n$ is bounded above, the sequence $(f^{-}(x_{n}, y_{n}))_n$ is uniformly integrable.
\end{Def}
The corollary below is an extension of the Ioffe's result \cite{70}.
\begin{coro} \label{thm 101} Suppose $(E,\tau)$ is a metrisable topological space,  $F$ is a reflexive separable Banach space and assumptions $(H^{'}_{1})$ and $(H_{2})$ hold. 
Let $f:\Omega\times E\times F\to \overline{{\R}}$ be a $\mathbb{T}\otimes \mathcal{B}(E\times F)$-measurable integrand satisfying for almost every $\omega\in \Omega$ the conditions:\\
$(i)$  the function $f_{\omega}$ is sequentially $\tau\times \sigma_{F}$-lower semicontinuous at each point of $\{x(\omega)\}\times F$.\\
$(j)$ $f^{**}_{\omega}(x(\omega),\;y(\omega))=f_{\omega}(x(\omega),\;y(\omega))$.\\
Then, whenever the integral functional $I_f$ satisfies the Ioffe's criterion at $(x,y)\in \mathcal{X}\times \mathcal{Y}$ it is $\mathcal{S}\times\mathcal{T}$-sequentially lower semicontinuous at this point.\\ Conversely, suppose that the measure $\mu$ is atomless, assumption $(H_{3})$ holds and $I_{f}(x,y)\neq -\infty$. Then the Ioffe's criterion at $(x,y)$ is necessary for the sequential lower semicontinuity of $\;I_f$ at this point. 
\end{coro}
Proof of Corollary \ref{thm 101}. Let us prove first the sufficiency part. When the conditions $(i)$ and $(j)$ hold, then due to the Remark \ref{rem10} the condition $(d)$ holds. Let $r\in {\R}$ and a sequence $((x_{n},y_{n}))_n$ $\mathcal{S}\times \mathcal{T}$-converging to $(x,y)$ such for every integer $n$
$$I_{f}(x_{n},y_{n})\leq r\;.\;$$
The Ioffe criterion ensures that the sequence $((x_{n},y_{n}))_n$ verifies the lower compactness property with respect to the integrand $f$.
Due to $(H_{2})$  and Proposition \ref{prop00} the sequence $(y_{n})_n$ is Nagumo tight. Since $(H^{'}_{1})$ holds and every sequence converging in local measure admits a subsequence converging almost everywhere (\cite{2}) then $(H_{1})$ holds. Therefore from Theorem \ref{thm 1011} we get:
$$I_{f}(x,y)\leq \liminf_{n}\; I_{f_{n}}(x_{n},y_{n})\leq r$$
This proves that the integral functional $I_f$ is sequentially $\mathcal{S}\times \mathcal{T}$-lower semicontinuous at  $(x,y)$.\\
The proof of the necessity part is very similar to the proof given by A. D. Ioffe in \cite{70} Theorem 1.  Indeed
if $I_f$ is $\mathcal{S}\times\mathcal{ T}$-sequentially lower semicontinuous at $(x,y)$ and the Ioffe's criterion is not true, there exist a real number $r$, a sequence $((x_{n},y_{n}))_n$ $\mathcal{S}\times \mathcal{T}$-converging to $(x,y)$ such the sequence $(f^{-}((x_{n},y_{n})))_n$ is not uniformly integrable and for every $n$, $I_{f}(x_{n},y_{n})\leq r$.
Since $I_f$ is sequentially lower semicontinuous at $(x,y)$ we get: $I_{f}(x,y)\leq r<\infty$, and by assumptions $I_{f}(x,y)>-\infty$, thus $f(x,y)$ is integrable. The sequence $(f^{-}(x_{n},y_{n}))_n$ being not uniformly integrable and the measure considered being atomless, with the help of Theorem \ref{prop991} $(b)$ (\cite{67}, Proposition 1.7), we have:
$$\delta^{+}((f^{-}(x_{n},y_{n}))_n)>0\;.$$
Therefore there exist $\epsilon>0$, a subsequence $((x_{n_{k}},y_{n_{k}}))_k$ and a decreasing sequence $(A_{k})_k$ with a negligible intersection satisfying:
$$\int^{*}_{A_{k}}f^{-}(x_{n_{k}},y_{n_{k}}) d\mu\geq \epsilon\;.$$
Set $B_{k}=\{f(x_{n_{k}},y_{n_{k}})\leq 0\}$, $C_{k}=A_{k}\cap B_{k}$, $(x^{'}_{k},y^{'}_{k})=(x,y)1_{C_{k}^{c}}+(x_{n_{k}},y_{n_{k}})1_{C_{k}}$. Due to  $H_3$, since $C_{k}\subset A_{k}$, the sequence $((x^{'}_{k},y^{'}_{k}))_k$
is a sequence of elements of $\mathcal{X}\times \mathcal{Y}$ which $\mathcal{S}\times\mathcal{ T}$-converges to $(x,y)$ and satisfies eventually:
$$I_{f}(x^{'}_{k},y^{'}_{k})-I_{f}(x,y)=\int^{*}_{C_{k}}(f(x_{n_{k}},y_{n_{k}})-f(x,y)) d\mu=-\int^{*}_{A_{k}}f^{-}(x_{n_{k}},y_{n_{k}})d\mu-\int_{C_{k}}f(x,y)\leq -\frac{\epsilon}{2}\;.$$
This contradicts the sequential $\mathcal{S}\times\mathcal{ T}$-lower semicontinuity of $I_f$ at $(x,y)$. The proof of Corollary \ref{thm 101} is complete. $\;\Box$ 
\begin{rem}\label{r1} Under the assumption $(Q)$ on the topology $\mathcal{S}\times\mathcal{ T}$ instead of  ($H_3$), the necessity part is also a consequence of Theorem \ref{thm6.2}.
\end{rem}
\begin{rem}\label{r2} It may happen that at a given point $(x,y)$ the integral functional $I_f$ is sequentially $\mathcal{S}\times\mathcal{T}$-lower semicontinuous at this point  and the sequential $\mathcal{S}\times\mathcal{T}$-lower semicontinuity of $I_{f^{**}}$ fails at the same  point $(x,y)$. 
\end{rem} 
 Indeed consider the integrand $f$ of Example \ref{ex1} and a point $(0,y)\in\{0\}\times L_{2}((0,1), {\R})$ and keep $\mathcal{X}=\mathcal{Y}=L_{2}((0,1), {\R})$. The integrand $f$ verifies the condition $(d)$ at $(0, y)$, let us prove that the integral functional $I_f$ is sequentially strong-weak continuous on $L_{2}((0,1), {\R}^{2})$ at $(0,y)$. Let us give a direct proof. Given a strong-weak converging sequence $(x_{n}, y_{n})_n$ to $(0,y)$, from the H\"older inequality the sequence $(\vert x_{n}\vert. \vert y_{n}\vert)_n$ is uniformly integrable since  $(x_{n})_n$ is $2$-equi-integrable and for every measurable set $A$:
 $$\int_{A}\vert x_{n}\vert. \vert y_{n}\vert ds\leq \sup_{n}\Vert y_{n}\Vert_{2} \sup_{n}\Vert x_{n}1_{A}\Vert_{2}\;.$$
Due to Lemma \ref{890} the sequence $(\vert x_{n}\vert. \vert y_{n}\vert)_n$ strongly converges to the origin in $L_{1}((0,1), {\R})$ and we get with $f^{-}=\inf(\vert s\vert.\vert t\vert,\;1)=\vert f\vert $,
 $$0\leq \limsup_{n}\int^{1}_{0}f^{-}(x_{n}(s), y_{n}(s))ds=\limsup_{n}\int^{1}_{0}\vert f\vert (x_{n}(s), y_{n}(s))ds\leq \lim_{n}\int^{1}_{0}\vert x_{n}\vert. \vert y_{n}\vert ds=0 \;.$$
This proves that the sequence $(f(x_{n}, y_{n}))_n$ strongly converges in $L_{1}((0,1), {\R})$ (and therefore the Ioffe criterion holds for $f$ and $-f$). Thus $I_f$ is sequentially strong weak continuous at $(0, y)$ because 
$$\lim_{n}\int^{1}_{0}f(x_{n}(s), y_{n}(s))ds=\lim_{n}\int^{1}_{0}f^{-}(x_{n}(s), y_{n}(s))ds=0=I_{f}(0,y)\;.$$
$I_{f^{**}}$ is not semicontinuous at $(0,y)$: since $f^{**}(s,t)=0$ if $s=0$ and $f^{**}(s,t)=-1$ if $s\neq0$ then 
$$I_{f^{**}}(0, y)=0>\liminf_{n} I_{f^{**}}(n^{-1}, y)=-1\;.$$

\section{Lower compactness properties, usual examples. }
The Banach space $E$ is supposed separable. Hereafter $(f_{n})_n$ is a sequence of extended real valued measurable integrands defined on $\Omega\times E$.
\begin{Def} \label{Def381} Let $\mathcal{X}$ be a set. A family $\mathcal{B}$ of subsets of $\mathcal{X}$ is said hereditary if for every element $X$ of $\mathcal{B}$, every subset $Y$ of $X$ is an element of $\mathcal{B}$ too. A bornology on $\mathcal{X}$ is an hereditary family $\mathcal{B}$ of subsets of $X$ stable by finite unions which covers $X$.
\end{Def}
When  $(\mathcal{X}, \mathcal{T})$ is a locally convex topological space, the most usual bornologies are the family of $ \mathcal{T}$-bounded sets it is called the Fr\'echet bornology, it is denoted by  $\mathcal{B}_{F}$; the family of sequentially $\mathcal{T}$-relatively compacts sets is called the  $ \mathcal{T}$-Hadamard bornology, it is denoted by $\mathcal{B}_{H}$, and if  more generally $\tau$ is a topology on $\mathcal{X}$ the family of $\tau$-relatively sequentially compacts sets (denoted by $\mathcal{B}_{\tau}$) is a bornology.  For more simplicity and clarity, in the sequel $\mathcal{X}$ will be a normed space of Orlicz type: more precisely a Lebesgue space endowed with its Fr\'echet bornology or its Hadamard bornologies associated to the strong or weak (star)topologies.
\begin{Def} \label{Def379} Let $\mathcal{X}$ be a normed space contained in $L_{0}(\Omega, E)$ endowed with a bornology  $\mathcal{B}$. A sequence  $(f_{n})_n$ of extended real valued integrands defined on $\Omega\times E$ is said to have the $\mathcal{B}$-lower compactness property (denoted $\mathcal{B}$-lcp) if for every element $X$ of $\mathcal{B}$, every sequence $(x_{n})_n$ of elements in $X$ has the lower compactness property respect to $(f_{n})_n$.
The $\mathcal{B}_{F}$-lcp is called  Fr\'echet-lcp, the  $\mathcal{B}_{H}$-lcp is called Hadamard-lcp. Given a topology $\mathcal{\tau}$ on $\mathcal{X}$, the  $\mathcal{B}_{\mathcal{\tau}}$-lcp is simply denoted by  $\mathcal{\tau}$-lcp.
\end{Def}
Before to give concrete examples let us prove the following property.
\begin{prop}\label{prop379} Let a bornology  $\mathcal{B}$ on a normed subspace $\mathcal{X}$ of $L_{0}(\Omega, E)$. A sequence  $(f_{n})_n$ of extended real valued integrands defined on $\Omega\times E$ has the $\mathcal{B}$-lcp if and only if each subsequence  $(f_{n_{k}})_k$ has the $\mathcal{B}$-lcp. When $\mathcal{\tau}$ is a topology on $\mathcal{X}$, a sequence  $(f_{n})_n$ of integrands has the $\mathcal{B}_{\mathcal{\tau}}$-lcp if and only if every $\mathcal{\tau}$-converging sequence $(x_{n})_n$ has the lcp respect to $(f_{n})_n$. 
\end{prop}
Proof of Proposition \ref{prop379}. Let $X$ be an element of $\mathcal{B}$ and $(x_{k})_k$ be a sequence of elements of $X$. Let a subsequence  $(f_{n_{k}})_k$ of $(f_{n})_n$. We want to prove that the sequence $(x_{k})_k$ has the lcp respect to the sequence  $(f_{n_{k}})_k$. Consider the sequence $(y_{m})_m$ defined for $m\geq n_{0}$ by: $y_{m}=x_{k}\;\mbox{if}\;  n_{k}\leq m<n_{k+1}$.
Since $(f_{n})_n$ has the $\mathcal{B}$-lcp, then the sequence $(f^{-}_{n}(y_{n}))_n$ is uniformly integrable. But $(f^{-}_{n_{k}}(x_{k}))_k$ is a subsequence of  $(f^{-}_{n}(y_{n}))_n$, therefore  $(f^{-}_{n_{k}}(x_{k}))_k$ is uniformly integrable too. Now suppose that the sequence  $(f_{n})_n$ has not the $\mathcal{B}_{\tau}$-lcp there exists a subsequence $(f_{n_{k}})_k$ a $\mathcal{\tau}$-converging sequence  $(x_{k})_k$ such the sequence  $(f^{-}_{n_{k}}(x_{k}))_k$ is not uniformly integrable. Define  $y_{m}=x_{k}\;\mbox{if}\;  n_{k}\leq m<n_{k+1}$. Then $(y_{m})_m$ $\mathcal{\tau}$-converges too and has not the lcp respect to $(f_{m})_m$. The proof of Proposition \ref{prop379} is complete.$\;\Box$\\\\
Given a Young integrand $\phi$ recall that the Orlicz space $L_{\phi}(\Omega, E, \mu )= {\R}_{+}I_{\phi}^{\leq 1}$ becomes a decomposable normed space endowed with
the Minkowski gauge associated to the sublevel set $I_{\phi}^{\leq 1}$, $\Vert x\Vert_{\phi}=\inf\{t>0: x\in t.I_{\phi}^{\leq 1}\}$. For $t>0$, set $\phi_{t}(e)=\phi(te)$.
\begin{prop} \label{coro554} Let  $\phi$ be a Young integrand and a sequence ${(f_{n})}_n$ of measurable integrands such that:\\
for every $\epsilon>0$ and $\lambda>0$ there exists a sequence $(u_{n})_n$ of non negative uniformly integrable functions satisfying eventually
$$f_{n}\geq -\epsilon.\phi_{\lambda}(.) -u_{n} \;.$$ 
Then the sequence of integrands ${(f_{n})}_n$ has the Fr\'echet-lcp on $L_{\phi}(\Omega, E)$.
\end{prop}
Proof of Proposition \ref{coro554} (When the measure is atomless, this condition is characteristic \cite{6}).
Let ${(x_{n})}_n$ be a norm bounded sequence by $m$ in $L_{\phi}(\Omega, E)$. Then  $\sup_{n}I_{\phi}(m^{-1} x_{n})\leq 1$.
If the groth condition of Proposition \ref{coro554} holds with $\lambda=m^{-1}$ then
$$\delta^{+}(({f_{n}^{-}(x_{n})})_n)\leq \epsilon+\delta^{+}((u_{n})_n)=\epsilon \;.$$
Therefore  $\delta^{+}(({f_{n}^{-}(x_{n})})_n)=0$. The sequence  $({f_{n}^{-}(x_{n})})_n$ being trivially bounded in $L_{1}(\Omega, {\R})$ it is,  from Theorem \ref{prop991}  $(a)$, eventually equi-integrable thus eventually uniformly integrable. We deduce that the sequence ${(f_{n})}_n$ has the Fr\'echet-lcp on  $L_{\phi}(\Omega, E)$ .$\;\Box$\\\\
From the above result if $\mathcal{X}=L_{p}(\Omega, E)$ $1\leq p< \infty$ and $\phi$ is the Young integrand $\phi(\omega, e)=p^{-1}\Vert e\Vert^{p}$, we obtain: 
\begin{coro} \label{coro157} Let $1\leq p<+\infty$ and a sequence ${(f_{n})}_n$ of measurable integrands such that:\\
for every $\epsilon>0$ there exists a sequence $(u_{n})_n$ of non negative uniformly integrable functions satisfying eventually
$$ f_{n}\geq -\epsilon \Vert.\Vert^{p}-u_{n} \;.$$\
Then the sequence of integrands ${(f_{n})}_n$ has the Fr\'echet-lcp on $L_{p}(\Omega, E)$.
\end{coro}
The case $L_{\infty}(\Omega, E)$ may be of interest:
\begin{coro} \label{coro154} Let $p=\infty$ and $B$ be the unit ball of $E$. Let a sequence ${(f_{n})}_n$ of measurable integrands verifying the property:\\
For each  $\lambda>0$ there exists a uniformly integrable sequence of functions $(u_{n})_n$ such that (eventually):
$$ \inf_{e\in \lambda^{-1} B} f_{n}(., e)\geq -u_{n}\;.$$
Then the sequence of integrands ${(f_{n})}_n$ has the Fr\'echet-lcp on $L_{\infty}(\Omega, E)$.
\end{coro}
Proof of Corollary \ref{coro154}. $L_{\infty}(\Omega, E)$ is the Orlicz space associated to the Young integrand:  $\phi(\omega, e)=\iota_{B}(e)$ for which $\phi_{\lambda}(\omega, e)=\iota_{\lambda^{-1}B}(e)$ $\;\Box$.
\begin{prop} \label{coro553} Let  $\phi$ be a Young integrand and a sequence ${(f_{n})}_n$ of measurable integrands such that:\\
for every $x\in L_{\phi}(\Omega, E)$, $\epsilon>0$ there exist $\lambda>0$,  a positive constant $c\geq 1$ a sequence $(u_{n})_n$ of non negative bounded integrable functions satisfying $\delta^{+}((u_{n})_n)<\epsilon$ and eventually
$$f_{n}\geq -c.\phi_{\lambda}(.-x) -u_{n} \;.$$ 
Then the sequence of integrands ${(f_{n})}_n$ has the (strong) Hadamard-lcp at $x\in L_{\phi}(\Omega, E)$.
\end{prop}
Proof of Proposition \ref{coro553}. Let ${(x_{n})}_n$ be a norm converging sequence to $x\in L_{\phi}(\Omega, E)$. If the growth condition of Proposition \ref{coro553} holds  there exists  positive constants  $\lambda$, $c>1$ and a bounded sequence of integrable functions $(u_{n})_n$ with $\delta^{+}((u_{n})_n)<\epsilon$ satisfying eventually:\\
$$f_{n}\geq -c.\phi_{\lambda}(.-x) -u_{n} \;.$$
Therefore eventually:
$$f_{n}(x_{n})^{-}\leq \frac{1}{2}c\phi_{\lambda}(2(x_{n}-x))+ u_{n} \;. (*)$$
Since ${(x_{n})}_n$ converges to $x$, we have for $n$ large enough 
$$\Vert x_{n}-x\Vert_{\phi}\leq (2\lambda \epsilon^{-1}c)^{-1}\;,\;\mbox{or eventually}\;\; I_{\phi}(\lambda(\epsilon^{-1} c)2( x_{n}-x))= I_{\phi_{\lambda}}((\epsilon^{-1} c)2( x_{n}-x))\leq 1$$
then eventually: 
$$I_{\phi_{\lambda}}( 2(x_{n}-x))=I_{\phi_{\lambda}}(\epsilon c^{-1}\epsilon^{-1} c)2( x_{n}-x))\leq \epsilon  c^{-1} I_{\phi_{\lambda}}((\epsilon^{-1} c)2( x_{n}-x))\leq\epsilon  c^{-1}\;.$$
From $(*)$ we get the boundedness of the sequence $({f_{n}^{-}(x_{n})})_n$ in $L_{1}(\Omega, {\R})$ and:
$$\delta^{+}(({f_{n}^{-}(x_{n})})_n)\leq \frac{1}{2}cI_{\phi_{\lambda}}( 2(x_{n}-x))+\delta^{+}((u_{n})_n)) \leq \frac{3}{2}\epsilon \;.$$
This last inequality being valid for every $\epsilon>0$ then from Theorem \ref{prop991}  $(a)$ we obtain the equi-integrability of the sequence $({f_{n}^{-}(x_{n})})_n$, thus the uniform integrability of this sequence. $\;\Box$ 
\begin{coro} \label{coro156} Let $1\leq p<+\infty$ and a sequence ${(f_{n})}_n$ of measurable integrands verifying the property:\\
For each $\epsilon>0$ there exists a positive constant $c_{\epsilon}$ and a non negative bounded sequence of integrable functions $(u_{n})_n$ with $\delta^{+}((u_{n})_n)<\epsilon$ and such that (eventually):\\
$$ f_{n}\geq -c_{\epsilon}\Vert.\Vert^{p}-u_{n} \;.$$
Then the sequence of integrands ${(f_{n})}_n$ has the Hadamard-lcp on $L_{p}(\Omega, E)$.
\end{coro}
Proof of Corollary \ref{coro156}. For every $x\in L_{p}(\Omega, E)$ use Proposition \ref{coro553} and the following inequality: $\Vert e\Vert^{p}\leq 2^{p-1}( \Vert e-x\Vert^{p}+\Vert x\Vert^{p}$).
\begin{coro} \label{coro155} Let $p=\infty$. Let a sequence ${(f_{n})}_n$ of measurable integrands verifying the property:\\
For every $x\in L_{\infty}(\Omega, E)$, $\epsilon>0$,  there exist $\lambda>0$, a bounded sequence of integrable functions $(u_{n})_n$ such that  $\delta^{+}((u_{n})_n)<\epsilon$ and eventually:
$$ \inf_{e\in x+\lambda^{-1} B}f_{n}(., e)\geq -u_{n} \;.$$
Then the sequence of integrands ${(f_{n})}_n$ has the strong Hadamard-lcp on $L_{\infty}(\Omega, E)$.
\end{coro}
Proof of Corollary \ref{coro155}.  It is a consequence of Proposition \ref{coro553},  remark that $L_{\infty}(\Omega, E)$ is the Orlicz space associated to the Young integrand:  $\phi(\omega, e)=\iota_{B}(e)$ for which $\phi_{\lambda}(\omega, e)=\iota_{\lambda^{-1}B}(e)$. $\;\Box$\\\\
In order to give another class of examples of sequences of integrands with the lower compactness property, let us consider a measurable integrand $f:\Omega\times E\to \overline{{\R}}$. We will use its differential quotient defined by the formula:
$$ [f](\omega, e_{0},e,r)=[f_{\omega}]( e_{0},e,r)=r^{-1}(f_{\omega}(e_{0}+re)-f_{\omega}(e_{0}))$$
If as usually, $I_{f}$ is the integral functional associated to the integrand $f$ and if $x_{0}$ is  such that $f(x_{0})$ is integrable remark that for every $x\in L_{0}(\Omega, E)$ we have the equality:
$$[I_{f}](x_{0}, x, r)=I_{[f](x_{0}\;, .\;r)} (x)\;.$$
\begin{Def} \label{defw.2} Given a Bornology $\mathcal{B}$ on a normed subspace $\mathcal{X}$ of $L_{0}(\Omega, E)$, an integrand $f$ is said to have the $\mathcal{B}$-differential lower compactness property on  $\mathcal{X}$ (denoted $\mathcal{B}$-dlcp) at $x_0$ if for every sequence ${(r_{n})}_n$ of positive numbers converging to $0$,  the sequence of differential quotients integrands $([f](x_{0},.,r_{n}))_n$ has the $\mathcal{B}$-lcp. 
If $\mathcal{B}$ is the Fr\'echet-Bornology (respectively the Hadamard bornology, respectively the bornology associated to a topology $\tau$) we call this notion Fr\'echet-differential lower compactness property (respectively Hadamard-dlcp, respectively  $\tau$-dlcp) at $x_0$.
\end{Def}
\begin{prop} \label{propw.200} Let $\phi$ be a Young integrand and $f$ be a measurable integrand satisfying the following condition: \\
 for every $\epsilon >0$, and every $\lambda>0$, there exists a family of non negative eventually integrable functions $\{u_{r}, r\in (0,1)\}$ uniformly integrable in $L_{1}(\Omega,{\R})$ and verifying eventually
$$[f](x_{0}, ., r)\geq -\epsilon\phi_{\lambda} -u_{r}\;.$$ 
Then $f$ has the Fr\'echet differential lcp on $L_{\phi}(\Omega, E)$ at $x_{0}\in L_{\phi}(\Omega, E)\;.$
\end{prop}
Proof. Let ${(r_{n})}_n$ be a sequence of positive numbers converging to $0$, using Proposition \ref{coro554} with $f_{n}=[f](x_{0}, ., r_{n})$ we deduce that  the sequence of integrands ${(f_{n})}_n$ has the Fr\'echet-lcp on $L_{\phi}(\Omega, E)$. The result being valid for every sequence of positive numbers converging to $0$, $f$ has the Fr\'echet-lcp on $L_{\phi}(\Omega, E)$. $\;\Box$\\\\
Given a  Young integrand, and an integrand $f$, we will use the following assumption:\\\\
{\it $(\mathcal{S}_{\phi, x_{0}})$ $x_{0}\in  L_{\phi}(\Omega, E)$, for almost every $\omega\in \Omega$, the function $f_{\omega}$ is Lipschitzian on every ball of $E$; there exist two positive constants $c$ and $\lambda_{0}$ such that for every $\lambda>\lambda_{0}$ there exist $\beta>0$ and an integrable function $u_{\lambda}\in  L_{1}(\Omega,{\R})$  verifying:
$$\sup_{e^{*} \in \partial^{C} f( x_{0}+e)}\phi^{*}(\lambda e^{*})\leq c \phi(\beta e)+u_{\lambda}\:.$$
Where $\partial^{C}$ is the Clarke subdifferential see \cite{130}. $(\mathcal{S}_{\phi})$ denotes  $(\mathcal{S}_{\phi, 0})$.}
\begin{prop} \label{propw.2} Let $\phi$ be a Young integrand and $f$ be a measurable integrand satisfying the condition $(\mathcal{S}_{\phi, x_{0}})$. 
Then $f$ has the Fr\'echet differential lcp on $L_{\phi}(\Omega, E)$ at $x_{0}\in L_{\phi}(\Omega, E)\;.$
\end{prop}
Proof. This is a consequence of  Proposition \ref{propw.200} and the following Lemma with $e^{'}=0$:
\begin{lem} \label{lem 819} A measurable integrand $f:\Omega\times E\to {\R}$ satisfying $(\mathcal{S}_{\phi, x_{0}})$ verifies also:
 for every $\epsilon>0$ and every $\lambda>0$, there exist $r_{\epsilon}>0$ an integrable function $u_{\epsilon, \lambda}$ such that for $0<r\leq r_{\epsilon}$ and $e, e^{'}:\in E$
$$\vert r^{-1}(f(x_{0}+re)-f(x_{0}+re^{'}))\vert\leq \epsilon \phi_{\lambda}(e-e^{'}) +u_{\epsilon, \lambda}$$
\end{lem}
Proof of Lemma \ref{lem 819}. Let $\epsilon>0$, $\lambda>0$. From Lebourg's Mean value theorem \cite{130}, for each $e, e^{'}\in E$, $0<r<1$, for each $\omega\in \Omega$ there exists $0<\theta_{n}(\omega)<1$,
 $y=x_{0}+\theta r(e-e^{'})$ and $x^{*}\in {\partial}^{C}f(y)$ such for every $\epsilon>0$ and $\lambda>0$:
$$ r_{n}^{-1}(f(x_{0}+re)-f(x_{0}+re^{'}))=\langle (\epsilon\lambda)^{-1} x^{*}, \epsilon\lambda (e-e^{'})\rangle\;. $$
Therefore due to the Young inequality, and condition $(\mathcal{S}_{\phi})$, when $(\epsilon\lambda)^{-1}\geq \lambda_{0}$ that is when  $\epsilon\leq \lambda^{-1}\lambda_{0}^{-1}$, there exist $\beta>0$ and an integrable function  $u_{(\epsilon \lambda)^{-1}}$ verifying eventually:
$$\vert r_{n}^{-1}(f(x_{0}+re)-f(x_{0}+re^{'}))\vert\leq \phi^{*}((\epsilon\lambda)^{-1} x^{*})+ \phi(\epsilon\lambda (e-e^{'}))\;,$$
$$\vert r_{n}^{-1}(f(x_{0}+re)-f(x_{0}+re^{'}))\vert\leq c \phi(\beta(\theta_{n} r(e-e^{'}))+u_{(\epsilon\lambda)^{-1}}+\epsilon\phi_{\lambda}(e-e^{'})\leq cr^{\frac{1}{2}} \phi(\beta r^{\frac{1}{2}}e)+\epsilon\phi_{\lambda}(e)+u_{(\epsilon\lambda)^{-1}}$$
 when we have eventually: $0<cr^{\frac{1}{2}}\leq \epsilon$ and $0<\beta r^{\frac{1}{2}}\leq \lambda$, we deduce:
$$\vert r^{-1}(f(x_{0}+re)-f(x_{0}+re^{'}))\vert\leq 2\epsilon \phi_{\lambda} +u_{(\epsilon\lambda)^{-1}}$$
This ends the proof of  Lemma \ref{lem 819}. $\;\Box$\\
End of the proof of Proposition \ref{propw.2}. Since $f$ satisfies the condition of Lemma \ref{lem 819}, it satisfies the growth condition of  Proposition \ref{propw.200} therefore $f$ has the Fr\'echet differential lcp on $L_{\phi}(\Omega, E)$ at $x_{0}\in L_{\phi}(\Omega, E)\;$. The proof of Proposition \ref{propw.2} is complete. $\;\Box$ 
\begin{coro} \label{rem 713} A measurable integrand $f:\Omega\times E\to {\R}$ satisfying the condition of  Lemma \ref{lem 819} (or clearly the condition $(\mathcal{S}_{\phi, x_{0}})$) at $x_{0}\in L_{\phi}(\Omega, E)$ verifies too eventually:\\
$$\sup_{\Vert x\Vert_{\phi}\leq 1, \Vert y\Vert_{\phi}\leq 1 }  \vert I_{f}(x_{0}+rx )- I_{f}(x_{0}+r y)\vert-mr\Vert x-y\Vert_{\phi}\leq 0\;.$$
\end{coro} 
 Proof. Let $B_{\phi}$ be the unit ball of $L_{\phi}(\Omega, E)$,  when $f$ satisfies the condition of  Lemma \ref{lem 819}, for every $\lambda>0$, there exist $r_{1}>0$ and $ u_{1,  \lambda^{-1}}$ such that  for $0<t<r_{1}$ and $x^{'}, y^{'}\in L_{\phi}(\Omega, E)$: 
$$ t^{-1}\vert(f(x_{0}+tx^{'})- f(x_{0}+ty^{'}))\vert \leq \phi(\frac{1}{2}\lambda^{-1}(x^{'}-y^{'})) +u_{1, \frac{1}{2}\lambda^{-1}}\;.$$
Therefore for every $x, y\in \lambda B_{\phi}$, taking $x^{'}= \Vert x-y\Vert_{\phi}^{-1}x$ and $y^{'}= \Vert x-y\Vert_{\phi}^{-1}y$ and $t=r\Vert x-y\Vert_{\phi}$ we get for $r\leq r_{1} \lambda^{-1}$
$$ \sup_{x, y\in  \lambda B_{\phi}}\Vert x-y\Vert^{-1}_{\phi}r^{-1}\vert (f(x_{0}+r x)- f(x_{0}+ry))\vert\leq\phi(\frac{1}{2}(x^{'}-y^{'}))+ u_{1, \frac{1}{2}\lambda^{-1}}$$
thus, for every $x, y\in \lambda B_{\phi}$, and $r\leq r_{1} \lambda^{-1}$,
 $$\Vert x-y\Vert^{-1}_{\phi}r^{-1}\vert (f(x_{0}+r x)- f(x_{0}+ry))\vert\leq \frac{1}{2} \phi(\frac{(x-y)}{\Vert x-y\Vert_{\phi}}) + u_{1, \frac{1}{2}\lambda^{-1}}\;.$$
thus by integration,  for every $x, y\in \lambda B_{\phi}$, and $r\leq r_{1} \lambda^{-1}$,
$$\Vert x-y\Vert_{\phi}^{-1}r^{-1}\vert(I_{f}(x_{0}+r x)- I_{f}(x_{0}+r y)\vert\leq \frac{1}{2}+\int_{\Omega}u_{1, \frac{1}{2}\lambda^{-1}} d\mu=m<\infty\;$$
therefore for $r\leq r_{1} \lambda^{-1}$:
$$ \sup_{x, y\in \in \lambda B_{\phi}}r^{-1}\vert(I_{f}(x_{0}+rx)- I_{f}(x_{0}+ry))\vert-\leq m \Vert x-y\Vert_{\phi}\;,$$
this proves that with $\eta=r_{1} \lambda^{-1}$  and $r\leq \eta$  we have:
$$\sup_{x, y\in \lambda B_{\phi}}  \vert I_{f}(x_{0}+rx )- I_{f}(x_{0}+r y)\vert- mr\Vert x-y\Vert_{\phi}\leq 0\;. \;\Box$$
\begin{prop} \label{coro170} Let $\phi$ be a Young integrand and $f$ be an integrand verifying $(\mathcal{S}_{\phi})$, then the condition $(\mathcal{S}_{\phi, x_{0}})$  holds at every point $x_{0}\in E_{\phi}(\Omega, E)=\{x\in L_{\phi}(\Omega, E): \forall \lambda>0, \phi(\lambda x)\in L_{1}(\Omega,{\R}, \mu)\}$. 
\end{prop}
Proof. Suppose that  $x_{0}\in E_{\phi}(\Omega, E)$ and that the integrand $f$ verifies $(\mathcal{S}_{\phi})$. Then
for every $\lambda>\lambda_{0}$:
$$\sup_{e^{*} \in \partial^{C} f(x_{0}+e)}\phi^{*}(\lambda e^{*})\leq c \phi(\beta(x_{0}+ e))+u_{\lambda}\leq \frac{1}{2}c\phi(2\beta e)+\frac{1}{2}c\phi(2\beta x_{0})+u_{\lambda}\:.$$
Since  $x_{0}\in E_{\phi}(\Omega, E)\;$,  then $\phi(2\beta x_{0})$ is integrable, this proves that the condition $(\mathcal{S}_{\phi, x_{0}})$ of Proposition \ref{propw.2} is true at $x_{0}$. $\;\Box$\\\\
For more clarity and simplicity, in the sequel of this section $\mathcal{X}=L_{p}(\Omega, E)$, $p^{-1}+q^{-1}=1$. Recall that when $1<p<\infty$ and $\phi=p^{-1}\Vert .\Vert^{p}$, then $\phi^{*}=q^{-1}\Vert .\Vert^{q}$, when $\phi=\Vert .\Vert$, then $\phi^{*}=\iota_{B^{*}}$ where $B^{*}$ is the unit ball of $E^{*}$. We will made the following assumptions on $f$:\\
 {\it $(\mathcal{S}_{p})$:  $1\leq p<\infty$ and there exist a positive constant $c$ and a $q$-integrable non negative function $a$ such that  for every $e\in E$, and $e^{*} \in \partial^{C} f_{\omega}( e)$, 
 $$\Vert e^{*}\Vert_{*}\leq c\Vert e\Vert^{p-1} +a\;.$$
  $(\mathcal{S_{\infty}})$ There exist $\eta>0$ and an integrable function $c$ such that: 
$$\Vert e\Vert\leq \eta\Rightarrow  f_{\omega}(x_{0}(\omega)+e)-f_{\omega}(x_{0}(\omega)) \geq -c(\omega)\Vert e\Vert\;.$$}
\begin{coro}\label{thm8990} Suppose $1< p<\infty$ and  $\phi=p^{-1}\Vert .\Vert^{p}$. Let a measurable integrand $f:\Omega\times E\to {\R}$. The following assertions are equivalent:\\
$(a)$ for every point $x_{0}\in L_{p}(\Omega, E)$ the integrand $f$ verifies $(\mathcal{S}_{\phi, x_{0}})$, \\
$(b)$ there exists a point $x_{0}\in L_{p}(\Omega, E)$  such that $f$ verifies $(\mathcal{S}_{\phi, x_{0}})$, \\
$(c)$ the integrand $f$ verifies  $(\mathcal{S}_{p})$.\\ 
Moreover when they hold, the integrand $f$ has the Fr\'echet-dlcp at every $x_{0}\in L_{p}(\Omega, E)$.
\end{coro}
Proof of Corollary \ref{thm8990}. First let us proves that $(b)\Rightarrow (c)$.  Suppose that $f$ verifies $(\mathcal{S}_{\phi, x_{0}})$ holds with $\phi=p^{-1}\Vert .\Vert^{p}$ and the point  $x_{0}$ is in $L_{p}(\Omega, E)$. For almost every $\omega\in \Omega$, the function $f_{\omega}$ is Lipschitzian on every ball of $E$, moreover keeping $\lambda=\lambda_0$ in the definition of $(\mathcal{S}_{\phi, x_{0}})$, we obtain the existence of a positive constant $c$ and an element $u$ of $L_{1}(\Omega,{\R})$ such that  for every $e\in E$, and $e^{*} \in \partial^{C} f_{\omega}( x_{0}+e)$, 
$$\Vert e^{*}\Vert^{q}_{*}\leq c\Vert e\Vert^{p} +u\leq 2^{p-1}c(\Vert x_{0}+e\Vert^{p}+\Vert x_{0}\Vert^{p})+u \;.$$
Therefore there exist a positive constant $d$, an integrable function $v$ such that
$$\Vert e^{*}\Vert^{q}_{*}\leq d\Vert x_{0}+e\Vert^{p}+v$$
since for every $0<s<1$, $t\to t^{s}$ is subadditive, with $s=q^{-1}$, we get for every $e^{*} \in \partial^{C} f_{\omega}( x_{0}+e)$,
$$\Vert e^{*}\Vert_{*}\leq d^{s}\Vert x_{0}+e\Vert^{ps} +v^{s}\;,$$
but $ps=p-1$, the function $v^{s}$ being $q$-integrable we obtain for every $e^{*} \in \partial^{C} f_{\omega}(e)$,
 $$\Vert e^{*}\Vert_{*}\leq d^{s}\Vert e\Vert^{p-1} +v^{s}\;.$$
 This proves that the integrand $f$ verifies  $(\mathcal{S}_{p})$ and proves that $(b)\Rightarrow (c)$. Now suppose that  the integrand $f$ verifies  $(\mathcal{S}_{p})$. Then there exists a positive constant $c$ a non negative $q$-integrable function $a$ such for $e^{*} \in \partial^{C} f_{\omega}(e)$:
 $$\Vert e^{*}\Vert_{*}\leq c\Vert e\Vert^{p-1} +a\;,$$
 therefore since $\Vert e+ e^{'}\Vert^{q}\leq  2^{q-1}(\Vert e\Vert^{q} +\Vert e\Vert^{q})$, we deduce for  $e^{*} \in \partial^{C} f_{\omega}(e)$:
 $$\Vert e^{*}\Vert_{*}^{q}= (c\Vert e\Vert^{p-1} +a)^{q}\leq 2^{q-1}(c^{q}(\Vert e\Vert^{q(p-1)} +a^{q})=2^{q-1}c^{q}\Vert e\Vert^{p)}+2^{q-1}a^{q}\;,$$
 thus we obtain the existence of a positive real number $d>0$ and an integrable function $u$ such for every  $e^{*} \in \partial^{C} f_{\omega}(e)$:
$$\Vert e^{*}\Vert_{*}^{q}\leq d\Vert e\Vert^{p}+u\;.$$
Then for every point  $x_{0}$ is in $L_{p}(\Omega, E)$, for every  $e^{*} \in \partial^{C} f_{\omega}(e)$:
$$\Vert e^{*}\Vert_{*}^{q}\leq d2^{p-1}(\Vert x_{0}+e\Vert^{p} +\Vert x_{0}\Vert^{p})+u)\:,$$
but for every $\lambda>0$ for  every $e^{*} \in \partial^{C} f_{\omega}(e)$, with the above inequality we reach:
$$\Vert\lambda e^{*}\Vert_{*}^{q}\leq d2^{p-1}\Vert \lambda^{qp^{-1}}e\Vert^{p}+\lambda^{q}(d2^{p-1}\Vert x_{0}\Vert^{p}+u)\;.$$
This proves that the integrand $f$ verifies $(\mathcal{S}_{\phi,  x_{0}})$ at any point  $x_{0}$ in $L_{p}(\Omega, E)$, with $\phi=p^{-1}\Vert .\Vert^{p}$, and $(c)\Rightarrow (a)$. 
The proof of the last assertion is an immediate consequence of Proposition  \ref{propw.2},$\Box $\\\\
Let us consider the case $p=\infty$.
\begin{coro}\label{coro5050} If the measurable integrand $f$ satisfies the condition $(\mathcal{S_{\infty}})$ 
then the integrand $f$ has the Fr\'echet-dlcp on  $L_{\infty}(\Omega,E)$ at the point $x_{0}\in L_{\infty}(\Omega,E)$.
\end{coro}
Proof of Corollary \ref{coro5050}. Given $\lambda>0$, and $r>0$ such 
$\lambda^{-1} r<\eta$, then:
$$\Vert e\Vert\leq \lambda^{-1}\Rightarrow \Vert e\Vert r\leq \lambda^{-1} r<\eta,\;\;\mbox{and:}\;\; [f]( x_{0}, e,r)\geq -\lambda^{-1} c$$
or equivalently the property of Proposition \ref{propw.200} holds with $\phi=\iota_{B}$, for every $\epsilon>0$ and $\lambda>0$ we have eventually:
$$[f]( x_{0}, e,r)\geq -\epsilon\phi_{\lambda}(e) -\lambda^{-1} c\;.\;\Box $$\\
Hereafter $L_{p}(\Omega, E)$ is endowed with  the  topology $\sigma_{p}=\sigma(L_{p}(\Omega,E) L_{q}(\Omega,E^{*}))$ $p^{-1}+q^{-1}=1$. 
\begin{prop}\label{prop100} Let $1< p\leq\infty$, the Banach $E$ being reflexive. When $p=\infty$, suppose $L_{1}(\Omega, E^{*})$ is separable. Let $x_{0}\in L_{p}(\Omega, E)$ with $f(x_{0})$ integrable and let us consider the following assertions:\\
$(a)$ the integrand $f$ has the $\sigma_p$-dlcp at $x_0$,\\
$(b)$ the integrand $f$ has the Fr\'echet-dlcp at $x_{0}$\\
$(c)$ for every sequence $(r_{n})_n$ of positive real numbers converging to $0$ the sequence of differential quotients $([f](x_{0}, ., r_{n}))_n$ satisfies the $\sigma_p$-Ioffe's criterion (see Definition \ref{def6.3}) at every point of $L_{p}(\Omega, E)$.\\
Then $(a)\Rightarrow (b)\Rightarrow(c)$. Moreover if the measure $\mu$ is atomless all the assertions are equivalent.
\end{prop}
Proof of Proposition \ref{prop100}. First remark if $1<p<\infty$ since $E$ is reflexive, then $L_{p}(\Omega, E)$ is reflexive too. Therefore the bounded sets of $L_{p}(\Omega, E)$ are the relatively sequentially $\sigma_p$-compacts sets. When $p=\infty$, since $L_{1}(\Omega, E^{*})$ is separable, the bounded sets of $L_{\infty}(\Omega, E)$ are the relatively sequentially $\sigma_\infty$-compacts sets. Therefore from Proposition \ref{prop379}, $(a)\Leftrightarrow (b)$. Then always $(a)\Leftrightarrow (b)\Rightarrow(c)$. Conversely, let us prove that $(c)\Rightarrow (b)$. If $(b)$ fails there exist  a sequence $(r_{n})_n$
of positive real numbers converging to $0$, and a bounded sequence  $(x_{n})_n$ such the sequence  $([f]^{-}(x_{0}, x_{n}, r_{n}))_n$ is not uniformly integrable. Since the measure
$\mu$ is atomless from Theorem \ref{prop991} $(b)$, due to the definition of the index of equi-integrability Definition \ref{defw6.102}, there exists a decreasing sequence $(A_{k})_k$ of measurable sets with a negligible intersection and a  $\sigma_p$-converging subsequence  $(x_{n_{k}})_k$ satisfying:
$$\lim_{k}\int_{A_{k}}[f]^{-}(x_{0}, x_{n_{k}}, r_{n_{k}}) d\mu>\epsilon>0\;.$$
Define $B_{k}=A_{k}\cap\{[f](x_{0}, x_{n_{k}}, r_{n_{k}})\leq 0\}$. Then since $[f](x_{0}, 0, r_{n})=0$,
$$[f]^{-}(x_{0}, x_{n_{k}}, r_{n_{k}})1_{A_{k}}=-[f](x_{0}, x_{n_{k}}, r_{n_{k}})1_{B_{k}}=-[f](x_{0}, x_{n_{k}}1_{B_{k}}, r_{n_{k}})\;,\;\mbox{and}\;$$
$$\int_{A_{k}}[f]^{-}(x_{0}, x_{n_{k}}1_{B_{k}}, r_{n_{k}}) d\mu=\int_{A_{k}}[f]^{-}(x_{0}, x_{n_{k}}, r_{n_{k}}) d\mu\;.$$
But  $1< p\leq\infty$ thus $1\leq q<\infty$, and since $(A_{k})_k$ has a negligible intersection, the sequence $( x_{n_{k}}1_{B_{k}})_k$ $\sigma_p$-converges to the origin and verifies for every integer $k$:
$$I_{[f](x_{0},  ., r_{n_{k}})}(x_{n_{k}}1_{B_{k}})\leq 0\;\;,\;\mbox{and}\;$$
$$\lim_{k}\int_{A_{k}} [f]^{-}(x_{0}, x_{n_{k}}1_{B_{k}}, r_{n_{k}})=\lim_{k}\int_{A_{k}}[f]^{-}(x_{0}, x_{n_{k}}, r_{n_{k}}) d\mu>\epsilon>0\;.$$
This proves with Theorem \ref{prop991} $(b)$, that the $\sigma_p$-Ioffe's criterion for the sequence  $([f](x_{0}, ., r_{n_{k}}))_n$ fails at $x=0$. 
We have proved the desired equivalences.$\;\Box$

\section{Solid and integral Bornologies, subdifferentials}
In this short section, we provide statements and proofs valid in a more general setting than the Lebesgue spaces.
Let $(\mathcal{X}, \Vert .\Vert)$ be a normed space. Recall that the differential quotient $[f]$ of $f\in \overline{{\R}}^{\mathcal{X}}$ at a point $x_{0}$ of its domain is defined by the formula
$$[f](x_{0}, x, r)=r^{-1}(f(x_{0}+rx)-f(x_{0}))\;.$$
Observe that for all $x^{*}\in \mathcal{X}^{*}$ we have: $[f-\langle x^{*}, .\rangle]=[f]-\langle x^{*}, .\rangle$.
Following the presentation given in \cite{110} section 4.1.6, for a function $f\in {\R}^{\mathcal{X}}$ finite at $x_0$, the subdifferential associated to a  bornology $\mathcal{B}$ on $\mathcal{X}$  is the set $\partial^{\mathcal{B}} f(x_{0})$ of $x^{*}\in \mathcal{X}^{*}$ such that for all $X\in \mathcal{B}$ one has: 
$$\liminf_{r\to 0_{+}} \inf_{x\in X}([f](x_{0}, x, r)-\langle x^{*}, x\rangle)\geq 0\;,$$
or equivalently for all $X\in \mathcal{B}$:
$$\limsup_{r\to 0_{+}} \sup_{x\in X}([f](x_{0}, x, r)-\langle x^{*}, x\rangle)^{-}= 0\;.$$
The elements of  $\partial^{\mathcal{B}} f(x_{0})$ are called the $\mathcal{B}$-subderivatives of $f$ at $x_0$.\\
If $\mathcal{B}$ is the bounded sets of $\mathcal{X}$ then $\partial^{\mathcal{B}} f(x_{0})$ is the Fr\'echet subdifferential $\partial^{F} f(x_{0})$ of $f$ at $x_0$. Classically an equivalent definition is the following:
$$\partial^{F} f(x_{0})=\{x^{*}\in \mathcal{X}^{*}:\; \liminf_{\substack{x\to x_{0}}}\frac{f(x)-f(x_{0})-\langle x^{*},x-x_{0}\rangle}{\Vert x-x_{0}\Vert}\geq 0\;\}.$$
When $\partial^{F} f(x_{0})$ is non empty, the function $f$ is said to be Fr\'echet subdifferentiable at the point $x_{0}\in X$ of its domain.\\
If $\mathcal{T}$ is a topology on $\mathcal{X}$ finer or equal than the weak topology of $\mathcal{X}$ and $\mathcal{B}$ is the bornology associated to the sequentially relatively compact sets,  then $\partial^{\mathcal{B}} f(x_{0})$ is the Hadamard subdifferential $\partial^{\mathcal{T}} f(x_{0})$ of $f$ at $x_0$. One can see that an equivalent definition is the following:
$$\partial^{\mathcal{T}} f(x_{0})=\{x^{*}\in{\mathcal{X}^*}:\forall x\in{X}, \langle x^{*}, x\rangle\leq f^{\mathcal{T}}(x_{0};x)\}\;,$$
where $f^{\mathcal{T}}(x_{0};.)$ is the sequential Hadamard directional subderivate of $f$ at $x_0$ defined by:
$$ f^{\mathcal{T}}(x_{0};x)= \displaystyle \inf_{(x_{n})\xrightarrow{\mathcal{T}}x,\; (r_{n}) \to 0_{+}}\liminf_{n} [f](x_{0},x_{n}, r_{n})=\displaystyle \inf_{ (r_{n}) \to 0_{+}} seq \mathcal{T}-li_{e} [f](x_{0},., r_{n})(x) .\;$$
In this section, $E$ is a separable Banach space and $\mathcal{X}$ will be a decomposable normed subspace of $L_{0}(\Omega, E)$ and sometimes the following (Rockafellar's decomposability) property is used:\\\\
{\it $(\mathcal{D})$ There exists an increasing covering $(\Omega_{n})_n$ of $\Omega$ such that for every integer $n$,\\
\centerline{$L_{\infty}(\Omega,E)1_{\Omega_{n}}=\{y1_{\Omega_{n}}, y\in L_{\infty}(\Omega,E)\}\subset\mathcal{X}$.}}
\begin{Def} A bornology $\mathcal{B}$ on $\mathcal{X}$ is said solid if for every element $X\in \mathcal{B}$, the subset $\mathbb{T}(X)=\{x1_{A}:\;x\in X,\; A\in \mathbb{T}\}$ is an element of $\mathcal{B}$.
\end{Def}
Since for every $A\in \mathbb{T}$ and $x\in L_{p}(\Omega,E)$ we have $\Vert x1_{A}\Vert_{p}\leq \Vert x\Vert_{p}$, we deduce that the Fr\'echet bornology on a Lebesgue space $L_{p}(\Omega,E)$ is solid. More generally, if $\phi: \Omega\times E \to \overline{\R}_+$ is a Young integrand, the associated Orlicz space ${\R}_{+} I_{\phi}^{\leq 1}=L_{\phi}(\Omega, E, \mu )$ generated by the sublevel set $I_{\phi}^{\leq 1}$ becomes a decomposable normed space endowed with
the Minkowski gauge associated to the sublevel set $I_{\phi}^{\leq 1}$, $\Vert x\Vert_{\phi}=\inf\{t>0: x\in t.I_{\phi}^{\leq 1}\}$. Moreover for every $A\in \mathbb{T}$ and $x\in L_{\phi}(\Omega,E)$ we have $\Vert x1_{A}\Vert_{\phi}\leq \Vert x\Vert_{\phi}$; therefore the Fr\'echet bornology on an Orlicz space is solid. In addition, it verifies always the property $(\mathcal{D})$ \cite{307}. 
When $E$ is reflexive, due to the weak compactness Dunfort-Pettis's criterion , the bornology on $L_{1}(\Omega,E)$ associated to the weakly compact sets is solid. 
\begin{thm} \label{prop 020} Given a  decomposable normed subspace  $\mathcal{X}$ of $L_{0}(\Omega, E)$ and $f$  a measurable integrand, if $\mathcal{B}$ is a solid bornology on  $\mathcal{X}$, and $x^{*}\in  \mathcal{X}^{*}\cap L_{0}(\Omega,E^{*}_{\sigma^{*}})$ then the following assertions are equivalent:\\
$(a)$ $x^{*}$ is a  $\mathcal{B}$-subderivative of $I_f$ at $x_0$.\\
$(b)$ for every element $B\in \mathcal{B}$, $\lim_{r\to 0_{+}}\sup_{x\in B}\int_{\Omega}[f-\langle x^{*},  .\rangle]^{-}(x_{0}, x, r)=0$.\\
\end{thm}
Proof of Theorem \ref{prop 020}.  Keeping the integrand $f-\langle x^{*}, . \rangle$ it suffices to consider the case $x^{*}=0$. 
If $(b)$ fails, there exist $\epsilon >0$, a sequence  $(r_{n})_n$ of positive real numbers converging to $0$, an element $X\in\mathcal{B}$ a sequence $(x_{n})_n$ in $X$, such eventually $\int_{\Omega}[f]^{-}(x_{0}, x_{n}, r_{n}) d\mu\geq \epsilon >0$. 
For each integer $n$ there exists a measurable set $A_{n}=\{[f](x_{0}, x_{n}, r_{n})\leq 0\}$ such that:
$$\int_{\Omega} [f]^{-}(x_{0}, x_{n}, r_{n}) d\mu=-\int_{A_{n}} [f](x_{0}, x_{n}, r_{n}) d\mu\;.$$
But since $[f](x_{0}, 0, r_{n})=0$,
$$[f](x_{0}, x_{n}, r_{n})1_{A_{n}}=[f](x_{0}, x_{n}1_{A_{n}}, r_{n})\;,\;$$
we get eventually:
$$\int_{\Omega}[f](x_{0}, x_{n}1_{A_{n}}, r_{n}) d\mu=-\int_{\Omega}[f]^{-}(x_{0}, x_{n}, r_{n}) d\mu\leq -\epsilon\;.$$
The sequence $(x_{n})_n$ being in $X$, since $\mathcal{B}$ is solid, the sequence $(x_{n}1_{A_{n}})_n$ is in $\mathbb{T}(X)\in \mathcal{B}$. Thus the above inequality shows that there exists an element $Y=\mathbb{T}(X)\in\mathcal{B}$ with
$$\liminf_{r\to 0_{+}} \inf_{x\in Y} [I_{f}](x_{0}, x, r)\leq -\epsilon\;.$$
Therefore $0$ is not a $\mathcal{B}$-subderivative of $I_f$ at $x_0$. We have proved that $(a)\Rightarrow(b)$. \\
Conversely suppose $(b)$ holds. Given an element $X\in \mathcal{B}$ we get:
$$\liminf_{r\to 0_{+}} \inf_{x\in X}([I_{f}-\langle x^{*}, .\rangle](x_{0}, x, r)\geq -\limsup_{r\to 0_{+}}\sup_{x\in X} \int_{\Omega}[f-\langle x^{*}, .\rangle]^{-}(x_{0}, x, r)  d\mu=0\;.$$
This proves that $x^{*}$ is a  $\mathcal{B}$-subderivative of $I_f$ at $x_0$ and $(b)\Rightarrow(a)$. $\Box$
\begin{Def} An integral bornology $\mathcal{B}$ on $\mathcal{X}$ is a bornology such there exists a family $\mathbb{F}=\{\phi_{i}, i\in I\}$ of non negative integrands verifying $\phi_{i}(0)=0$ for all $i\in I$ and such that $\mathcal{B}$ is the family of sets $B$ contained in a sublevel set of an element of the family of integral functionals $\{I_{\phi_{i}}, i\in I\}$.
\end{Def}
Notice that any integral bornology is solid. The Fr\'echet bornology of  $L_{p}(\Omega,E)$, $1\leq p<\infty$ is an integral bornology associated to the singleton $\mathbb{F}=\{\Vert .\Vert^{p}
\}$. The Fr\'echet bornology of  $L_{\infty}(\Omega,E)$, is an integral bornology associated to the family $\mathbb{F}=\{I_{\phi_{\lambda}}, \;\lambda>0\;\}$ where $\phi_{\lambda}(\omega,e)=\iota_{B_{E}}(\lambda e)=\iota_{\lambda^{-1}B_{E}}(e)$ and $B_{E}$ is the unit ball of $E$.  More generally, considering a Young integrand $\phi$ and the associated Orlicz space $L_{\phi}(\Omega, E)$ then the  Fr\'echet bornology of  $L_{\phi}(\Omega,E)$ is the family of sets $B$ contained in a sublevel set of an element of the family of integral functionals $\{I_{\phi_{\lambda}}, \;\lambda>0\}$ where $\phi_{\lambda}(\omega,e)=\phi(\omega, \lambda e)$. 
Since for every weakly compact set $X$ of  $L_{1}(\Omega,E)$, due to the Dunfort-Pettis and de la Vall\'ee Poussin criterions  Theorem \ref{def33.2} $(b)$, there exists an integrand $\phi$ of $\alpha$-Young type such that $\sup_{x\in X} I_{\phi}(x)\leq 1$, the bornology associated to the weakly relatively compact sets of $L_{1}(\Omega,E)$ called the weak Hadamard bornology is the integral bornology on  $L_{1}(\Omega,E)$ associated to the family of  $\alpha$-Young integrands.
\begin{thm} \label{prop 021} Let a decomposable normed subspace $\mathcal{X}$ of $L_{0}(\Omega, E)$ verifying in addition the decomposability property  $(\mathcal{D})$, and let $f$ be a measurable integrand. If $\mathcal{B}$ is an integral bornology on  $\mathcal{X}$ associated to the family $\mathbb{F}=\{\phi_{i}, i\in I\}$, let us consider the following assertions:\\
$(a)$ $x^{*}\in \mathcal{X}^{*}\cap L_{0}(\Omega,E^{*}_{\sigma^{*}})$ is a  $\mathcal{B}$-subderivative of $I_f$ at $x_0$.\\
$(b)$ for every $i\in I$, for every $\epsilon >0$, there exists a family of non negative eventually integrable functions $\{u_{r}, r\in (0,1)\}$ norm converging to $0$ in $L_{1}(\Omega,{\R})$ and verifying eventually
$$[f-\langle x^{*},  .\rangle](x_{0}, ., r)\geq -\epsilon\phi_{i} -u_{r}\;.$$ 
Then always $(b)\Rightarrow (a)$. If in addition the measure is atomless then these assertions are equivalent.
\end{thm}
Proof of Theorem \ref{prop 021}. Proof of the sufficiency part. Let $X\in \mathcal{B}$, there exists $(i, s)\in I\times {\R}_+$ such $\sup_{x\in X} I_{\phi_{i}}(x)\leq s$. Using condition $(b)$, for every $\epsilon >0$, there exists a family of non negative eventually integrable functions $\{u_{r}, r\in (0,1)\}$ norm converging to $0$ and verifying eventually for 
 $x\in X$: $[f-\langle x^{*},  .\rangle]^{-}(x_{0}, x, r)\leq \epsilon\phi_{i}(x)+u_{r}$, therefore for every $\epsilon>0$:
 $$\limsup_{r\to 0_{+}}\sup_{x\in X}\int_{\Omega}[f-\langle x^{*},  .\rangle]^{-}(x_{0}, x, r) d\mu\leq \epsilon.s\;,$$
thus
$$\lim_{r\to 0_{+}}\sup_{x\in X}\int_{\Omega}[f-\langle x^{*},  .\rangle]^{-}(x_{0}, x, r) d\mu=0\;,$$
and the sufficency part is proved with Theorem \ref{prop 020} $(b)$. Conversely suppose that $x^{*}\in \mathcal{X}^{*}\cap L_{0}(\Omega,E^{*}_{\sigma^{*}})$ is a  $\mathcal{B}$-subderivative of $I_f$ at $x_0$.
Then for every $i\in I$:
$$\limsup_{r\to 0_{+}}\sup_{I_{\phi_{i}}(x)\leq 1; x\in \mathcal{X}}\int_{\Omega}[f-\langle x^{*},  .\rangle]^{-}(x_{0}, x, r) d\mu=0\;.$$
Fix $i\in I$, and denotes $\phi_{i}$ by $\phi$. Setting $g_{r}(\omega,e)=[f_{\omega}-\langle x^{*}(\omega), e\rangle]^{-}(x_{0}(\omega), e, r)$, we have:
$$\mbox{for every}\;\epsilon>0\;\mbox{there exists} \;r_{\epsilon}>0 :\;\sup_{0<r\leq r_{\epsilon}} \sup_{I_{\phi}(x)\leq 1; x\in \mathcal{X}} I_{g_{r}}(x)\leq \epsilon.$$
The assumptions of  \cite{87} Corollary 5.7 with the constant multifunction $M(\omega)=E$ are satisfied. Indeed since $(\mathcal{D})$ holds $\mathcal{X}$ is rich (in sense of  \cite{87} Definition 3.7) in $S_{M}=L_{0}(\Omega,E)$. Using \cite{87} Corollary 5.7, for every $0<r<r_{\epsilon}$ we get the existence of a multiplier $y^{*}_{r} \geq 0$ such that the function
$$ v_{r}(\omega)=\inf\{-g_{r}(\omega,e)+y^{*}_{r} \phi(\omega,e), e\in E\}$$
is integrable and verifies
$$-\epsilon + y^{*}_{r} \leq \int_{\Omega} v_{r} d\mu\;.$$
Moreover since $\phi(0)=0$ and $g_{r}(0)=0$ the function $ v_{r}$ is non positive and we get for every $0<r<r_{\epsilon}$ that $y^{*}_{r} \leq \epsilon$ and that $\int_{\Omega} -v_{r} d\mu \leq \epsilon -y^{*}_{r}\leq \epsilon$. 
Define
$$u_{r, \epsilon}(\omega)=\inf\{-g_{r}(\omega,e)+\epsilon \phi(\omega,e), e\in E\geq v_{r}$$
then
$$g_{r}\leq -v_{r}+y^{*}_{r}\phi\leq -v_{r}+\epsilon\phi\leq -u_{r, \epsilon}+\epsilon\phi\;.$$
 We get: for  $0<r<r_{\epsilon}$, $u_{r, \epsilon}\geq 0$, $g_{r}\leq \epsilon\phi_{i} +u_{r, \epsilon}$ and $\int_{\Omega} -u_{r, \epsilon} d\mu\leq \epsilon$. Then $\lim_{r\to 0_{+}} \int_{\Omega} u_{r, \epsilon} d\mu=0$. Indeed take $\epsilon>0$.
 Given $0<\epsilon^{'}<\epsilon$, for $0<r<r_{\epsilon^{'}}$, we have $ u_{r, \epsilon^{'}}\leq u_{r,\epsilon}$ and  
 $$\int_{\Omega}-u_{r, \epsilon} d\mu\leq \int_{\Omega}- u_{r, \epsilon^{'}} d\mu\leq \epsilon^{'}.$$
 This  proves that $\lim_{r\to 0} \int_{\Omega} u_{r, \epsilon} d\mu=0$ and moreover $g_{r}\leq \epsilon\phi -u_{r, \epsilon}$     $\;\Box$\\\\
Taking the Fr\'echet bornology on an Orlicz space, we obtain immediatly the following criterion of subdifferentiability: 
\begin{coro} \label{coro613} Let $\phi$ be a Young integrand and  $f: \Omega\times E\to \overline{{\R}}$ be a measurable integrand. Consider a function $x^{*}\in  L_{\phi^{*}}(\Omega,E^{*}_{\sigma^{*}})$ and the following assertions:\\
$(a)$ $x^{*}$ is a  Fr\'echet subderivative of $I_f$ at $x_0$ on $L_{\phi}(\Omega, E)$,\\
$(b)$ for every $0<\epsilon$, $0<\lambda$, there exists a family of non negative eventually integrable functions $(u_{r})_{r>0}$ such that $\lim_{r\to 0}\Vert u_{r} \Vert_{1}=0$ and eventually: 
$$[f-\langle x^{*}, .\rangle](x_{0}, ., r)\geq -\epsilon\phi_{\lambda}(.)-u_{r}\;.$$ 
Then always $(b)\Rightarrow (a)$. If in addition the measure is atomless then these assertions are equivalent.
\end{coro}
The following particular consequence in case $1<p<\infty$ is nothing else than a rephrase of the J. P Penot's characterization \cite{216} Theorem 22. But not only, since  it contains also the  N. H. Chieu's and J. P Penot's characterization \cite{707} Theorem and \cite{216} Theorem 12 of the Fr\'echet subdifferential in case $p=1$ (see Corollary \ref{coro523}).
\begin{coro} \label{coro713} Suppose $1\leq p<\infty$ and $f: \Omega\times E\to \overline{{\R}}$ be a measurable integrand. Given $x^{*}\in  L_{q}(\Omega,E^{*}_{\sigma^{*}})$ let us consider the following assertions:\\
$(a)$ $x^{*}$ is a Fr\'echet subderivative of $I_f$ at $x_0$ on $L_{p}(\Omega, E)$,\\ 
$(b)$ For every $\epsilon>0$ there exists a family of non negative eventually integrable functions $(u_{r})_{r>0}$ such that $\lim_{r\to 0}\Vert u_{r} \Vert_{1}=0$ and eventually: 
$$[f-\langle x^{*}, .\rangle](x_{0}, ., r)\geq -\epsilon\Vert .\Vert^{p}-u_{r}\;.$$
Then always $(b)\Rightarrow (a)$ and if addition the measure is atomless these assertions are equivalent.
\end{coro}
Proof. The relations between the two assertions are an immediate consequence of Corollary \ref{coro613} with $\phi=p^{-1}\Vert .\Vert^{p}$. $\Box$
\begin{coro} \label{coro523} Let $p=1$,  $f: \Omega\times E\to \overline{{\R}}$ be a measurable integrand and suppose the measure is atomless.  Given $x^{*}\in  L_{\infty}(\Omega,E^{*}_{\sigma^{*}})$ the following assertions are equivalent:\\
$(a)$ $x^{*}$ is a Fr\'echet subderivative of $I_f$ at $x_0$ on $L_{1}(\Omega, E)$,\\
$(b)$ The integrand $f$ verifies the assertion $(b)$ of Corollary \ref{coro613} with $p=1$.\\
$(c)$ $x^{*}$ is a Moreau-Rockafellar subderivative of $I_f$ at $x_0$ on $L_{1}(\Omega, E)$.
\end{coro}
Proof. First remark that it suffices to prove this Corollary in case $x^{*}=0$. $(a)\Rightarrow (b)$ is a consequence of Corollary \ref{coro613} with $\phi=\Vert .\Vert$. Now suppose that condition $(b)$ of Corollary \ref{coro613} holds with  $x^{*}=0$. Fix $\epsilon>0$. For every $r>0$, define $v_{r}=\inf_{e\in E} [f](x_{0}, e, r)-\epsilon\Vert e\Vert$. Then  $v_{r}$ satisfies $-u_{r}\leq v_{r}\leq 0$, therefore the family $(v_{r})_{r>0}$ of non positive eventually integrable functions norm converges to the origin of $L_{1}(\Omega, {\R})$ when $r\to 0_+$. But
$$v_{r}=\inf_{e\in E} [f](x_{0}, e, r)]-\epsilon\Vert e\Vert=\inf_{e\in E} r^{-1}([](x_{0}, re)-\epsilon\Vert re\Vert)= r^{-1}v_{1}\;.$$
Since the  family $(v_{r})_{r>0}$ of eventually integrable functions norm converges to the origin of $L_{1}(\Omega, {\R})$ when $r\to 0_+$, neccesarily eventually $v_{1}=0 =v_{r}$.
Therefore the integrand $[f]$ verifies eventually at $x_{0}$:
$$\inf_{e\in E} [f](x_{0}, e, r)-\epsilon\Vert e\Vert\geq 0\;.$$
Or equivalently,
$$\inf_{e\in E} f(x_{0}+ e)-f(x_{0})-\epsilon\Vert e\Vert\geq 0\;.$$
This last property being valid  for every $\epsilon>0$,
one get for every $e\in E$,  $f(x_{0}+ e)-f(x_{0})\geq 0$, therefore $0$ is a Moreau-Rockafellar subderivative of $I_f$ at $x_0$ on $L_{1}(\Omega, E)$. Thus  $(b)\Rightarrow (c)$ is true.  $(c)\Rightarrow (a)$ is immediate. $\Box$\\\\
In the same spirit let us consider the case $p=\infty$. Using Corollary \ref{coro613} it can be obtained the following characteristic condition (when $f$ is only measurable) when the measure is atomless; let us give a direct proof without this assumption but when $f$ is normal.
\begin{coro} \label{coro73} Suppose $p=\infty$, and $f: \Omega\times E\to \overline{{\R}}$ is a normal integrand. An integrable function $x^{*}\in  L_{1}(\Omega,E^{*}_{\sigma^{*}}) $ is a Fr\'echet subderivative of $I_f$ at $x_0$ on  $L_{\infty}(\Omega,E)$, if and only if there exists a family of eventually non negative integrable functions $(u_{r})_{r>0}$ such that $\lim_{r\to 0_{+}}\Vert u_{r} \Vert_{1}=0$ and eventually: 
$$\inf_{\Vert e\Vert\leq 1}[f-\langle x^{*}, .\rangle](x_{0}, e, r)\geq -u_{r}\;.$$
\end{coro}
Proof of Corollary \ref{coro73}. Remark that it suffices to gives the proof in case $x^{*}=0$. Let $B_\infty$ be the unit ball of $L_{\infty}(\Omega,E)$. The condition is sufficient. Indeed: $$\lim_{r\to {+}}\sup_{x\in B_{\infty}}\int_{\Omega}[f]^{-}(x_{0}, x, r)d\mu\leq \lim_{r\to 0_{+}}\int_{\Omega} u_{r} d\mu=0\;,$$ 
then we obtain for $s>0$ $\lim_{r\to 0_{+}}\sup_{x\in sB_{\infty}}\int_{\Omega}[f]^{-}(x_{0}, x, r)d\mu=0$, this proves with Theorem \ref{prop 020} $(b)$ that $x^{*}$ is a Fr\'echet subderivative of $I_f$ at $x_0$ on  $L_{\infty}(\Omega,E)$. Conversely let us show the necessity. If $0$ is a Fr\'echet subderivative of $I_f$ at $x_0$ then with Theorem \ref{prop 020} $(b)$ we get 
$$\lim_{r\to 0_{+}}\sup_{x\in B_{\infty}}I_{[f]^{-}}(x_{0}, x, r)=0\;.$$ 
But for every $r>0$ if $u_{r}=\sup_{\Vert e\Vert\leq 1}[f]^{-}(x_{0}, e, r)$ then \cite{8} Theorem 2.2 asserts that
$$\sup_{x\in B_{\infty}} I_{[f]^{-}}(x_{0}, x, r)=\int_{\Omega} u_{r} d\mu\;.$$
Therefore $\lim_{r\to 0_{+}}\int_{\Omega} u_{r} d\mu=0$. 
The proof of Corollary \ref{coro73} is complete. $\Box$\\\\
Using Theorem \ref{prop 021} one obtain the following criterion for the weak Hadamard subdifferentiability:
\begin{coro} \label{coro822} Suppose $p=1$ and $f: \Omega\times E\to \overline{{\R}}$ is a measurable integrand. Given a function $x^{*}\in  L_{\infty}(\Omega,E^{*}_{\sigma^{*}})$ let us consider the following assertions:\\
$(a)$ $x^{*}$ is a weak Hadamard  subderivative of $I_f$ at $x_0$ on $L_{1}(\Omega, E)$:\\ 
$(b)$ for every (respectively there exists an) integrable positive valued function $\alpha$, for every $\alpha$-Young integrand $\phi$ and  for every
 $\epsilon>0$ there exists a family of eventually non negative integrable functions $(u_{r})_{r>0}$ such that $\lim_{r\to 0_{+}}\Vert u_{r} \Vert_{1}=0$ and eventually: 
$$[f-\langle x^{*}, .\rangle](x_{0}, ., r)\geq -\epsilon \phi(.)-u_{r}\;.$$
Then always $(b)\Rightarrow (a)$ and if the measure is atomless these assertions are equivalent.
\end{coro}
Proof. We have yet seen with  De la Vall\'e-Poussin's Theorem \ref{def33.2} that for every positive valued integrable function  $\alpha$, the weak Hadamard bornology is the integral bornology associated with the family of  $\alpha$-Young integrands and the result is a consequence of Theorem \ref{prop 021}. $\;\Box$

\section{Additional results on Fr\'echet subdifferentiability}
In this section, when $E$ is separable, the study of some properties of the Fr\'echet subdifferentiability is mainly made in relation with the  Fr\'echet differential compactness property. 
\begin{Def}\label{def120} Let an integrand $f: \Omega\times E\to \overline{{\R}}$, and  $x_{0}\in L_{0}(\Omega, E)$. The integrand $f$ is said Fr\'echet subdifferentiable   along $x_{0}$ when for almost every $\omega\in \Omega$, the function $f_\omega$ is Fr\'echet subdifferentiable at $x_{0}(\omega)$. An element  $x^{*}\in L_{0}(\Omega, E_{\sigma^{*}}^{*})$ is a Fr\'echet  subderivative of $f$ along $x_{0}$
when for almost every $\omega\in \Omega$,  $x^{*}(\omega)$ is a Fr\'echet subderivative of $f_\omega$ at $x_{0}(\omega)$.
\end{Def}
The following result gives a practical sufficient criterion for the  Fr\'echet subderivability of an integral functional.
\begin{thm} \label{thm4.8} Let us consider a decomposable normed space $(\mathcal{X}, \Vert .\Vert_{\mathcal{X}})$ topologically contained in some Lebesgue space $L_{p}(\Omega, E, \beta\mu)$ for some measurable positive valued function $\beta$. Suppose that the Fr\'echet bornology on  $\mathcal{X}$ is solid.  Let $f$ be a measurable integrand and the following assertions:\\
$(a)$ $x^{*}\in \mathcal{X}^{*}\cap L_{0}(\Omega,E^{*}_{\sigma^{*}})$ is a  Fr\'echet subderivative of $I_f$ at $x_0$ on $\mathcal{X}$.\\
$(b)$ The integrand $f-\langle x^{*}, \rangle$ has the Fr\'echet-dlcp at $x_0$.\\
Then $(a)\Rightarrow(b)$, if moreover $x^*$ is a Fr\'echet subderivative of $f$ along $x_0$, then $(b)\Rightarrow (a)$.
\end{thm}
Proof of Theorem \ref{thm4.8}.  Since the Fr\'echet bornology of $\mathcal{X}$ is solid then $(a)\Rightarrow(b)$ is a consequence of Theorem \ref{prop 020} in case where $\mathcal{B}$ is the Fr\'echet bornology. Conversely suppose $(b)$ holds and $x^*$ is a Fr\'echet  subderivative of $f$ along $x_{0}$. Setting $g(\omega,e)= f(\omega, x_{0}(\omega)+e)-f( x_{0}(\omega))-\langle x^{*}(\omega), e \rangle$, one may suppose that $x^{*}=0$ is a  Fr\'echet subderivative of $g$ along the origin,  thus $g(\omega,e)=\Vert e\Vert\epsilon(\omega,e)$ with $\displaystyle\liminf_{e\to 0}\epsilon_{\omega}(e)\geq 0$, and that $g$ has the  Fr\'echet-dlcp.
From Theorem \ref{prop 020} it suffices to prove that for every sequence  $(r_{n})_n$ of positive real numbers converging to $0$, for every bounded sequence  $(x_{n})_n$ in $\mathcal{X}$, the sequence $([g]^{-}(0, x_{n}, r_{n}))_n$ strongly converges to the origin in $L_{1}(\Omega, {\R})$.\\
The following Lemma is  proved in \cite{704}:
\begin{lem} \label{lem413} If $0$ is a Fr\'echet subderivative of $f$ along $x_0$, then the function $\epsilon(\omega,e)=\epsilon_{\omega}(e)$ if $e\neq 0$, $\epsilon(\omega,0)=0$, is measurable on $\Omega\times E$, and verifies for every $(\omega,e)\in \Omega\times E$: $f_{\omega}(x_{0}(\omega)+e)-f_{\omega}(x_{0}(\omega))=\Vert e\Vert\epsilon(\omega,e)\;.$
The integrand $\;{\epsilon}^{-}=-\min(\epsilon, 0)$ is measurable on $\Omega\times E$ and
moreover $\displaystyle\lim_{e\to 0}{\epsilon}^{-}(\omega,e)=0$.
\end{lem}
From the preceding Lemma $[g]^{-}(0, x_{n}, r_{n})= \Vert x_{n} \Vert\epsilon^{-}(r_{n}x_{n})$, moreover since the sequence
$(x_{n})_n$ is bounded  in $\mathcal{X}$ then by assumptions $(x_{n})_n$ it is bounded in some  $L_{p}(\Omega, E, \beta\mu)$ therefore the sequence $(r_{n}x_{n})_n$ norm converges to $0$ in $L_{p}(\Omega,E, \beta\mu)$, thus in $\beta\mu$- measure (\cite{2} section 4.7 (or \cite{38} Lemma 16.4)) and since $\beta$ is positive valued, extracting subsequences almost everywhere converging to $0$ we obtain that the convergence is in in local $\mu$-measure , so is the convergence of the sequence $(\epsilon^{-}(r_{n}x_{n})_n$. Let $v_{n}= \Vert x_{n} \Vert\epsilon(r_{n}x_{n})$. Since $f$ has the Fr\'echet-dlcp at $x_0$, the integrand $g$ too at the origin. Therefore the sequence $(v^{-}_{n})_n$ is uniformly integrable. We have $v_{n}^{-}=\Vert x_{n}\Vert {\epsilon}^{-}(r_{n}x_{n})$ applying Lemma \ref{890} with $E={\R}$, $y_{n}=\Vert x_{n}\Vert$, $u_{n}={\epsilon}^{-}(r_{n}x_{n})$, we deduce that $(v_{n}^{-})_n$ converges strongly to $0$ in $L_{1}(\Omega,{\R}, \mu)$.
The proof of Theorem \ref{thm4.8} is complete. $\;\Box$\\\\
As a consequence of \cite{306} Theorem 2.4 (vii), the statement of Theorem \ref{thm4.8} is valid on every Orlicz space $L_{\phi}(\Omega, E, \mu)$. Let us give a practical criterion for the  Fr\'echet-subdifferentiability. 
\begin{coro} \label{propw.199} Let $\phi$ be a Young integrand, and $f$ be a measurable integrand such that $x^{*}\in L_{\phi^{*}}(\Omega, E^{*})$ is a Fr\'echet-subderivative of $f$ along $x_{0}\in L_{\phi}(\Omega, E)$. Suppose that the following condition holds: \\
 for every $\epsilon >0$, and every $\lambda>0$, there exists a family of non negative eventually integrable functions $\{u_{r}, r\in (0,1)\}$ uniformly integrable in $L_{1}(\Omega,{\R})$ and verifying eventually
$$[f](x_{0}, ., r)-\langle x^{*}, .\rangle\geq -\epsilon\phi_{\lambda} -u_{r}\;.$$
Then $ x^{*}$ is a  Fr\'echet-subderivative of $I_f$ at $x_0$ on $ L_{\phi}(\Omega, E)$.
\end{coro}
Proof. due to Proposition \ref{propw.200}, $f-\langle x^{*}, .\rangle$ has the Fr\'echet differential lcp on $L_{\phi}(\Omega, E)$ at $x_{0}\in L_{\phi}(\Omega, E)\;.$ The result is then an immediate consequence  of Theorem \ref{thm4.8}. $\;\Box$\\\\
Given a Young function $\phi$ we will consider the following subspace of $L_{\phi}(\Omega, E, \mu)$:\\ 
\centerline{$E_{\phi}(\Omega, E)=\{x\in L_{\phi}(\Omega, E): \forall \lambda>0, \phi(\lambda x)\in L_{1}(\Omega,{\R}, \mu)\}$.}
Given an  $E^*$-valued multifunction $M$,  $L_{\phi^{*}}(M)_{\sigma^{*}}$, (respectively $E_{\phi^{*}}(M)_{\sigma^{*}})$ denotes the set of almost everywhere selections of $M$ which are in $L_{\phi^{*}}(\Omega,E^{*}_{\sigma^{*}})$ (respectively $E_{\phi^{*}}(\Omega,E^{*}_{\sigma^{*}}))$.
\begin{prop} \label{coro17} Let $\phi$ be a Young integrand and $f$ be an integrand  such $x^{*}\in E_{\phi^{*}}(\Omega, E^{*})$ is a Fr\'echet-subderivative of $f$ along $x_{0}\in L_{\phi}(\Omega, E)$. If the condition $(\mathcal{S}_{\phi,  x_{0}})$ of Theorem \ref{propw.2} holds then $E_{\phi^{*}}(\partial^{F}f(x_{0}))_{\sigma^{*}}\subset \partial^{F} I_{f}(x_{0})$.
\end{prop}
Proof.  First remark that for any  $x^{*}\in E_{\phi^{*}}(\Omega, E^{*})$ the integrand $f-\langle x^{*}, .\rangle$ has the  Fr\'echet-dlcp at $x_{0}$. Indeed let $f$ verifying $(\mathcal{S}_{\phi,  x_{0}})$. Since for every $x^{*}\in E_{\phi^{*}}(\Omega, E^{*})$,  due to the Young inequality, for every $e\in E$,  and every $\lambda>0$, $\vert\langle x^{*}, e\rangle\vert\leq \phi(\lambda e)+\phi^{*}(\lambda^{-1}x^{*})$ with $\phi^{*}(\lambda^{-1}x^{*})$ integrable, we obtain with Lemma \ref{lem 819}, that $f-\langle x^{*}, .\rangle$ verifies the condition of Corollary \ref{propw.199}. If moreover  $x^{*}\in E_{\phi^{*}}(\partial^{F}f(x_{0}))_{\sigma^{*}}$ is a Fr\'echet-subderivative of $f$ along $x_{0}$, then  due to Corollary \ref{propw.199}, $x^{*}$ is a  Fr\'echet-subderivative of $I_f$ at $x_0$ on $ L_{\phi}(\Omega, E)$. $\;\Box$
\begin{coro} \label{coro171} Let $\phi$ be a Young integrand, $f$ be a real valued integrand with a Fr\'echet-derivative $x^{*}\in E_{\phi^{*}}(\Omega, E^{*})$ along $x_{0}\in L_{\phi}(\Omega, E)$. If in addition the condition $(\mathcal{S}_{\phi, x_{0}})$ holds, then $I_f$ is Fr\'echet-differentiable at $x_{0}$ on $L_{\phi}(\Omega, E)$. Moreover $ I_{f}^{'}(x_{0})= f^{ '}(x_{0})$.
\end{coro}
Proof. Remark that if $f$ satisfies $(\mathcal{S}_{\phi}, x_{0})$, since for almost every $\omega\in \Omega$, the function $f_{\omega}$ is Lipschitzian on every ball of $E$, thus  $\partial^{C} f(e)=-\partial^{C} -f(e)$, and the Young integrands being even then $-f$ satisfies $(\mathcal{S}_{\phi,  x_{0}})$ too.  Moreover due to Corollary \ref{rem 713},  $I_f$ is finite on a ball of $L_{\phi}(\Omega, E)$ centered at $x_0$, thus locally $I_{f}(x)=-I_{-f}(x)$.
Proposition \ref{coro17} ensures that $f^{ '}(x_{0})\in \partial^{F} I_{f}(x_{0})\cap -\partial^{F} -I_{f}(x_{0})$ this proves that $I_f$ Fr\'echet-differentiable at $x_{0}\in E_{\phi}(\Omega, E)\;$ on $L_{\phi}(\Omega, E)$ and  $ I_{f}^{'}(x_{0})= f^{ '}(x_{0})$. $\;\Box$\\\\
For clarity now we will restrict ourselves to the case of Lebesgue spaces. The following result is a permanence property.
\begin{prop} \label{coro174} Let $1<p\leq \infty$. If the integrand $f$ has the Fr\'echet-dlcp at $x_{0}\in L_{p}(\Omega, E)$ and $x^{*}\in L_{q}(\Omega, E_{\sigma^{*}}^{*})$, then the function $f-\langle x^{*}, \rangle$ has the Fr\'echet-dlcp at $x_0$.
\end{prop} 
Proof of Proposition \ref{coro174}. For every measurable set $A$, and every $x\in L_{p}(\Omega, E)$, $x^{*}\in L_{q}(\Omega, E^{*})$, we have the upper bound: 
$$\int_{A}\vert \langle x^{*}, x\rangle\vert d\mu\leq \Vert x^{*}1_{A}\Vert_{q}\Vert x\Vert_{p}\;.$$
Since $1\leq q<\infty$, for every norm bounded sequence ${(x_{n})}_n$, the sequence $(\langle x^{*}, x_{n}\rangle)_n$ is uniformly integrable.
Moreover setting $g=f-\langle x^{*}, .\rangle$, we have: $[g]=[f]-\langle x^{*}, .\rangle$ and since $f$ has the  Fr\'echet-dlcp at $x_{0}\in X$, the integrand $g$ has the Fr\'echet-dlcp at $x_0$.$\;\Box$
\begin{prop}\label{prop 5.10}  Let $1\leq p\leq \infty$. Let $x_{0}\in L_{p}(\Omega, E) $ with $f(x_{0})$ integrable. 
Suppose that the integral functional $I_f$ is  Fr\'echet-subdifferentiable at $x_{0}$ with a  Fr\'echet-subderivative  $x^{*}_{0}\in L_{q}(\Omega, E^{*}_{\sigma^{*}})$. Then $f-\langle x^{*}_{0} , \rangle$ has the  Fr\'echet-differential lower compactness property. Moreover if $p\neq 1$, for every $x^{*}\in L_{q}(\Omega, E_{\sigma^{*}}^{*})$  the integrand $f-\langle x^{*}, \rangle$ has the Fr\'echet-dlcp at $x_0$.
\end{prop}
Proof of Proposition \ref{prop 5.10}. If $I_f$ is  Fr\'echet-subdifferentiable at $x_0$, with a Fr\'echet-subderivative $x_{0}^{*}\in L_{q}(\Omega, E_{\sigma^{*}}^{*})$ then Theorem \ref{thm4.8} asserts that $f-\langle x^{*}_{0} , \rangle$ has the  Fr\'echet-differential lower compactness property. Moreover if $p\neq 1$, due to Proposition \ref{coro174} for every $x^{*}\in L_{q}(\Omega, E_{\sigma^{*}}^{*})$  the integrand $f-\langle x^{*}, \rangle$ has the Fr\'echet-dlcp at $x_0$. $\;\Box$

\begin{thm} \label{corow11} Suppose $1< p \leq \infty$, if  the integrand $f$ has the Fr\'echet-dlcp at $x_{0}\in L_{p}(\Omega,E)$, then: 
$$L_{q}(\partial^{F}f(x_{0}))_{\sigma^{*}}\subseteq \partial^{F}I_{f}(x_{0})\;.$$
\end{thm}
Proof of Theorem \ref{corow11}.  From Proposition \ref{coro174}, since $p\neq 1$, for every $x^{*}\in L_{q}(\partial^{F}f(x_{0}))_{\sigma^{*}}$ the integrand $f-\langle x^{*}, \rangle$ has the Fr\'echet-dlcp at $x_{0}$ and Theorem \ref{thm4.8} ensures that $x^{*}$ is a Fr\'echet subderivative of $I_f$ at $x_0$. $\;\Box$ 
\begin{coro} \label{coroww11} Suppose $1< p\leq \infty$. If the integral functional $I_f$ is Fr\'echet subdifferentiable at $x_{0}\in L_{p}(\Omega,E)$ with a  Fr\'echet-subderivative  $x^{*}_{0}\in L_{q}(\Omega, E^{*}_{\sigma^{*}})$, then the integrand $f$ has the  Fr\'echet-dlcp at $x_0$ and:
$$L_{q}(\partial^{F}f(x_{0}))_{\sigma^{*}}\subseteq \partial^{F}I_{f}(x_{0})\;.$$
\end{coro}
Proof of Corollary \ref{coroww11}. Due to Proposition \ref{prop 5.10} the integrand $f$ has the  Fr\'echet-dlcp at $x_{0}\in L_{p}(\Omega,E)$, and the result is a consequence of Theorem \ref{corow11}.$\;\Box$ 
\begin{coro} \label{prop450} Suppose $1<p<\infty$. Given $x_{0}\in L_{p}(\Omega, E)$ and a measurable integrand $f:\Omega\times E\to {\R}$ satisfying condition $(\mathcal{S}_{p})$. Then $$L_{q}(\partial^{F}f (x_{0}))_{\sigma^{*}}\subset \partial^{F} I_{f}(x_{0})\;.$$ 
\end{coro}
Proof.  It is an immediate consequence of Proposition \ref{coro17}. $\;\Box$
\begin{coro} \label{prop451} Suppose $p=\infty$, given $x_{0}\in L_{\infty}(\Omega, E)$ and a measurable integrand $f:\Omega\times E\to \overline{{\R}}$ satisfying the condition $(\mathcal{S}_{\infty})$. Then $$L_{1}(\partial^{F}f (x_{0}))_{\sigma^{*}}\subset \partial^{F} I_{f}(x_{0})\;.$$ 
\end{coro}
Proof. From Corollary \ref{coro5050}  the integrand $f$ has the Fr\'echet-dlcp, and from Proposition \ref{coro174}, for every element  $y^{*}\in L_{1}(\partial^{F}f (x_{0}))_{\sigma^{*}}$, the integrand $f-\langle y^{*}, \rangle$ has the Fr\'echet-dlcp at the point $x_0$, and Theorem \ref{thm4.8}  allows to conclude. $\;\Box$

\section{More on weak Hadamard subdifferentiability.}
The author in \cite{49}, makes a first study of the strong Hadamard subdifferentiability of integral functionals on Lebesgue spaces, related with the properties of the differential quotients. In a initial version of this article, when $E$ is reflexive, all the results on Fr\'echet-subdifferentiability of the section 10 has been proved with the results of section 6 by considering the weak topology on $L_{p}(\Omega,E)$ and  the weak star topology when $p=\infty$.
On reflexive spaces the Fr\'echet bornology coincide with the weak Hadamard bornology therefore the  Fr\'echet subdifferential coincide with the weak Hadamard subdifferential.
Notice that this last result can be extended (the proof is omitted) to the case where $X$ has a predual $Y$:
\begin{lem} \label{propw11a} Suppose that $Y$ is a separable predual of the Banach $(X, \Vert.\Vert)$ and $x_{0}\in X$, then with $\sigma_{Y}=\sigma(X, Y)$:\\
\centerline{$\;\displaystyle\partial^{F}f(x_{0})\cap Y= \partial^{\sigma_{Y}}f(x_{0})\cap Y\;.$}
\end{lem}
Moreover, if $1< p\leq\infty$, when $E$ is reflexive, due to the Proposition \ref{prop100}, Lemma \ref{propw11a}, the results of the sections 9 and 10, Corollary \ref{coro713}, Corollary \ref{coro73},  Proposition \ref{prop 5.10},   Theorem \ref{corow11}, Theorem \ref{thm4.8} with its corollaries can be rephrased in terms of the weak dlcp (weak star if $p=\infty$) and of the weak (weak star if $p=\infty$) Hadamard subdifferentiability. Therefore the study of the weak (weak star) Dini Hadamard subdifferentiability is reduced to the study of Fr\'echet subdifferentiability in the cases  $1< p\leq\infty$. As a consequence, only the case $p=1$ may be of interest. In the sequel, $E$ is a separable reflexive space and we consider a $\mathbb{T}\otimes \mathcal{B}(E)$-measurable extended real valued integrand $f$ and its differential quotient $[f]$. Corollary \ref{coro73} gives a complete characterization of weak Hadamard subdifferentiability when the measure is atomless. The following statement is the analog of both results Theorem \ref{corow11} and Corollary \ref{coroww11} and gives, with Theorem \ref{thmw100}, practical sufficient conditions  for the weak Hadamard subdifferentiability.
\begin{thm} \label{thm6.17} Let $p=1$. If the integral functional $I_f$ is weakly Hadamard subdifferentiable at $x_0$ on $L_{1}(\Omega,E)$ then the integrand $f$ has the $\sigma$- dlcp at $x_0$. If the integrand $f$ has the $\sigma$-dlcp at $x_0$, then $\displaystyle L_{\infty}(\partial^{\sigma}f(x_{0}))\subseteq \partial^{\sigma}I_{f}(x_{0})$.
\end{thm}
Proof of Theorem \ref{thm6.17}.  Suppose that $x^*$ is a weak Hadamard subderivative of $I_f$ at $x_0$. Since the weak Hadamard bornology is solid, using Theorem \ref{prop 020} we deduce that the integrand $f-\langle x^{*}, .\rangle$ has the weak Hadamard dlcp; but due to the Dunford-Pettis criterion, for every relatively weakly compact sequence $(x_{n})_n$ the sequence $(\langle x^{*},x_{n}\rangle)_n$ is uniformly integrable, therefore $f$ has the weak Hadamard dlcp. Conversely, let ${(r_{n})}_n$ be a sequence of positive real numbers converging to the origin. Define $f_{n}(\omega,e)=[f]_{\omega}(x_{0}(\omega), e, r_{n})$. For each $x^{*}\in L_{\infty}(\partial^{\sigma}f(x_{0}))$ we have  for almost every $\omega\in \Omega$ the inequalities:
$$\langle x^{*}(\omega),\;.\rangle\leq f^{\sigma}_{\omega}(x_{0}(\omega),\; .)^{**}\leq (seq\;\sigma- li_{e}{f_{n}}_{\omega})^{**}\;.$$ 
Applying  Corollary \ref{thm6.1} to the sequence $\displaystyle (I_{f_{n}})_n$ we obtain for every  $x\in L_{1}(\Omega,E)$,
$$ I_{\langle x^{*},\;x\rangle}\leq I_{f^{\sigma}(x_{0},\; .)^{**}}(x)\leq I_{(seq\;\sigma- li_{e}f_{n})^{**}}(x)\leq\displaystyle seq\;\sigma-li_{e} I_{f_{n}}(x)\;.$$
This proves that for each $x\in L_{1}(\Omega,E)$, $\displaystyle \int_{\Omega}\langle x^{*}, x\rangle d\mu \leq I^{\sigma}_{f}(x_{0};x)$ or equivalently $x^{*}\in \partial^{\sigma}I_{f}(x_{0})$. $\;\Box$\\\\ 
\begin{thm}\label{thmw100} Let $p=1$, $x_{0}\in L_{1}(\Omega, E)$, and consider the following assertions\\
$(a)$ The integrand $f$ has the $\sigma$-dlcp at $x_0$,\\
$(b)$ There exists a positive constant $c$ such that the integrand $f$ satisfies out of a negligible set, for every $e\in E$,\\
\centerline{ $f(x_{0}+e)\geq f(x_{0})-c\Vert e\Vert $.}\\ 
Then $(b)\Rightarrow(a)$. If the measure $\mu$ is atomless these assertions are equivalent.
\end{thm}
Proof of Theorem \ref{thmw100}. 
By the Dunford-Pettis criterion, every sequentially weakly compact sequence ${(x_{n})}_n$ is uniformly integrable. Therefore when $(b)$ holds, for every sequence ${(r_{n})}_n$ of positive real numbers converging to the origin, the sequence ${([f]^{-}(x_{0}, r_{n}, x_{n}))}_n$ is uniformly integrable. This proves that  $(b)\Rightarrow(a)$. If the measure is atomless and $(a)$ is true, then the integrand $f$ has the (strong) Hadamard-dlcp at $x_0$,  thus \cite{49} Theorem 4.1 and \cite{49} Corollary 4.2 or \cite{216} Proposition 9, show that  $(a)\Rightarrow (b)$.$\;\Box$ \\\\
The author thanks Professor Jean Paul Penot for his help. And moreover for his suggestions and remarks about this article. Some consequences of these results in the study of the Clarke and limiting subdifferentials are listed in \cite{707}.

\end{document}